\documentclass[letterpaper,11pt]{article}
\usepackage{amsmath}
\usepackage{amsthm,amsfonts,amssymb}
\allowdisplaybreaks[1]

\usepackage{graphicx}
\usepackage{indentfirst}
\usepackage{mathrsfs}
\usepackage{amscd}
\usepackage{amsfonts}
\usepackage{epsfig}
\usepackage{color}
\usepackage{graphics}
\usepackage{dsfont}
\usepackage{bbm}
\usepackage{fullpage}
\usepackage[numbers]{natbib}
\usepackage{authblk}
\usepackage{setspace}
\usepackage{algpseudocode}
\usepackage{algorithm}
\floatstyle{plaintop}
\restylefloat{algorithm}
\usepackage{caption}
\usepackage{rotating}
\usepackage{multirow}
\usepackage{tabularx}
\usepackage{lscape}
\usepackage{array}
\usepackage{float}
\usepackage{url}
\numberwithin{equation}{section}

\newtheorem{theorem}{\bf Theorem}[section]
\newtheorem{proposition}{\bf Proposition}[section]
\newtheorem{lemma}{\bf Lemma}[section]

\theoremstyle{assumption}
\newtheorem{assumption}{\bf Assumption}[section]
\theoremstyle{definition}

\newtheorem{remark}{\bf Remark}[section]

\newcommand{\be}{\begin{equation}}
\newcommand{\ee}{\end{equation}}
\newcommand{\bes}{\begin{equation*}}
\newcommand{\ees}{\end{equation*}}
\newcommand{\ba}{\begin{aligned}}
\newcommand{\ea}{\end{aligned}}
\newcommand{\bi}{\begin{itemize}}
\newcommand{\ei}{\end{itemize}}

\newcommand{\hatx}{\hat{X}}
\newcommand{\barx}{\bar{X}}

\newcommand{\hatn}{\hat{N}}
\newcommand{\hatw}{\hat{W}}

\newcommand{\E}{\mathbb{E}}
\newcommand{\RR}{\mathbb{R}}
\newcommand{\ink}{$\blacksquare$}

\newcommand{\go}{\rightarrow}
\newcommand{\Go}{\Rightarrow}
\newcommand{\barn}{\bar{N}}

\newcommand{\map}{\mathcal{M}}
\newcommand{\textt}{\ \ \mbox}
\newcommand{\clf}{\mathcal{F}}
\newcommand{\NN}{\mathbb{N}}
\newcommand{\PP}{\mathbb{P}}
\newcommand{\var}{\mbox{Var}}
\newcommand{\tilz}{\tilde{Z}}
\newcommand{\hatmbx}{\hat{\mathbb{X}}}

\newcommand{\hata}{\hat{A}}

\title{\bf Diffusion Models for Double-ended Queues with Renewal Arrival Processes}
\author[1]{Xin Liu\thanks{Email: xliu9@clemson.edu.}}
\author[2]{Qi Gong\thanks{Email: qgong@email.unc.edu.}}
\author[2]{Vidyadhar G. Kulkarni\thanks{Email: vkulkarn@email.unc.edu.}}
\affil[1]{Department of Mathematical Sciences, Clemson University, Clemson, SC 29634.}
\affil[2]{Department of Statistics and Operations Research, University of North Carolina,
            \authorcr Chapel Hill, NC 27599.}

\begin{document}
\maketitle

\begin{abstract}
We study a double-ended queue where buyers and sellers
arrive to conduct trades. When there is a pair of buyer and seller in the
system, they immediately transact a trade and leave. Thus there
cannot be non-zero number of buyers and sellers simultaneously in
the system. We assume that sellers and buyers arrive at the system
according to  independent renewal processes, and they would leave the system after independent exponential patience times. We establish fluid and diffusion approximations for the queue length process under a suitable asymptotic regime. The fluid limit is the solution of an ordinary differential equation, and the diffusion limit is a time-inhomogeneous asymmetric Ornstein-Uhlenbeck process (O-U process). A heavy traffic analysis is also developed, and the diffusion limit in the stronger heavy traffic regime is a time-homogeneous asymmetric O-U process. The limiting distributions of both diffusion limits are obtained. We also show the interchange of the heavy traffic and steady state limits. \\

\end{abstract}

\section{Introduction}

Consider a simple trading market where sellers and buyers arrive according to independent renewal processes. When a seller is matched with a buyer, a trade occurs and they both leave the system. The trading follows first-come-first-served principle. If an arriving seller (buyer) cannot be matched with a buyer (seller), he/she will stay in a queue and wait for the upcoming buyers (sellers), and so there cannot be non-zero number of buyers and sellers simultaneously in the system. We further assume that traders (sellers and buyers) are impatient, that is, if they do not see a matching trader within a trader-specific random time (called the trader's patience time) they leave without completing the trade. Such system forms a double-ended queueing system, which is schematically shown in Figure \ref{fig1}. It is assumed that the arrival processes for sellers and buyers are independent renewal processes, and the patience times are independently exponentially distributed. A direct study of such system becomes challenging. In this work, we establish fluid and diffusion approximations for such double-ended queue in an appropriate asymptotic regime. The fluid limit is the solution of an ordinary differential equation, and the diffusion limit is a time-inhomogeneous asymmetric O-U process. A heavy traffic diffusion approximation is also studied, and the diffusion limit is a time-homogeneous asymmetric O-U process. We also show the validity of heavy traffic steady state approximation, i.e., the interchange of the heavy traffic and steady state limits. 

Double-ended queues arise in many applications, such as taxi-service system, buyers and sellers in a common market, assembly systems, organ transplant systems, to name a few. The first work on {\it double-ended queue} is by Kashyap \cite{kashyap1966double} for a taxi service example. Kashyap considers the taxi queueing system as a double-ended queue with limited waiting space. Under the assumptions that arrival processes of taxies and passengers are Poisson processes, he derives the analytical results about the steady state distribution of the system state. Conolly et al. \cite{conolly2002double} study the effect of impatience behavior primarily in the context of double-ended queues under the assumption of Poisson arrivals and exponential patience times. Researchers also find many other practical applications of the double-ended queues, such as networks with synchronization nodes (Prabhakar et al. \cite{prabhakar2000synchronization}), and perishable inventory system (Perry et al. \cite{perry1999perishable}). When renewal arrivals are considered, the explicit form of the limiting distribution becomes intractable. Degirmenci \cite{talayasymptotic} studies the asymptotic behavior of the limiting distribution of the double-ended queue using algebraic approximation methods. Several researchers study the double-ended queue using simulation methods, see Zenios \cite{zenios1999modeling} and Kim et al. \cite{kim2010simulation}. In this work, we develop rigorous diffusion approximations for double-ended queues under appropriate asymptotic regime.

There is a rich literature on diffusion approximations for (one-sided) queueing systems with abandonment in heavy traffic. Two heavy traffic regimes -- conventional heavy traffic regime and Halfin-Whitt regime -- have been extensively studied. Loosely speaking, in both regimes, the queueing system is roughly balanced. In conventional heavy traffic regime, one considers queueing systems with fixed number of servers, while in Halfin-Whitt regime, the number of servers approaches to infinity. Ward and Glynn \cite{WardGlynn03} study the $M/M/1 + M$ model ($+M$ denotes independent exponential patience times), and the result is extended to the $G/GI/1+GI$ model ($+GI$ denotes generally distributed independent patience times) in Ward and Glynn \cite{WardGlynn05}. Later on, Reed and Ward \cite{Reed08} develop a more stable hazard rate scaling for patience time distribution for the $G/GI/1 + GI$ model. Such scaling is extended more generally for single server queue in Lee and Weerasinghe \cite{Lee11}. In the Halfin-Whitt regime, the $M/M/n+M$ model is considered in Garnett et al. \cite{Garnet02}, Zeltyn and Mandelbaum \cite{Man05} study the $M/M/n+G$ model, and Dai et al. \cite{Dai10} work on the $G/PH/n + GI$ model. Recently, Mandelbaum and Momcilovic \cite{Man12} and Dai and He \cite{DaiHe} both develop diffusion approximations for $G/GI/n+GI$ model. The hazard rate scaling has been applied to study the $G/M/n+GI$ model by Reed and Tezcan \cite{ReedTezcan12}. In this work, we will focus on exponential patience times. In a forthcoming paper Liu \cite{liu}, we study generally distributed independent patience times with hazard rate scaling, and the (conventional) heavy traffic diffusion limit is expected to be an asymmetric O-U process with drift given by an appropriate hazard rate scaling limit of the patience time distributions.

Asymmetric O-U processes developed in our work have piecewise-linear state dependent drift and possibly time dependent diffusion coefficient. The drift function has linear pieces on $[0,\infty)$ and $(-\infty, 0)$, which are corresponding to different reneging rates of sellers and buyers (see \eqref{limit} and \eqref{diff-ht-eqn}). Because of such drift structure, the asymmetric O-U processes admit unique limiting distributions (see Theorem \ref{limitphi}). Such process is a special case of the so-called piecewise O-U processes, which have arisen as diffusion approximations for queueing systems, e.g. Garnett et al. \cite{Garnet02} and Dai et al. \cite{Dai10}. In Browne and Whitt \cite{BrowneWhitt95}, a closed form of the stationary distributions of one-dimensional time-homogeneous piecewise-linear diffusion processes is established. Recently, Dieker and Gao \cite{DiekerGao13} study the positive recurrence property of multidimensional time-homogeneous piecewise O-U processes, which arises as diffusion approximation in Dai et al. \cite{Dai10}. 

Furthermore, we study the validity of the heavy traffic steady state approximation (see Theorem \ref{limit-change}). More precisely, we show that the pre-limit queue length process admits a limiting distribution, and in heavy traffic such limiting distribution converges to that of the limit diffusion process. A pioneer work on the study of such interchange of limit operations for generalized Jackson networks is by Gamarnik and Zeevi \cite{GamarnikZeevi06}, the essiential idea of which is to construct an appropriate Lyapunov function. Budhiraja and Lee \cite{BudhirajaLee09}, on the other hand, study the interchange of limits by establishing a uniform bound for the growth of the moments of the diffusion scaled pre-limit processes. Recently, Dai et al. \cite{Dai13} study the interchange of heavy traffic and steady state limits for the $G/PH/n + M$ model, following the approach in Gamarnik and Zeevi \cite{GamarnikZeevi06}, and establish an interchange limit theorem under some sufficient conditions. In this work, we will follow the approach in Chapter 4 of Bramson \cite{bramson08} and Budhiraja and Lee \cite{BudhirajaLee09} to establish the positive recurrence of the pre-limit processes and the interchange of heavy traffic and steady state limits under very mild assumptions \eqref{positivity} and \eqref{uniform-int}.

The main contributions of this work are (1) we study a two-sided queue, which is different in structure from the one-sided queue; (2) we not only study the (conventional) heavy traffic diffusion approximation, but also study a diffusion approximation for the centered and scaled queue length process without assuming heavy traffic condition; (3) the limiting behaviors of diffusion models are developed; (4) we establish the interchange of the heavy traffic and steady state limits for the queue length process; (5) we conduct numerical experiments to study the goodness of the fluid and diffusion models. 

The rest of the paper is organized as follows. In the next section we present the model of the double-ended queue with renewal arrivals and exponential patience times, and introduce the relevant notation. In Section \ref{sec:PA}, we collect the results about the special case when the arrival processes are Poisson processes. Some of these results are known, while some are new. We use these results to construct a Poisson approximation model in Section \ref{sec:NE}. In Section \ref{sec:DA} we study the fluid and diffusion approximations for the queue length process, and the main results are presented in Theorems \ref{fluid-d} -- \ref{limit-change}. More precisely, under suitable conditions (Assumptions \ref{main-d}), the fluid limit is provided in Theorem \ref{fluid-d}, the first diffusion approximation result appears in Theorem \ref{diffusion-d}, and the heavy traffic diffusion limit is obtained in Theorem \ref{diff-ht}. We further provide the exact solution of the fluid equation in Lemma \ref{fluid-s}, and study the moments and limiting distributions of the diffusion limits in Theorem \ref{limitphi}. Finally, we establish the interchange limit theorem in Theorem \ref{limit-change}. In Section \ref{sec:NE}, we study several numerical examples, and compare goodness of three approximations: the Poisson approximation, and the two diffusion approximations. We make comments on extensions of this model in Section \ref{sec:C}. All the proofs are provided in Appendix A, and the numerical results can be found in Appendix B. 

We use the following notation. Denote by $\mathbb{R},$ $\mathbb{R}_+,$ $\mathbb{Z}$, and $\mathbb{N}$ the sets of real numbers, nonnegative real numbers, integers, and positive integers, respectively. For a real number $a$, define $a^+ = \max\{a, 0\}$ and $a^- = \max\{0,-a\}.$ Similarly, for a real function $f$ defined on $[0,\infty)$, define $f^+(t) = \max\{0, f(t)\}$ and $f^-(t) = \max\{0, -f(t)\}, \ t\ge 0.$ Denote by $\mathcal{D}([0,\infty): \RR)$ the space of right continuous functions with left limits defined from $[0,\infty)$ to $\RR$ with the usual Skorohod topology. For $x\in \mathcal{D}([0,\infty): \RR)$, let $\|x\|_t = \sup_{s\in[0, t]} |x(s)|, \ t\ge0.$
 A mapping $F: \mathcal{D}([0,\infty): \RR) \to \mathcal{D}([0,\infty): \RR)$ is called Lipschitz continuous if for any $t\in [0,\infty),$ there exists $\kappa\in (0,\infty)$ (may depending on $t$) such that for $x_1, x_2 \in \mathcal{D}([0:\infty), \RR)$,
\bes
\|F(x_1) - F(x_2)\|_t \le \kappa \|x_1 - x_2\|_t.
\ees
The normal distribution with mean $\mu$ and variance $\sigma^2$ is denoted by $N(\mu,\sigma^2),$ and its density and distribution functions are denoted by $\phi(\cdot; \mu, \sigma^2)$ and $\Phi(\cdot; \mu, \sigma^2)$, respectively. For a Markov process $\{X(t): t\ge 0\}$ with stationary distribution $\pi$, denote by $X(\infty)$ a random variable with distribution $\pi$. Finally, we will denote generic positive constants by $c_1, c_2, \ldots$. Their values may change from one proof to another.

\section{Model formulation} \label{sec:MF}

\begin{figure}[h]
\centering
  \includegraphics[width=12cm]{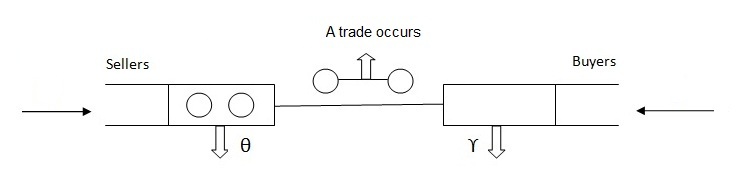}
  \caption{A double-ended queue.} \label{fig1}
\end{figure}

Consider a double-ended queue as in Figure \ref{fig1}, where the sellers and buyers arrive
according to independent renewal processes. A trade occurs when there is a pair of seller and buyer in the system, and the pair leaves the system instantaneously. The trading follows first-come-first-served principle. So there cannot be nonzero sellers and buyers in the system simultaneously. We assume the customers are impatient and the patience times are independently exponentially distributed. 
Let $X(t)$ be the length of the double-ended queue at time $t$. Then $X(t)$ takes values in $\mathbb{Z}$. 
If $X(t)>0$, there are $X(t)$ sellers waiting
in the system, and if $X(t)<0$, there are $-X(t)$ buyers waiting in
the system. Let $(\Omega, \mathcal{F}, \PP, \{\mathcal{F}_t\}_{t\geq 0})$ be a filtered
probability space satisfying the usual conditions. All the random variables and stochastic processes in Sections \ref{sec:MF} and \ref{sec:PA} are assumed to be defined on this space. We assume the inter-arrival times of sellers and buyers are independent sequences of i.i.d. random
variables $\{U_k: k\in \mathbb{N}\}$ and $\{V_k: k\in
\mathbb{N}\}$, respectively. The random variable $U_1$ has mean $1/\alpha$ and standard deviation $\sigma$, and $V_1$ has mean $1/\beta$ and standard deviation $\varsigma$.
Define
\begin{equation*}
\begin{aligned}
N_s(t) &= \max\left\{k: \sum_{i=1}^k U_i\leq t\right\},\\
N_b(t) &= \max\left\{k: \sum_{i=1}^k V_i\leq t\right\}.
\end{aligned}
\end{equation*}
The renewal processes  $N_s$ and  $N_b$ can be interpreted as the arrival processes for sellers and buyers, respectively. The patience times of sellers and buyers are independent sequences of i.i.d. exponential random variables with rate $\theta$ and $\gamma$, respectively. Let $N_{sr}$ and $N_{br}$ be two independent unit-rate Poisson processes, which are independent of $N_s$ and $N_b$.  Then we have the following evolution equation for $\{X(t): t\ge 0\}$. For $t\geq 0$,
\begin{equation}\label{eq:1}
X(t)=X(0)+N_s(t)-N_b(t)-N_{sr}\left(\theta \int_0^t
X^{+}(s)ds\right)+N_{br}\left(\gamma\int_0^t
X^{-}(s)ds\right),
\end{equation}
where $X(0)$ denotes the initial number of sellers or buyers in the system, which is assumed to be independent of $N_s, N_b, N_{sr}$ and $N_{br}.$

\section{A special case: Poisson arrivals} \label{sec:PA}

When the arrival processes are Poisson processes, it is easy to see that
$\{X(t): t \ge 0\}$ is a birth and death process on $\mathbb{Z}$
with birth parameters $\lambda_i = \alpha + i^-\gamma$ and death
parameters $\mu_i = \beta + i^+\theta$ for $i\in\mathbb{Z}$. Using the standard theory
(see Kulkarni \cite{kulkarni1996modeling}), we see that this
birth and death process is:
\begin{itemize}
\item positive recurrent, if $\theta > 0$ and $\gamma > 0$;
\item null recurrent, if $\theta = \gamma = 0$ and $\alpha = \beta$;
\item transient, if $\theta = \gamma = 0$ and $\alpha \neq \beta$.
\end{itemize}
In the analysis, we assume $\theta > 0$ and $\gamma > 0$.
Therefore, this continuous time Markov chain (CTMC) $X$ has a unique limiting distribution (which is also the unique stationary distribution). Let $\{\pi_i: \ i\in \mathbb{Z}\}$ denote this limiting distribution. Using the standard theory of balance equations, the limiting distribution is given by the following:
\begin{eqnarray}
{\pi _i} &=& \frac{{{\alpha ^i}}}{{\prod\limits_{j = 1}^i
 {(\beta+j\theta )} }}{\pi _0},\;\;i\in\NN, \label{limitdis1} \\
{\pi _{ - i}} &=& \frac{{{\beta ^i}}}
{{\prod\limits_{j = 1}^i {(\alpha+j\gamma) } }}{\pi
_0}, \;\;i\in\NN, \label{limitdis2}\\
 {\pi _0} &=& \left( {1 +
\sum\limits_{i = 1}^\infty {\frac{{{\alpha ^i}}}{{\prod\limits_{j =
1}^i {(\beta+j\theta) } }} + \sum\limits_{i = 1}^\infty
{\frac{{{\beta ^i}}}{{\prod\limits_{j = 1}^i {(\alpha+j\gamma) } }}}
} } \right)^{-1}. \label{limitdis3}
\end{eqnarray}

The next lemma studies some ergodicity properties of $X.$ Denote by $\E_x$ the expectation conditioning on the process $X$ starting from $x\in\mathbb{Z}$.  

\begin{lemma}\label{limit-moments} For $x\in \mathbb{Z}$ and $s\ge 0,$
\be\label{mgf}
\lim_{t\to\infty } \E_x\left[e^{s|X(t)|}\right] = \sum_{i=-\infty}^\infty e^{s|i|}\pi_i. 
\ee
In particular, for measurable function $f: \mathbb{Z} \to \RR$ such that for some $s\ge 0,$ $|f(x)| \le e^{s|x|}, \ x\in \mathbb{Z}$, we have 
\be\label{mfs}
\lim_{t\to\infty } \E_x\left[f(X(t))\right] = \sum_{i=-\infty}^\infty f(i)\pi_i. 
\ee
\end{lemma}

Using Lemma \ref{limit-moments}, various limiting performance of $X$ can be studied by using the limiting distribution $\pi.$  In the following, we simply focus on the first two moments of $X(t)$. From Lemma \ref{limit-moments}, the first two limiting (or steady state) moments can be expressed in terms of $\pi$. We next simplify these moments using Gamma functions and incomplete Gamma functions. Define for $t\ge 0,$
\begin{equation*}
\ba
m(t) &= \E(X(t)), \ m_+(t) = \E(X^+(t)),\ m_-(t) = \E(X^-(t)),\\
s(t) &= \E(X(t))^2, \ s_+(t) = \E(X^+(t))^2,\ s_-(t) = \E(X^-(t))^2.
\ea
\end{equation*}
Clearly, $m(t)=m_+(t)-m_-(t)$ and $s(t)=s_+(t)+s_-(t)$. We also note that in general $m_+(t) \neq m^+(t), \ m_-(t) \neq m^-(t),$ and so $|m(t)|\neq m_+(t) + m_-(t)$. For $t>0,$ the Gamma function $\Gamma(t) = \int_0^\infty x^{t-1}e^{-x} dx$, and for $t>0, y>0,$ the incomplete Gamma function $\Gamma(t,y)=\int_0^y x^{t-1}e^{-x}dx$. 

\begin{lemma}\label{simp}
Let 
\[ p_1 = \frac{\beta}{\theta} e^{\alpha/\theta}  \left(\frac{\alpha}{\theta}\right)^{-\beta/\theta} \Gamma(\beta/\theta, \alpha/\theta)- 1, \]
\[ p_2 = \frac{\alpha}{\gamma}  e^{\beta/\gamma} \left(\frac{\beta}{\gamma}\right)^{-\alpha/\gamma} \Gamma(\alpha/\gamma, \beta/\gamma)- 1.\]
Then
\[\pi_0  = (1+p_1+p_2)^{-1},\]
and
\bes\ba
\lim_{t\to\infty}m_+(t) &= \left[\frac{\alpha-\beta}{\theta} p_1 +\frac{\alpha}{\theta}\right]\pi_0, \\
\lim_{t\to\infty}m_-(t) &=\left[\frac{\beta-\alpha}{\gamma} p_2 +\frac{\beta}{\gamma}\right]\pi_0, \\
\lim_{t\to\infty}s_+(t) & = \frac{\alpha-\beta}{\theta} \lim_{t\to\infty} m_+(t) + \frac{\alpha}{\theta}(p_1 +1)\pi_0,\\
\lim_{t\to\infty}s_-(t)&= \frac{\beta-\alpha}{\gamma} \lim_{t\to\infty} m_-(t) + \frac{\beta}{\gamma}(p_2 +1)\pi_0. 
\ea\ees
\end{lemma}

Lemma \ref{simp} will be used in Section \ref{sec:NE} for Poisson approximation. We next introduce $m(t)$ and $s(t)$ as the solutions of two ordinary differential equations. 
 
\begin{lemma} \label{th:mom}
Assume that $X(0)$ has finite first two moments. Then the moment functions $m(t)$ and $s(t)$ satisfy  the following differential equations. For $t\ge0,$
\begin{equation}\label{m1}
\frac{{dm(t)}}{{dt}} = (\alpha  - \beta ) -  \theta {m_ + }(t) +
\gamma {m_ - }(t),
\end{equation}
and
\begin{equation}\label{m2}
\frac{{ds(t)}}{{dt}} =  - 2\theta {s_+ }(t) - 2 \gamma {s_ - }(t) +
(2\alpha  - 2\beta  + \theta ){m_ + }(t) + ( - 2\alpha  + 2\beta  +
\gamma ){m_ - }(t) + \alpha  + \beta,
\end{equation}
with initial conditions $m(0) = \E(X(0))$ and $s(0)=\E(X(0))^2.$
\end{lemma}

From Lemma \ref{th:mom}, when $\theta = \gamma,$ \eqref{m1} and \eqref{m2} are simplified to be
\begin{equation}\label{m11-ode}
\frac{{dm(t)}}{{dt}} = (\alpha  - \beta ) -  \theta {m }(t),
\end{equation}
and
\begin{equation}\label{m22-ode}
\frac{{ds(t)}}{{dt}} =  - 2\theta {s }(t) + (2\alpha  - 2\beta) m(t)  + \theta (m_+(t) + {m_ - }(t)) + \alpha  + \beta.
\end{equation}
Solving \eqref{m11-ode}, we have 
\be\label{m1-ode-soln}
m(t) = \left(m(0)-\frac{\alpha-\beta}{\theta}\right) e^{-\theta t}+\frac{\alpha-\beta}{\theta}, \ t\ge 0, \ \mbox{and} \  \lim_{t\to\infty} m(t) =  \frac{\alpha - \beta}{\theta}.
\ee
However, observing that $m_+(t) + m_-(t) \neq |m(t)|$, \eqref{m22-ode} cannot be solved directly. Consider the following ODE by replacing $m_+(t) + m_-(t)$ with $|m(t)|$,
\be\label{new-s}
\frac{d\tilde s(t)}{dt}  =  - 2\theta \tilde s(t) + (2\alpha - 2\beta )m(t) + \theta | {m}(t)| + \alpha  +
\beta, \ \mbox{and $\tilde s(0) = s(0)$.}
\ee
Noting that $m_+(t) + m_-(t) \ge |m(t)|$ and $s(0)=\tilde s(0)$, we have that
\be\label{lower-bound}
s(t) \ge \tilde s(t), \ t\ge 0.
\ee
Using \eqref{m1-ode-soln} to solve the ODE in \eqref{new-s}, we have that
\be\label{s11}
\lim_{t\to\infty} s(t) \ge \lim_{t\to\infty} \tilde s(t) = \left(\frac{\alpha-\beta}{\theta}\right)^2 + \frac{\max\{\alpha, \beta\}}{\theta},
\ee
which provides a lower bound for $\lim_{t\to\infty} s(t).$ Nevertheless, we can always use Lemma \ref{simp} to study $\lim_{t\to\infty} s(t)$ as follows:
\be\label{MC-sec-moment}\ba
\lim_{t\to\infty} s(t)& = \lim_{t\to\infty}s_+(t) + \lim_{t\to\infty}s_-(t)\\
& = \frac{\alpha-\beta}{\theta} \lim_{t\to\infty} m(t) + \frac{\alpha}{\theta}(p_1+1)\pi_0 + \frac{\beta}{\gamma}(p_2+1)\pi_0 \\ 
& = \frac{\alpha+\beta}{2\theta} + \frac{(\alpha-\beta)^2}{\theta^2} + \frac{ \lim_{t\to\infty}m_+(t) +  \lim_{t\to\infty} m_-(t)}{2}.
\ea\ee
In Remark \ref{compare-remark}, we construct proper estimates for $m(t)$ and $s(t)$ using the fluid and diffusion approximations established in Section \ref{sec:DA}.  It is shown that, under Assumption \ref{main-d} (additional heavy traffic condition needed for heavy traffic diffusion model), the fluid and diffusion models provide simple and good approximations for $\lim_{t\to\infty}m(t)$ and $\lim_{t\to\infty}s(t).$

\section{Fluid and diffusion approximations} \label{sec:DA}

In this section we study double-ended queues with renewal arrivals and exponential patience times. In this setting, $\{X(t): t\ge 0\}$ is no longer a Markov process. We will focus on establishing fluid and diffusion approximations for $\{X(t): t\ge 0\}$ under appropriate conditions (see Assumption \ref{main-d}). To describe the asymptotic region where such approximations are valid, we consider a sequence of double-ended queues indexed by $n\in\mathbb{N}$. For the $n$-th system, all the notation introduced in Section \ref{sec:MF} is carried forward except that we append a superscript $n$ to all quantities to indicate the dependence of parameters, random variables, and stochastic processes on $n$. In particular, on the space $(\Omega^n, \mathcal{F}^n, P^n, \{\mathcal{F}^n_t\}_{t\geq 0})$, $\{U^n_k: k\in\NN\}$ and $\{V^n_k: k\in\NN\}$ are the sequences of interarrival times, $N_s^n$ and $N_b^n$ are the arrival processes,  $\theta^n$ and $\gamma^n$ are the reneging rates, and $N_{sr}^n$ and  $N_{br}^n$ are the unit-rate Poisson processes used to formulate the abandonment processes. Also $1/\alpha^n, \sigma^n$ and $1/\beta^n, \varsigma^n$ are the means and standard deviations of the inter-arrival times of sellers and buyers, respectively. The expectation operator with respect to $\PP^n$ will be denoted by $\E^n$, but frequently we will suppress $n$ from the notation.
We further assume the following strict positivity and uniform integrability on $\{U^n_1: n\in\NN\}$ and $\{V^n_1: n\in\NN\}$.
\be\label{positivity}
\PP^n(U^n_1>0) = \PP^n(V^n_1>0) =1\ \mbox{for all $n\in\NN$}.
\ee
\be\label{uniform-int}
\mbox{$\{(U^n_1)^2: n\in\NN\}$ and $\{(V^n_1)^2: n\in\NN\}$ are uniformly integrable.}
\ee
Finally, the queue length process $X^n$ can be described as follows:
For $t\geq 0$,
\begin{equation}\label{queue}
X^n(t)=X^n(0)+N_s^n(t)-N_b^n(t)-N_{sr}^n\left(\theta^n \int_0^t
X^{n,+}(s)ds\right)+N_{br}^n\left(\gamma^n\int_0^t
X^{n,-}(s)ds\right).
\end{equation}
The following assumption describes the asymptotic regime of the parameters. We will assume Assumption \ref{main-d} holds for the entire section. 
\begin{assumption}\label{main-d}\hfill
\begin{enumerate}
\item[\rm (i)] There exist $\alpha, \beta, \sigma, \varsigma\in (0,\infty)$ such that
 $$\alpha^n\rightarrow\alpha, \ \beta^n\rightarrow\beta, \
  \sigma^n\rightarrow\sigma, \ \varsigma^n\rightarrow\varsigma.$$
\item[\rm (ii)] For $\theta, \gamma\in (0,\infty),$ we have that
 $$n\theta^n\rightarrow\theta, \ n\gamma^n\rightarrow\gamma.$$
\end{enumerate}
\end{assumption}

\begin{remark} \hfill
\bi
\item Assumption \ref{main-d} says the means and variances of the arrival processes are $\mathcal{O}(1).$ However, the reneging rates are $\mathcal{O}(n^{-1}).$ So comparing with the arrival rates, the reneging rates are very small. However, such reneging is nonnegligible in both fluid and diffusion approximations. Similar scalings are widely considered in the diffusion approximations for queueing systems with reneging, e.g. \cite{WardGlynn03, Garnet02}. 
\item The diffusion analysis of queueing systems often assumes heavy traffic conditions (under which the queueing system is roughly balanced) and establish an approximation for the diffusion scaled queue length process (or other interesting processes). In this work, we first study the diffusion scaled queue length process centered at the fluid limit, and develop a suitable scaling theorem (see Theorem \ref{diffusion-d}) only under Assumption \ref{main-d}. We then establish the heavy traffic analysis in Theorem \ref{diff-ht} with an additional heavy traffic condition \eqref{htc}.   
\ei
\end{remark}

\subsection{Fluid approximation}\label{fluid-app}

We begin by defining the fluid scaled processes. Loosely speaking, we scale up the time by the factor $n$ and scale down the space by the same factor $n$.  More precisely, for $t\geq 0$, define
\bes\label{fluid-scale}
 \barx^n(t) = \frac{X^n(nt)}{n}, \; \barn^n_s(t) = \frac{N_s^n(nt)}{n}, \; \barn^n_b(t) = \frac{N_b^n(nt)}{n}, \ \barn^n_{sr}(t) = \frac{N_{sr}^n(nt)}{n}, \; \barn^n_{br}(t) = \frac{N_{br}^n(nt)}{n}.
\ees
Recall that for a stochastic process $\{Y(t),t\ge0\}$, $
{\left\| Y \right\|_t} = \mathop {\sup }\limits_{0 \le u \le t} \left| {Y(u)} \right|, \ t\in [0,\infty).$
We first obtain the limit of $\bar X^n$ as $n \rightarrow \infty$ in Theorem \ref{fluid-d}. The solution of the fluid limit equation is then given in Lemma \ref{fluid-s}. 
\begin{theorem}\label{fluid-d}
Assume that for some $x_0\in \RR$, $\E(|\bar X^n(0) - x_0|)\to 0$ as $n \to \infty$. Then we have that for $t\in[0,\infty),$
\be\label{con}
\E\left(\|\bar X^n -x \|_t\right) \to 0, \textt{as $n\to \infty$},
\ee
where $x$ is the solution of the following integral equation
\be\label{ode}
x(t) = x_0 + (\alpha-\beta)t - \theta \int_0^t x^+(s) ds + \gamma \int_0^t x^-(s) ds, \ t\ge 0.
\ee

\end{theorem}

\begin{lemma}\label{fluid-s}
Consider the integral equation in \eqref{ode}. 
\bi
\item[\rm (i)] If $\alpha \ge \beta$ and $x_0\ge 0$, then
\bes
x(t) = \left(x_0 - \frac{\alpha -\beta}{\theta}\right) e^{-\theta t} + \frac{\alpha - \beta}{\theta}, \ t\in[0, \infty).
\ees
\item[\rm (ii)] If $\alpha \ge \beta$ and $x_0 < 0$, then
\bes
x(t)  = \begin{cases}
\left(x_0 - \frac{\alpha -\beta}{\gamma}\right) e^{-\gamma t} + \frac{\alpha - \beta}{\gamma}, & t\in[0, t_1], \\
\frac{\alpha -\beta}{\theta}\left(1-e^{-\theta (t-t_1)}  \right), &  t\in[t_1,\infty),
 \end{cases}
\ees
where
\bes
t_1 = \gamma^{-1}\log\left(\frac{\alpha -\beta -\gamma x_0}{\alpha - \beta}\right)
\ees
is the first time for $x$ to reach $0.$ 

\item[\rm (iii)] If $\alpha < \beta$ and $x_0 \le 0$,  then
\bes
x(t)  = \left(x_0 -  \frac{\alpha -\beta}{\gamma} \right) e^{-\gamma t} + \frac{\alpha -\beta}{\gamma}, \ t\in[0,\infty).
\ees
\item[\rm (iv)] If $\alpha < \beta$ and $x_0 > 0$,  then
\bes
x(t)  = \begin{cases}
\left(x_0 - \frac{\alpha -\beta}{\theta}\right) e^{-\theta t} + \frac{\alpha - \beta}{\theta}, & t\in[0, t_2], \\
\frac{\alpha -\beta}{\gamma}\left(1-e^{-\gamma (t-t_2)}  \right), &  t\in[t_2,\infty),
 \end{cases}
\ees
where
\bes
t_2 = \theta^{-1}\log\left(\frac{\alpha -\beta -\theta x_0}{\alpha - \beta}\right)
\ees
is the first time for $x$ to reach $0.$
\ei

\end{lemma}

\begin{remark} From Lemma \ref{fluid-s}, after finite time ($t_1$ and $t_2$ in Lemma \ref{fluid-s}), the fluid limit $x$ attains zero. The larger the reneging rate and the deviation between the two arrival rates, the faster $x$ attains zero. After reaching zero, $x(t)$ has the same sign as $\alpha-\beta$, and eventually approaches a stable point. The existence of the positive reneging rates $\theta$ and $\gamma$ guarantees such stability of $x(t)$ as $t\to\infty.$ Indeed, we have
\be\label{limit-fluid}
\lim_{t\to\infty} x(t) = \begin{cases} \frac{\alpha - \beta}{\theta}, & \alpha \ge \beta, \\  \frac{\alpha - \beta}{\gamma}, & \alpha < \beta. \end{cases}
\ee

\end{remark}

\subsection{Diffusion approximations}

In this subsection, we study diffusion approximations and define the diffusion scaled processes. This time we scale up the time by the same factor $n$ and scale down the space by factor $\sqrt{n}$.  To be  precise, for $t\geq 0$, define
\be\label{diff-scale}\ba
&\hatx^n(t) = \frac{X^n(nt)}{\sqrt{n}}, \; \hatn^n_{s}(t) = \frac{N^n_{s}(nt)-n\alpha^nt}{\sqrt{n}}, \; \hatn^n_{b}(t) = \frac{N^n_{b}(nt)-n\beta^nt}{\sqrt{n}}, \\
& \hatn^n_{sr}(t) = \frac{N^n_{sr}(nt)-nt}{\sqrt{n}}, \; \hatn^n_{br}(t) = \frac{N^n_{br}(nt)-nt}{\sqrt{n}}.
\ea\ee

We now state our main results. Define a diffusion process $Z$ as follows. Recall the fluid limit $x$ in Section \ref{fluid-app}. For a given random variable $Z(0)$
with law $\nu$, and a standard Brownian motion $W$, let $Z$ be the
unique solution to the following stochastic integral equation
\be\label{limit} Z(t)= Z(0) + \int_0^t \sqrt{{\alpha ^3}{\sigma ^2} + {\beta ^3}{\varsigma ^2} + \theta x^+(u) + \gamma x^-(u) } dW(u)  - \theta\int_0^t   Z^+(u) du +\gamma\int_0^t  Z^-(u)du.  \ee 
The existence and uniqueness of $Z$ is guaranteed by the following lemma.

\begin{lemma}[Reed and Ward \citep{reed2004diffusion}] \label{lip-mapping}
Let $\phi: \mathcal{D}([0,\infty), \RR)\rightarrow \mathcal{D}([0,\infty), \RR)$ be Lipschitz continuous. Then for any given $w\in D([0,\infty), \RR)$, there exists a unique $y\in D([0,\infty), \RR)$ that satisfies the integral
\bes
y(t) = w(t) + \int_0^t \phi(y)(u)du,
\ees
and $y(0) = w(0)$. Moreover, define the mapping $\map^\phi: D([0,\infty), \RR)\rightarrow D([0,\infty), \RR)$ by $\map^\phi(w) = y$, and then $\map^\phi$ is Lipschitz continuous.
\end{lemma}

Recall the fluid limit $x$ in Theorem \ref{fluid-d}, and define 
\be\label{pre-fluid}
x^n(t) = x(0) + (\alpha^n - \beta^n)t - n\theta^n \int_0^t x^{n,+}(s)ds + n\gamma^n \int_0^t x^{n,-}(s)ds. 
\ee
It is clear that $x^n(t) \to x(t)$ for all $t\ge 0.$
Let 
\be\label{centered}
Z^n(t) = \sqrt{n} (\barx^n(t) - x^n(t)) = \hatx^n(t) - \sqrt{n}x^n(t), \ t\ge 0.
\ee
 
The following are the main results for the diffusion approximations.

\begin{theorem}\label{diffusion-d} Assume that $Z^n(0)$ converges weakly to  a
probability measure $\nu$. Then
$Z^n\Rightarrow Z$, where $Z$ is defined by \eqref{limit}.
\end{theorem}

\begin{theorem}\label{diff-ht}
Assume that $\E|\barx^n(0)| \to 0$, and there exists $c\in \RR$ such that 
\be\label{htc}
\sqrt{n}(\alpha^n-\beta^n) \to c, \ \mbox{ as $n\to\infty.$}
\ee
Then we have $x\equiv 0,$ and, if $\hatx^n(0)$ converges weakly to a probability measure $\mu$, then $\hatx^n \Go \hatx$, where
\be\label{diff-ht-eqn}
\hatx(t) = \hatx(0) + \sqrt{{\alpha ^3}{\sigma ^2} + {\beta ^3}{\varsigma ^2}} W(t) + c t  - \theta\int_0^t   \hatx^+(u) du +\gamma\int_0^t  \hatx^-(u)du, \ t\ge 0,
\ee
and here $W$ is a standard Brownian motion and $\hatx(0)$ has law $\mu$.
\end{theorem}

Condition \eqref{htc} is well known as the heavy traffic condition in queueing theory, and Theorem \ref{diff-ht} provides a heavy traffic diffusion analysis of double-ended queues. 

\begin{remark}\label{remark-diff}
Under Assumption \ref{main-d} and the heavy traffic condition \eqref{htc}, the fluid limit $x\equiv 0,$ and the diffusion limit in Theorem \ref{diffusion-d} is reduced to be 
\be\label{reduced-diff}
Z(t) = Z(0) + \sqrt{{\alpha ^3}{\sigma ^2} + {\beta ^3}{\varsigma ^2}} W(t)  - \theta\int_0^t   Z^+(u) du +\gamma\int_0^t  Z^-(u)du, \ t\ge 0.
\ee
We next observe that \eqref{pre-fluid} becomes
\bes
x^n(t) = (\alpha^n-\beta^n)t - n\theta^n \int_0^t x^{n,+}(s)ds + n\gamma^n \int_0^t x^{n,-}(s)ds, \ t\ge0.
\ees
Noting that $\sqrt{n}(\alpha^n -\beta^n) \to c$, we have  
\bes
\sqrt{n}x^n(t) = \sqrt{n}(\alpha^n-\beta^n)t - n\theta^n \int_0^t \sqrt{n}x^{n,+}(s)ds + n\gamma^n \int_0^t \sqrt{n}x^{n,-}(s)ds, \ t\ge0,
\ees
and $\sqrt{n}x^n(t) \to \hat{x}(t)$ for all $t\ge 0,$ where
\be\label{scaled-fluid}
\hat{x}(t) = ct - \theta \int_0^t \hat{x}^+(s)ds + \gamma \int_0^t \hat{x}^-(s)ds, \ t\ge0.
\ee
Recalling that $Z^n(t) = \hatx^n(t) - \sqrt{n}x^n(t), t\ge 0$, we note that under Assumption \ref{main-d} and the heavy traffic condition \eqref{htc},  
\be\label{relation-X-Z}
\hatx(t) = Z(t) + \hat{x}(t), \ t\ge 0,
\ee
where $\hatx$ is the heavy traffic diffusion limit, and $Z$ and $\hat{x}$ are defined in \eqref{reduced-diff} and \eqref{scaled-fluid}, respectively. 

\end{remark}

When $\theta =\gamma,$ the diffusion process $Z$ is a time-inhomogeneous O-U process, and the process $\hatx$ is a time-homogeneous O-U process. The following lemma collects some well known results for O-U processes (see e.g. \cite{Karatzas}). 

\begin{lemma}\label{o-u}
Assume $\theta =\gamma.$ Then
\bi
\item[\rm (i)] 
\bes
Z(t) = e^{-\theta t} Z(0) + \int_0^t e^{-\theta (t- u)}\sqrt{{\alpha ^3}{\sigma ^2} + {\beta ^3}{\varsigma ^2} + \theta x^+(u) + \gamma x^-(u) } dW(u).
\ees
\item[\rm (ii)] The unique limiting distribution of $Z$ is a normal distribution with mean $0$ and variance $({\alpha ^3}{\sigma ^2} + {\beta ^3}{\varsigma ^2} +|\alpha-\beta|)/2\theta$.
\item[\rm (iii)] 
\bes
\hatx(t) = e^{-\theta t} \hatx(0) + c(1-e^{-\theta t}) + \int_0^t e^{-\theta (t- u)}\sqrt{{\alpha ^3}{\sigma ^2} + {\beta ^3}{\varsigma ^2} } dW(u).
\ees
\item[\rm (iv)] The unique limiting distribution (also stationary distribution) of $\hatx$ is a normal distribution with mean $c/\theta$ and variance $({\alpha ^3}{\sigma ^2} + {\beta ^3}{\varsigma ^2})/2\theta$.
\ei
\end{lemma}

Using Lemma \ref{o-u} and Ito's formula, we can study the moment of the $k^{th}$ order of $Z(t)$ and $\hatx(t)$ provided that $\E(|Z(0)|^k)<\infty$ and $\E[(\hatx(0))^k] <\infty.$ For example, assuming that $\E(|Z(0)|^2)<\infty$ and $\E[(\hatx(0))^2]<\infty,$
\be\label{diffusion-moments}\ba
\E(Z(t)) & = \E(Z(0)) e^{-\theta t}, \\
\E(Z(t))^2 & =\E[(Z(0))^2] e^{-2\theta t} + \int_0^t e^{-2\theta (t-u)} [ {\alpha ^3}{\sigma ^2} + {\beta ^3}{\varsigma ^2} + \theta x^+(u) + \gamma x^-(u) ] du, \\
\E(\hatx(t)) & = \left(\E(\hatx(0)) - \frac{c}{\theta}\right) e^{-\theta t} + \frac{c}{\theta}, \\
\E(\hatx(t))^2 & = \left(\E[(\hatx(0))^2] - \frac{2c}{\theta}\left(\E(\hatx(0)) -\frac{c}{\theta}\right) -\left( \frac{c}{\theta}\right)^2 - \frac{\alpha^3\sigma^2+\beta^3\varsigma^2}{2\theta}\right) e^{-2\theta t}  \\
&\quad + \frac{2c}{\theta}\left(\E(\hatx(0)) -\frac{c}{\theta}\right) e^{-\theta t}+ \left( \frac{c}{\theta}\right)^2+ \frac{\alpha^3\sigma^2+\beta^3\varsigma^2}{2\theta}.
\ea\ee

When $\theta\neq \gamma$, $Z$ or $\hatx$ cannot be solved directly. In the following, we consider the limiting distributions of $Z$ and $\hatx$, and compute the limiting moments. Let $a =\sqrt{ {\alpha ^3}{\sigma ^2} + {\beta ^3}{\varsigma ^2}}$ and $b= \sqrt{a^2 + |\alpha - \beta|},$ and denote by $\phi(\cdot; \xi, \eta)$ and $\Phi(\cdot; \xi, \eta)$ the density and distribution functions of $N(\xi, \eta)$.

\begin{theorem}\label{limitphi}  Define 
\begin{equation}\label{psi}
\psi(x; \kappa, \mu, \sigma) = \left\{ \begin{array}{l}
\frac{C}{\sqrt{\theta}}\exp \left\{\frac{\kappa}{\theta}
\right\}\phi \left( {x;\frac{\mu}{\theta },\frac{\sigma^2}{2\theta}}
\right),\;\;\;x \ge 0,\\
\frac{C}{\sqrt{\gamma}}\exp \left\{\frac{\kappa}{\gamma}
\right\}\phi \left( {x;\frac{\mu}{\gamma },\frac{{{\sigma^2}}}{{2\gamma }}}
\right),\;\;\;x < 0,
\end{array} \right.
\end{equation}
where $C$ is given by
\begin{equation}\label{C}
C = C(\kappa, \mu, \sigma) = \frac{1}{{ \frac{1}{\sqrt{\theta}} \exp \left\{
{\frac{\kappa}{{\theta}}} \right\}\left( {1 - \Phi \left(
{0;\frac{\mu}{\theta },\frac{{{\sigma^2}}}{{2\theta }}} \right)} \right) +
\frac{1}{\sqrt{\gamma}}\exp \left\{ {\frac{\kappa}{{\gamma}}}
\right\} { \Phi \left( {0;\frac{\mu}{\gamma },\frac{{{\sigma^2}}}{{2\gamma
}}} \right)} }}.
\end{equation}
Then 
\bi
\item[\rm(i)] the density of the unique limiting distribution of the diffusion process $Z$ is given by
\begin{equation}\label{psi-Z}
\psi_Z (x) = \psi(x; 0, 0, b),
\end{equation}
\item[\rm (ii)] the density of the unique limiting distribution (also stationary distribution) of the diffusion process $\hatx$ is given by
\begin{equation}\label{psi-X}
\psi_{\hatx} (x) = \psi(x; c^2/a^2, c, a).
\end{equation}

\ei
\end{theorem}

In the following, we calculate the first two moments of the distribution $\psi(x;\mu^2/\sigma^2,\mu,\sigma), x\in\RR$. First, note that if $X \sim N(\xi, \eta)$, the density of a
truncated normal random variable on $(x_1,x_2)$ is given by$\frac{{\phi
(x;\xi ,{\eta})}}{{\Phi (x_2;\xi ,{\eta}) - \Phi (x_1;\xi
,{\eta})}}$. Let $X_1$ be a truncated $N(\frac{\mu}{\theta}, \frac{\sigma^2}{2\theta})$ random variable on $(0,+\infty)$, and $X_2$ be a truncated $N(\frac{\mu}{\gamma}, \frac{\sigma^2}{2\gamma})$ random variable on $(-\infty,0)$. Let $V$ be a mixture of these two truncated normal random variables as follows:
\begin{equation*}
V = \left\{ \begin{array}{ll} X_1 & \;\; \mbox{w.p. $d_1,$}\\
X_2 & \;\; \mbox{w.p. $d_2,$}
\end{array} \right.
\end{equation*}
 where
\begin{eqnarray}
{d_1} &=& \frac{C}{{\sqrt \theta  }}\exp \left\{ {\frac{{{\mu^2}}}{{\theta {\sigma^2}}}} \right\}\left( {1 - \Phi
\left( {0;\frac{\mu}{\theta },\frac{{{\sigma^2}}}{{2\theta }}} \right)} \right),\label{d1}\\
{d_2} &=&
\frac{C}{{\sqrt \gamma }}\exp \left\{ {\frac{{{\mu^2}}}{{\gamma {\sigma^2}}}} \right\}\Phi \left(
{0;\frac{\mu}{\gamma },\frac{{{\sigma^2}}}{{2\gamma }}} \right), \label{d2}
\end{eqnarray}
and $C=C(\mu/\sigma^2, \mu,\sigma)$ is defined as in \eqref{C}. Then $V$ has density function $\psi(x;\mu^2/\sigma^2,\mu,\sigma), x\in\RR$. Now, the first and second moments of the truncated normals are given by
\begin{eqnarray*}
\E(X_1) &=& \frac{\mu}{\theta } + \frac{\sigma}{{\sqrt {2\theta } }}
\frac{{\phi \left( { - \frac{\mu}{\sigma}\sqrt
{\frac{2}{\theta }} }; 0,1 \right)}}{{1 - \Phi \left( { - \frac{\mu}{\sigma}
\sqrt {\frac{2}{\theta }} }; 0, 1 \right)}},\\
\E(X_2) &=& \frac{\mu}{\gamma } - \frac{\sigma}{{\sqrt {2\gamma } }}
\frac{{\phi \left( { - \frac{\mu}{\sigma}\sqrt {\frac{2}{\gamma }} }; 0,1
\right)}}{{\Phi \left( { - \frac{\mu}{\sigma}\sqrt {\frac{2}{\gamma }} }; 0, 1
\right)}}, \\
\E(X_1^2) &=& {\left( {\frac{\mu}{\theta }} \right)^2} +
\frac{{{\sigma^2}}}{{2\theta }} + \frac{{\sqrt 2 }}{2}\frac{{\mu\sigma}}{{\theta
\sqrt \theta  }}\frac{{\phi \left( { - \frac{\mu}{\sigma}\sqrt
{\frac{2}{\theta }} }; 0 ,1
\right)}}{{1 - \Phi \left( { - \frac{\mu}{\sigma}\sqrt {\frac{2}{\theta }} }; 0,1
 \right)}},\\
\E(X_2^2) &=& {\left( {\frac{\mu}{\gamma }} \right)^2} +
\frac{{{\sigma^2}}}{{2\gamma }} - \frac{{\sqrt 2 }}{2}\frac{{\mu\sigma}}{{\gamma
\sqrt \gamma  }}\frac{{\phi \left( { - \frac{\mu}{\sigma}\sqrt
{\frac{2}{\gamma }} }; 0,1 \right)}}{{\Phi \left( { - \frac{\mu}{\sigma}\sqrt
{\frac{2}{\gamma }} }; 0,1 \right)}}.
\end{eqnarray*}
Hence the first and second moments of $V$ are given by
\begin{equation}\label{EV}
\E(V) = d_1 \E(X_1) + d_2 \E(X_2),
\end{equation}
and
\begin{equation}\label{EV2}
\E(V^2) = d_1 \E(X_1^2) + d_2 \E(X_2^2),
\end{equation}
where $d_1$ and $d_2$ is given by equation \eqref{d1} and \eqref{d2}. Note that when $\theta = \gamma,$
\be\label{EV-same}
\psi(x;\mu^2/\sigma^2,\mu,\sigma) = \phi\left(x; \frac{\mu}{\theta}, \frac{\sigma^2}{2\theta}\right),
\ee
and
\be\label{EV2-same}
\E(V) = \frac{\mu}{\theta}, \ \E(V^2) = \frac{\mu^2}{\theta^2} + \frac{\sigma^2}{2\theta}.
\ee
Finally, the first two limiting moments of $Z$ and $\hatx$ can be given by \eqref{EV} and \eqref{EV2} by replacing $(\mu, \sigma)$ with $(0, b)$ and $(c,a)$, respectively. When $\theta=\gamma$, \eqref{EV-same} and \eqref{EV2-same} give the same results as in Lemma \ref{o-u}.

\subsection{Interchange of the heavy traffic and steady state limits}

Let $A^n_s(t)$ denote the residual time of the interarrival time of sellers at time $t$, i.e., the time between $t$ and the first arrival after time $t$, and $A^n_b(t)$ is defined as the residual time of the interarrival time of buyers at time $t$. Denote by $\mathbb{X}^n$ the triplet $(X^n, A^n_s, A^n_b).$ Recall from Section \ref{sec:PA} and the beginning of Section \ref{sec:DA} that $\{U^n_k: k\in\NN\}$ and $\{V^n_k: k\in\NN\}$ are the sequences of interarrival times of sellers and buyers, respectively. For simplicity, we let $A^n_s(0) = U^n_1$ and $A^n_b(0) = V^n_1.$ Note that all the results in this section hold even if the distributions of $A^n_s(0)$ and $A^n_b(0)$ are different from those of $U^n_k$ and $V^n_k, k\ge 2.$ One can check that $\mathbb{X}^n$ is a Markov process with state space $\mathbb{Z}\times \RR_+\times \RR_+.$ Define the diffusion scaled process 
\[
\hat{\mathbb{X}}^n(t) = \frac{\mathbb{X}^n(nt)}{\sqrt{n}}, \ t\ge 0.
\]

\begin{theorem}\label{limit-change}
For $n\in\NN,$ $\hat{\mathbb{X}}^n$ is positive recurrent and admits a stationary distribution $\Pi^n$. Let $\pi^n$ be the marginal distribution of the first coordinate, i.e. $\pi^n(B) = \Pi^n(B\times \RR_+\times \RR_+)$ for $B\subset \{\frac{x}{\sqrt{n}}: x\in \mathbb{Z}\}.$ Then $\pi^n \Go \psi_{\hatx},$ where $\psi_{\hatx}$ is the limiting distribution (also the stationary distribution) of $\hatx$ as in Theorem \ref{limitphi}.

\end{theorem}

\begin{remark} \hfill
\bi
\item[\rm (i)] Theorem \ref{limit-change} essentially establishes the following interchange of limits 
\bes
\lim_{t\to\infty}\lim_{n\to\infty} \hatx^n(t) = \lim_{n\to\infty}\lim_{t\to\infty} \hatx^n(t),
\ees
which shows the validity of the heavy traffic steady state approximation for double-ended queues (see the numerical examples in the next section). 
\item[\rm (ii)] We only establish the interchange limit theorem for the heavy traffic diffusion approximation. However, following from \eqref{relation-X-Z} in Remark \ref{remark-diff}, under Assumption \ref{main-d} and the heavy traffic condition \eqref{htc}, the interchange limit result also holds for $Z^n.$
\ei
\end{remark}

\begin{remark}\label{compare-remark} Recall that in Section 3, $m(t), s(t)$ are the first two moments of the queue length process $X(t)$ of the double-ended queue with Poisson arrivals. At the end of this section, we compare $m(t), s(t)$ with proper estimators constructed from the fluid and diffusion limits. It will be seen that, under Assumption \ref{main-d} (additional heavy traffic condition needed for heavy traffic diffusion model), the fluid and diffusion models provide simple and good approximations for $\lim_{t\to\infty}m(t)$ and $\lim_{t\to\infty}s(t).$
\bi
\item[\rm (i)] We first compare $m(t)$ with the fluid limit $x(t)$ obtained in Theorem \ref{fluid-d} and Lemma \ref{fluid-s}. When $\theta = \gamma,$ the fluid equation is the same as the ODE for $m(t)$. However, when $\theta \neq \gamma$, noting that  in general $m_+(t) \neq m^+(t)$ and $m_-(t) \neq m^-(t)$, the fluid equation doesn't match with the ODE for $m(t).$ However, we note that when $\alpha, \beta$ are much larger than $\gamma, \theta$, and $\alpha$ and $\beta$ have the same order (see Assumption \ref{main-d}), using (8) in \cite{Amore05} and Stirling approximation,
\be\label{asym}
p_1 \approx \sqrt{\frac{\beta}{\theta}} e^{\alpha/\theta} \left(\frac{\beta}{\alpha e}\right)^{\beta/\theta}, \ \mbox{and} \ p_2 \approx \sqrt{\frac{\alpha}{\gamma}} e^{\beta/\gamma} \left(\frac{\alpha}{\beta e}\right)^{\alpha/\gamma}.
\ee
Then the limiting mean in Lemma \ref{simp} has the following approximations. When $\alpha > \beta$,
\be\label{asym-first-1}
\lim_{t\to\infty}m(t) = \frac{\alpha -\beta}{\theta} \frac{p_1 + \frac{\theta}{\gamma} p_2 + \frac{\alpha}{\alpha-\beta} - \frac{\theta}{\gamma}\frac{\beta}{\alpha-\beta}}{1+p_1+p_2} \approx \frac{\alpha -\beta}{\theta} = \lim_{t\to\infty} x(t),
\ee
and when $\alpha < \beta,$
\be\label{asym-first-2}
\lim_{t\to\infty}m(t) = \frac{\beta-\alpha}{\gamma} \frac{p_2 + \frac{\gamma}{\theta} p_1 + \frac{\beta}{\beta-\alpha} - \frac{\gamma}{\theta}\frac{\alpha}{\beta-\alpha}}{1+p_1+p_2} \approx \frac{\beta-\alpha}{\gamma} = \lim_{t\to\infty} x(t).
\ee

\item[\rm (ii)] We compare $\var(X(t)) =s(t) - (m(t))^2$ with the second moment of the diffusion limit $Z(t)$ in Theorem \ref{diffusion-d}. For simplicity, assume $\theta=\gamma.$ From \eqref{MC-sec-moment} and Lemma \ref{simp},
\[
\lim_{t\to\infty} s(t) = \frac{\alpha+\beta}{2\theta} + \frac{(\alpha-\beta)^2}{\theta^2} + \frac{ \lim_{t\to\infty}m_+(t) +  \lim_{t\to\infty} m_-(t)}{2},
\]
and 
\[
\lim_{t\to\infty}m_+(t) +  \lim_{t\to\infty} m_-(t) = \left[\frac{\alpha-\beta}{\theta}p_1 + \frac{\alpha}{\theta} + \frac{\beta-\alpha}{\gamma}p_2 + \frac{\beta}{\gamma}\right]\pi_0.
\]
Under the asymptotic regime in Assumption \ref{main-d}, from \eqref{asym-first-1} and \eqref{asym-first-2}, we have
\[
\lim_{t\to\infty}m_+(t) +  \lim_{t\to\infty} m_-(t)  \approx \frac{|\alpha-\beta|}{\theta}.
\]
Thus for the queueing system with Poisson arrivals with parameters in the asymptotic regime in Assumption \ref{main-d}, 
\be\label{s1-ode-soln}
\lim_{t\to\infty} s(t) \approx \frac{\alpha+\beta}{2\theta} + \frac{(\alpha-\beta)^2}{\theta^2} + \frac{|\alpha-\beta|}{2\theta} =  \frac{(\alpha-\beta)^2}{\theta^2}+\frac{\max\{\alpha, \beta\}}{\theta},
\ee
and so combining \eqref{s1-ode-soln} with \eqref{asym-first-1} and \eqref{asym-first-2}, we have
\be\label{limit-variance}
\lim_{t\to\infty}\var(X(t)) \approx \frac{\max\{\alpha, \beta\}}{\theta}.
\ee
We now consider the diffusion limit $Z$. From \eqref{diffusion-moments}, we have
\[
\lim_{t\to\infty} \E(Z(t))^2 = \frac{{\alpha ^3}{\sigma ^2} + {\beta ^3}{\varsigma ^2} +|\alpha-\beta|}{2\theta}.
\]
When the arrivals are Poisson processes, the variances of interarrival times become $\alpha$ and $\beta$, and so
\be
\lim_{t\to\infty} \E(Z(t)^2) = \frac{\alpha+\beta +|\alpha-\beta|}{2\theta} = \frac{\max\{\alpha,\beta\}}{\theta},
\ee
which is the same as \eqref{limit-variance}.  

\item[\rm (iii)] We finally consider the heavy traffic diffusion approximation $\hatx$ in Theorem \ref{diff-ht}. Again we assume $\theta = \gamma.$ We consider a sequence of double-ended queues with Poisson arrivals, indexed by $n\in\NN.$ So the parameters of the $n$-th system are $\alpha^n, \beta^n, \theta^n$. Under Assumption \ref{main-d} and heavy traffic condition \eqref{htc}, we have 
\be\label{all-cond}
\alpha^n \to \alpha, \ \beta^n \to \beta, \ \sqrt{n}(\alpha^n-\beta^n) \to c, \ n\theta^n \to \theta, \ n\gamma^n \to\gamma. 
\ee
From Theorem \ref{limit-change}, we have
\bes\ba
\hatx(\infty) = \lim_{t\to\infty}\hatx(t) = \lim_{t\to\infty} \lim_{n\to\infty} \frac{X^n(nt)}{\sqrt{n}} = \lim_{n\to\infty} \lim_{t\to\infty} \frac{X^n(nt)}{\sqrt{n}} = \lim_{n\to\infty} \frac{X^n(\infty)}{\sqrt{n}}.
\ea\ees
Thus for large enough $n\in\NN,$ we have 
\bes
X^n(\infty) \approx \sqrt{n} \hatx(\infty),
\ees
and so using Lemmas \ref{limit-moments} and \ref{o-u} (iv), 
\be\label{ht-m1}
\lim_{t\to\infty} \E(X^n(t)) = \E(X^n(\infty)) \approx \sqrt{n} \E(\hatx(\infty)) =  \frac{\sqrt{n} c}{\theta},
\ee
and 
\be\label{ht-s1}
\lim_{t\to\infty}\E[(X^n(t))^2] = \E(X^n(\infty))^2 \approx n \E(\hatx(\infty))^2 = \frac{n c^2}{\theta^2} + \frac{n(\alpha+\beta)}{2\theta}.
\ee
From the convergence in \eqref{all-cond}, we have for large $n\in\NN,$
\be\label{ht-ms}
\frac{\sqrt{n} c}{\theta} \approx \frac{\alpha^n - \beta^n}{\theta^n}, \ \frac{n c^2}{\theta^2} \approx \left(\frac{\alpha^n - \beta^n}{\theta^n}\right)^2, \ \mbox{and} \  \frac{n(\alpha+\beta)}{2\theta} \approx \frac{\alpha^n + \beta^n}{2\theta^n}.
\ee
Combining \eqref{ht-m1}, \eqref{ht-s1}, and \eqref{ht-ms}, we have for large $n\in\NN,$
\be\label{ht-approx-1}
\lim_{t\to\infty} \E(X^n(t))  \approx  \frac{\alpha^n - \beta^n}{\theta^n},
\ee
and 
\be\label{ht-approx-2}
\lim_{t\to\infty}\E(X^n(t)^2) \approx \left(\frac{\alpha^n - \beta^n}{\theta^n}\right)^2 + \frac{\alpha^n + \beta^n}{2\theta^n} \approx \left(\frac{\alpha^n - \beta^n}{\theta^n}\right)^2 + \frac{\max\{\alpha^n, \beta^n\}}{\theta^n},
\ee
where the last approximation follows from the heavy traffic condition, i.e. $\alpha^n$ and $\beta^n$ are roughly equal when $n$ is large.

Consider now a double-ended queueing system with parameters $\alpha, \beta, \theta$. Suppose the parameters satisfy Assumption \ref{main-d} and heavy traffic condition \eqref{htc}, i.e. the arrival rates $\alpha$ and $\beta$ are very close to each other and the reneging rates $\theta$ and $\gamma$ are very small comparing with $\alpha, \beta$ and $|\alpha-\beta|$. The above approximations in \eqref{ht-approx-1} and \eqref{ht-approx-2} say that the first two stationary moments of the queue length process can be approximated by 
\bes
\frac{\alpha - \beta}{\theta}, \ \mbox{and} \ \left(\frac{\alpha - \beta}{\theta}\right)^2 + \frac{\max\{\alpha, \beta\}}{2\theta},
\ees
which are the same as the approximations in \eqref{asym-first-1}, \eqref{asym-first-2}, and \eqref{s1-ode-soln}. 

\ei

\end{remark}

\section{Numerical Examples} \label{sec:NE}

In this section, we use simulations to evaluate  the performance of the Poisson and the diffusion approximations under different arrival processes. We now consider a double-ended queue length process $\{X(t): t\geq0\}$ with seller inter-arrival time distribution $F_s(\cdot)$ and buyer inter-arrival time distribution $F_b(\cdot)$. Let $m_s$, $m_b$, $sd_s$ and $sd_b$ be the means and standard deviations of the inter-arrival time for sellers and buyers. We consider the following inter-arrival time distributions:

\begin{itemize}
\item Exponential: $F_s(x) = 1 - e^{-\alpha x}$,
$F_b(x) = 1 - e^{-\beta x}$, $m_s=\frac{1}{\alpha}$,
$m_b=\frac{1}{\beta}$, $sd_s=\frac{1}{\sqrt{\alpha}}$, $sd_b=\frac{1}{\sqrt{\beta}}$
\item Uniform: $F_s(x) = \frac{\alpha x}{2}\;
(x\in[0,\frac{2}{\alpha}])$, $F_b(x) = \frac{\beta
x}{2}\;(x\in[0,\frac{2}{\beta}])$, $m_s=\frac{1}
{\alpha}$, $m_b=\frac{1}{\beta}$,
$sd_s=\frac{1}{\sqrt{3}\alpha}$, $sd_b=\frac{1}{\sqrt{3}\beta}$
\item Erlang(2): $F_s(x) = 1 - e^{-2\alpha x} -
2\alpha x e^{-2\alpha x}$, $F_b(x) = 1 - e^{-2\beta x} - 2\beta x e^{-2\beta x}$, $m_s=\frac{1}{\alpha}$,
$m_b=\frac{1}{\beta}$, $sd_s=\frac{1}{\sqrt{2}\alpha}$, $sd_b=\frac{1}{\sqrt{2}\beta}$
\item Hyper-exponential: $F_s(x) = \frac{1}{3}(1 -
e^{-\frac{1}{2}\alpha x}) + \frac{2}{3}(1 - e^{-2 \alpha x})$,
$F_b(x) = \frac{1}{3}(1 - e^{-\frac{1}{2}\beta x})  + \frac{2}{3}(1
- e^{-2 \beta x})$, $m_s=\frac{1}{\alpha}$, $m_b=\frac{1}{\beta}$,
$sd_s=\frac{\sqrt{2}}{\alpha}$, $sd_b=\frac{\sqrt{2}}{\beta}$
\end{itemize}

We consider the following arrival rates $(\alpha, \beta) =$  $(1, 1)$,
$(1, 1.5)$ and $(1, 2)$, and choose the following reneging rates
$(\theta, \gamma) =$ $(\alpha, \beta)$, $0.1(\alpha, \beta)$ and
$0.01(\alpha, \beta)$. For example, when $(\theta, \gamma) = 0.1(\alpha,
\beta)$, the sellers' (buyers') expected patience time is 10 times
the sellers' (buyers') expected inter-arrival time. Thus we consider
a total of $4\times3\times3 = 36$ different parameter settings. Also note that the means of inter-arrival time of the above distributions are the same, while their standard
deviations are different, with the following ordering: Hyper-exponential $>$ Exponential $>$ Erlang $>$
Uniform.

In each parameter setting, we use simulation, Poisson approximation, and two diffusion approximations (which will be made precise below) to estimate the following performance measures:
\be\label{moments}\ba
{L_1} &= \mathop {\lim }\limits_{t \to \infty } E(X(t)), \ \ {L_2} &= \mathop {\lim }\limits_{t \to \infty } E[(X{(t))^2}].
\ea\ee

{\bf Simulation.} We compute the performance measure by using $N$  replications of
the simulation by Matlab. Each replication consists of simulating
the system for $0 \leq t \leq T$ and the estimates are computed by
using the sample paths over $t \in [\tau, T]$, where $\tau<T$ is a
given warmup period. Let $X^k(t)$ be the state of the system at time
$t$ in the $k$-th replication, $k = 1,2,\cdots,N$, $0\leq t\leq T$.
Using these sample paths, we compute
\begin{equation*}
 \pi _i^s = \frac{1}{N}\sum\limits_{k = 1}^N {\frac{1}{{T - \tau }}}
\int\limits_\tau ^T {{1_{\left\{ {{X^k}(t) = i} \right\}}}} {\rm{
}}dt,\;\;\;\;  - 1000 \le i \le 1000.
\end{equation*}
Using these we compute the following simulation estimates of the
first and second moments of the queue length:
\be\ba
L_1^s &= \sum\limits_{i =  - 1000}^{1000} {i\pi _{_i}^s}, \ \ L_2^s = \sum\limits_{i =  - 1000}^{1000} {{i^2}\pi _{_i}^s}.
\ea\ee

{\bf Poisson approximation.} In this approximation we replace the renewal
arrival processes by Poisson arrival processes with the same arrival
rates. Clearly, this approximation is exact in the exponential case. Let $L_1^p$ and $L_2^p$ be the Poisson  approximation of $L_1$ and
$L_2$ respectively. From equations (\ref{limitdis1}) -- (\ref{limitdis3}), and Lemma \ref{simp}, we have:
\begin{equation}\label{p-approx-1}
\ba
L_1^p &=  \left[\frac{\alpha-\beta}{\theta} p_1 +\frac{\alpha}{\theta}-\frac{\beta-\alpha}{\gamma} p_2 -\frac{\beta}{\gamma}\right]\pi_0, \\
L_2^p &=  \left[\frac{\alpha-\beta}{\theta} \left(\frac{\alpha-\beta}{\theta} p_1 +\frac{\alpha}{\theta}\right) + \frac{\alpha}{\theta}(p_1 +1)+\frac{\beta-\alpha}{\gamma} \left(\frac{\beta-\alpha}{\gamma} p_2 +\frac{\beta}{\gamma}\right) + \frac{\beta}{\gamma}(p_2 +1)\right]\pi_0,
\ea
\end{equation}
where 
\bes\ba
p_1 & = \frac{\beta}{\theta} e^{\alpha/\theta}  \left(\frac{\alpha}{\theta}\right)^{-\beta/\theta} \Gamma(\beta/\theta, \alpha/\theta)- 1,\\
p_2 & = \frac{\alpha}{\gamma}  e^{\beta/\gamma} \left(\frac{\beta}{\gamma}\right)^{-\alpha/\gamma} \Gamma(\alpha/\gamma, \beta/\gamma)- 1,\\
\pi_0 & = (1+p_1+p_2)^{-1}.
\ea\ees

We compute the relative error of the above moments to the ones from
simulation method. To be precise, the relative error of $L_1^p, L_2^p$ to
$L_1^s, L_2^p$ are given by $\left(|L_1^p - L_1^s|/L_1^s\right) \times 100\%$ and $\left(|L_2^p - L_2^s|/L_2^s\right) \times 100\%$, respectively.

{\bf Two diffusion approximations (see Section \ref{diff-model} for details).} We consider two diffusion models arising from Theorems \ref{diffusion-d} and \ref{diff-ht}. The first diffusion model is constructed from Theorem \ref{diff-ht}. We note that the arrival rates should be roughly equal in the heavy traffic regime. However, the first diffusion model can be applied for any parameter regimes. More precisely, we approximate $L_1$ and $L_2$ (defined in \eqref{moments}) by using the stationary moments of $X_1,$ where $X_1$ is defined in \eqref{approx-model-1}, i.e.,
\begin{equation*}
X_1(t) = X(0) + \sqrt{\alpha^3\sigma^2+\beta^3\varsigma^2}  W(t) + \left( {\alpha  - \beta }
\right)t - \theta \int_0^t X_1^+(u)du  + \gamma \int_0^t X_1^-(u)du, \ t\ge 0.
\end{equation*} 
We compute the first two limiting moments of $X_1$ using \eqref{EV} and \eqref{EV2} with $(\mu,\sigma)$ replaced by $(\alpha-\beta, \sqrt{ {\alpha ^3}{\sigma ^2} + {\beta ^3}{\varsigma ^2}})$, and denote them by $L_1^{d,1}$ and $L_2^{d,1}$. Namely, 
\be\label{d-approx-1}
L_1^{d,1} = \E\left( \lim_{t\to\infty} X_1(t)\right), \ L_2^{d,1} = \E\left(\lim_{t\to\infty} (X_1(t))^2\right). 
\ee
The second diffusion model is constructed from Theorems \ref{fluid-d} and \ref{diffusion-d}. We directly use the limiting point of $x_2$ defined in \eqref{f-approx} (which is the same as the fluid limit in \eqref{limit-fluid}) to approximate the limiting mean $L_1$, and use the second limiting moment of $Z_2$ defined in \eqref{d-approx} to approximate the limiting variance of $X$, namely, $L_2 - L_1^2.$ More precisely, we have  
\bes
x_2(t) =  (\alpha - \beta) t - \theta \int_0^t x^+_2(s)ds + \gamma \int_0^t x^-_2(s)ds, \ t\ge 0,
\ees
\bes
Z_2(t)= X(0) + \int_0^t \sqrt{{\alpha ^3}{\sigma ^2} + {\beta ^3}{\varsigma ^2} + \theta x^+_2(u) + \gamma x^-_2(u) } dW(u)  - \theta\int_0^t   Z_2^+(u) du +\gamma\int_0^t  Z_2^-(u)du, \ t\ge 0,
\ees
and the approximations are as follows:
\be\label{d-approx-2}
L_1^{d,2} = \lim_{t\to\infty} x_2(t), \ L_2^{d,2} = \left(\lim_{t\to\infty} x_2(t)\right)^2 + \E\left(\lim_{t\to\infty} Z_2(t)^2\right). 
\ee
We compute the second limiting moment of $Z_2$ using \eqref{EV2} with $(\mu,\sigma)$ replaced by $(0,b)$. We also compute the relative error of the above moments to the ones from simulation method.

The comparisons of the approximations of $L_1$ and $L_2$ are shown in Tables 1-8. The results about the first moment are shown in the tables 1-4 and those about the
second moment are shown in tables 5-8. In the columns of $L_1^p$ and $L_2^p$ (see \eqref{p-approx-1}), we evaluate the performance measures by Poisson approximation method, and obtain the relative error of each performance measure to the one from simulation method. In the columns of $L_1^{d,1}, L_1^{d,2}$ and $L_2^{d,1}, L_2^{d,2}$ (see \eqref{d-approx-1} and \eqref{d-approx-2}), we obtain the performance measures by diffusion approximation methods, and also obtain the relative error of each performance measure to the one from simulation method. The comparisons of limiting density are shown in Figure 2-5. In the figures, we compare the density graphs derived from the simulation method, the stationary distribution of the Poisson model, and the stationary distribution of the heavy traffic diffusion model. When using simulation method, we evaluate the performance measures using the parameter $(N, \tau, T) = (400, 1000, 4000)$ and obtain the $90\%$ confidence interval.

From the numerical examples, we have the following conclusions. 
\bi
\item From Table 1-4, consider the limiting mean $L_1$ and its approximations $L_1^p, L_1^{d,1}$, and $L_1^{d,2}$.
\bi
\item Both diffusion approximations $L_1^{d,1}$ and $L_1^{d,2}$ improve when the reneging rates become smaller.
\item The heavy traffic diffusion approximation $L_1^{d,1}$ behaves better than the Poisson approximation in all non-exponential cases. 
\item The fluid approximation $L_1^{d,2}$ becomes very close to the heavy traffic diffusion approximation $L_1^{d,1}$, when the reneging rates are small, which in particular suggests that {\em the fluid limit is a good approximation of the transient mean of the heavy traffic diffusion limit}.
\ei
\item From Table 5-8, consider the second limiting moment $L_2$ and its approximations $L_2^p, L_2^{d,1}$, and $L_2^{d,2}$.
\bi
\item Both diffusion approximations $L_2^{d,1}$ and $L_2^{d,2}$ improve when the reneging rates become smaller.
\item The heavy traffic diffusion approximation $L_2^{d,1}$ behaves better than the Poisson approximation in all non-exponential cases. 
\item The diffusion approximation $L_2^{d,2}$ is not as good as the heavy traffic diffusion approximation $L_2^{d,1}$. However, as the reneging rates get smaller, $L_2^{d,2}$ gets closer to $L_2^{d,1}.$ It seems that the convergence rate of $L_2^{d,2}$ is smaller than that of $L_2^{d,1}$. Roughly speaking, the diffusion process $Z$ is centered at the fluid limit $x$ (see \eqref{centered}). From Table 1-4, we observe that the fluid first moment approximation $L_1^{d,2}$ behaves well only when the reneging rates are small, which make the second moment approximation $L_2^{d,2}$ works well only with small reneging rates. 
\ei
\ei
In summary, when the interarrival times are not exponential distributed, the heavy traffic diffusion model performs well for general parameters. In particular, Theorem \ref{limit-change} guarantees the validity of such approximation. When the reneging rates are small comparing with the arrival rates, the fluid limit approximation $L_1^{d,2}$ is a good simple approximation of the limiting mean $L_1.$

%%%%%-------------------------------section 7----------------------------------------------------------
\section{Extensions} \label{sec:C}

We end this paper with suggestions for three extension.
\bi
\item[1.] In this paper we have assumed that the patience times of buyers and sellers are exponentially distributed. It would be interesting to study the situation when the distributions are general, and to establish diffusion approximations for $\{X(t): t \ge 0\}$ under similar parameter regime.

\item[2.] In this work we assume that the arrival processes of buyers and sellers are independent of the state of the system. It would be interesting to consider an extension where the arrival processes are counting processes whose intensity parameters depend on the state of the
double-ended queue. For example, the parameters could simply depend upon the sign of $X(t)$. Thus we could capture the situation where the arrival rate of the buyers exceeds that of the sellers when there are sellers waiting (it's a buyers' market), and the arrival rate of the sellers exceeds that of the buyers when there are buyers waiting (it's a sellers' market). We think this extension is doable, and the methods developed in this paper will be useful in the study.

\item[3.] In Cont et al. \cite {Cont10}, the dynamics of a limit order book (i.e. the number of limit ask and bid orders) with $n$ different prices is modeled as a $n$-dimensional double-ended queueing system with Markovian primitives. The authors then use simple matrix computations and Laplace transform methods to study interesting behaviors. Later on, Cont and de Larrard \cite{Cont12} propose a diffusion model for the limit order book by considering the joint dynamics of the (one-sided) bid and ask queues. We think it would be interesting to consider a (multidimensional) double-ended queue to model the bid and ask queues, and study the diffusion approximations. However, this extension promises to be hard, since the state-space of the multidimensional double-ended queue is typically not convex.

\ei

\section*{Acknowledgement}
We thank the anonymous referees for valuable comments and suggestions.

\bibliographystyle{plain}
\bibliography{references}

\newpage

\appendix

\section{Proofs} \label{appendix:proofs}

{\bf Proof of Lemma \ref{limit-moments}:} From \eqref{limitdis1} -- \eqref{limitdis3}, there exists a random variable $X(\infty)$ with distribution $\{\pi_i: i\in\mathbb{Z}\}$ such that $X(t) \Go X(\infty)$ as $t\to\infty.$ Using the continuous mapping theorem, we have 
\be\label{weak-con}
e^{s|X(t)|} \Go e^{s|X(\infty)|}, \ \mbox{as $t\to\infty.$}
\ee
In the following, we show $e^{s|X(t)|}$ is uniformly integrable. Define for $t\ge 0,$
\bes
\ba
Y_1(t) & = X^+(0) + N_s(t) - N_{sr}\left(\theta \int_0^t Y_1(s)ds\right), \\
Y_2(t) & = X^-(0) + N_b(t) - N_{br}\left(\theta \int_0^t Y_2(s)ds\right).
\ea
\ees
Noting that
\bes
\ba
X^+(t) & \le X^+(0) + N_s(t) - N_{sr}\left(\theta \int_0^t X^+(s)ds\right), \\
X^-(t) & \le X^-(0) + N_b(t) - N_{br}\left(\theta \int_0^t X^-(s)ds\right),
\ea
\ees
we have for $t\ge 0,$
\[
X^+(t) \le Y_1(t), \ \mbox{and} \ X^-(t) \le Y_2(t).
\]
We next observe that both $Y_1$ and $Y_2$ are birth and death processes on $\mathbb{N}\cup\{0\}$, with stationary distributions $\pi_1(j) = \frac{(\alpha/\theta)^j}{j!}, j\in \mathbb{N}\cup\{0\}$ and $\pi_2(j) = \frac{(\beta/\gamma)^j}{j!}, j\in \mathbb{N}\cup\{0\}$, respectively. From Section 8.4 (page 286) in \cite{anderson91}, and noting that $Y_1(t)$ and $Y_2(t)$ are independent given that $X(0)=x\in\mathbb{Z},$ we have 
\bes\ba
\E_x\left[e^{2s|X(t)|}\right] & \le \E_x\left[e^{2s[Y_1(t) + Y_2(t)]}\right] = \E_x\left[e^{2sY_1(t)}\right] \E_x\left[e^{2sY_2(t)}\right] \\
& = \left(1-e^{-\theta t} + e^{2s-\theta t}\right)^x e^{ (1- e^{-\theta t}) (e^{2s}-1)\alpha/\theta} \left(1-e^{-\gamma t} + e^{2s-\gamma t}\right)^x e^{ (1- e^{-\gamma t}) (e^{2s}-1)\beta/\gamma} \\
& \le (1+e^{2s})^{2x} e^{(\alpha/\theta + \beta/\gamma) (e^{2s}-1)} < \infty. 
\ea\ees
This shows the uniform integrability of $e^{s|X(t)|}$. Using \eqref{weak-con} and Theorem 5.4 in Chapter I.5 of \cite{billingsley1999convergence}, \eqref{mgf} follows. The convergence in \eqref{mfs} can be shown similarly. In particular, the uniform integrability of $f(X(t))$ follows from that of $e^{s|X(t)|}$. \ink

{\bf Proof of Lemma \ref{simp}:} From \cite{AnckerGafarian62}, we can simplify the sums in $\pi_0$ as follows.
\bes\ba
p_1 = \sum_{i=1}^\infty \frac{\alpha^i}{\prod_{j=1}^i (\beta + j\theta)} & = \sum_{i=1}^\infty \frac{(\alpha/\theta)^i}{\prod_{j=1}^i (\beta/\theta + j)} = \Gamma(\beta/\theta+1)\sum_{i=1}^\infty \frac{(\alpha/\theta)^i}{\Gamma(\beta/\theta + i+1)} \\
& =\frac{\beta}{\theta} e^{\alpha/\theta}  \left(\frac{\alpha}{\theta}\right)^{-\beta/\theta} \Gamma(\beta/\theta, \alpha/\theta)- 1,
\ea\ees
and by symmetry,
\bes
p_2 = \sum_{i=1}^\infty \frac{\beta^i}{\prod_{j=1}^i (\alpha + j\gamma)} = \frac{\alpha}{\gamma}  e^{\beta/\gamma} \left(\frac{\beta}{\gamma}\right)^{-\alpha/\gamma} \Gamma(\alpha/\gamma, \beta/\gamma)- 1.
\ees
Similarly, we have 
\bes\ba
m_1=\sum_{i=1}^\infty \frac{i\alpha^i}{\prod_{j=1}^i (\beta + j\theta)} & = \sum_{i=1}^\infty \frac{i (\alpha/\theta)^i}{\prod_{j=1}^i (\beta/\theta + j)} = \sum_{i=1}^\infty \frac{(\alpha/\theta)^i}{\prod_{j=1}^{i-1} (\beta/\theta + j)} - \beta/\theta \sum_{i=1}^\infty \frac{(\alpha/\theta)^i}{\prod_{j=1}^i (\beta/\theta + j)} \\
& = \frac{\alpha-\beta}{\theta} p_1 +\frac{\alpha}{\theta},
\ea\ees
and
\bes\ba
s_1=\sum_{i=1}^\infty \frac{i^2\alpha^i}{\prod_{j=1}^i (\beta + j\theta)} & = \sum_{i=1}^\infty \frac{i^2 (\alpha/\theta)^i}{\prod_{j=1}^i (\beta/\theta + j)}  \\
& = \sum_{i=1}^\infty \frac{i (\alpha/\theta)^i}{\prod_{j=1}^{i-1} (\beta/\theta + j)} - \beta/\theta \sum_{i=1}^\infty \frac{i (\alpha/\theta)^i}{\prod_{j=1}^i (\beta/\theta + j)}\\
& = \frac{\alpha-\beta}{\theta} m_1 + \frac{\alpha}{\theta}(p_1 +1). 
\ea\ees
Then by symmetry,
\bes\ba
m_2 = \sum_{i=1}^\infty \frac{i\beta^i}{\prod_{j=1}^i (\alpha + j\gamma)} & = \frac{\beta-\alpha}{\gamma} p_2 +\frac{\beta}{\gamma},
\ea\ees
and
\bes\ba
s_2 =\sum_{i=1}^\infty \frac{i^2\beta^i}{\prod_{j=1}^i (\alpha + j\gamma)}=   \frac{\beta-\alpha}{\gamma} m_2 + \frac{\beta}{\gamma}(p_2 +1). 
\ea\ees
The lemma follows, on noting that $\pi_0 = (1+p_1+p_2)^{-1}, \lim_{t\to\infty} m_+(t) = m_1\pi_0, \lim_{t\to\infty} m_-(t) = m_2\pi_0, \lim_{t\to\infty} s_+(t) = s_1\pi_0,$ and $\lim_{t\to\infty} s_-(t) = s_2\pi_0.$ \ink

{\bf Proof of Lemma \ref{th:mom}:} We first consider  $m(t) = \E(X(t))$. Taking
expectation of equation (\ref{eq:1}), we get
\begin{equation*}
m(t) = m(0) + \alpha t - \beta t - \theta \int_0^t  {{m_ + }} (s)ds
+ \gamma \int_0^t {{m_ - }} (s)ds.
\end{equation*}
Taking derivative on both sides of above equation  we get equation
(\ref{m1}).

Next we consider the second moment of $X(t)$. Using  the
infinitesimal analysis, for a small $h > 0$, we get,
\begin{equation*}
\left(X(t + h)\right)^2 = \left\{ \begin{array}{l}
\left( {X{{(t)}} + 1} \right)^2,\;\;\;w.p.\;\;\;
\left( {\alpha  + \gamma {X^ - }(t)} \right)h + o(h)\\
X{(t)^2},\;\;\;\;\;\;\;\;\;\;\;\;\;w.p.\;\;\;
1 -\left(\alpha  + \gamma {X^ - }(t) + \beta  + \theta
{X^ + }(t) \right)h + o(h) \\
\left( {X{{(t)}} - 1} \right)^2,\;\;\;w.p.\;\;\;\left(
{\beta  + \theta {X^ + }(t)} \right)h + o(h).
\end{array} \right.
\end{equation*}
Therefore,
\begin{eqnarray*}
E(\left(X(t + h)\right)^2|X(t)) &=& (X(t))^2 + 2X(t)\left( {\alpha  -
\beta  - \theta {X^ + }(t) + \gamma {X^ - }(t)} \right)h\\
&& + \left( {\alpha  + \beta  + \theta {X^ + }(t)  + \gamma {X^ -
}(t)} \right)h + o(h).
\end{eqnarray*}
Since $X(t)=X^+(t)-X^-(t)$, we have
\begin{eqnarray*}
E((X(t + h))^2 |X(t)) &=& (X(t))^2 - 2\theta (X^+(t))^2h -  2\gamma (X^ -
(t))^2h + (2\alpha  - 2\beta  + \theta
){X^ + }(t)h\\
&&+ ( - 2\alpha  + 2\beta  + \gamma ){X^ - }(t)h  + \alpha h  +
\beta h + o(h).
\end{eqnarray*}
Taking expectation on both sides of above equation, we get
\begin{eqnarray*}
\frac{{s(t + h) - s(t)}}{h} &=&  - 2\theta {s_ + }(t)  - 2\gamma {s_
- }(t) + (2\alpha  - 2\beta  + \theta ){m_
+ }(t)\\
&& + ( - 2\alpha  + 2\beta  + \gamma ){m_ - }(t) +  \alpha  + \beta
+ \frac{{o(h)}}{h}.
\end{eqnarray*}
Taking limit $h \to 0$, the equation (\ref{m2}) follows. \ink

{\bf Proof of Theorem \ref{fluid-d}:} We note from \eqref{queue} that for $t\geq 0$,
\bes
\barx^n(t)=\barx^n(0)+\barn_s^n(t)-\barn_b^n(t)-\barn_{sr}^n\left(n\theta^n \int_0^t
\barx^{n,+}(s)ds\right)+\barn_{br}^n\left(n\gamma^n\int_0^t
\barx^{n,-}(s)ds\right).
\ees
For $t\in [0,\infty),$ let $N^n(t) = N^n_{sr}(t) + N^n_{br}(t)$ and $O^n(t) = |\barn_s^n(t)-\alpha^nt
- \barn_b^n(t)+\beta^nt|$.
Then we have that
\bes\ba
|\barx^n(t)|  &\le |\barx^n(t)- (\barn_s^n(t)-\alpha^nt
- \barn_b^n(t)+\beta^nt)| + O^n(t) \\
&\quad \le |\barx^n(0)| +|\alpha^n-\beta^n|t + O^n(t) + n^{-1}N^n \left(n^2(\gamma^n+\theta^n) \int_0^{t} |\barx^{n}(u)|du \right).
\ea\ees
Define for $t\in [0,\infty),$
\bes
Y^n(t) = |\barx^n(0)|+|\alpha^n-\beta^n|t + O^n(t) + n^{-1}N^n \left(n^2(\gamma^n+\theta^n) \int_0^{t} Y^{n}(s) ds \right).
\ees
Then $$|\barx^n(t)| \le Y^n(t), \ t\in [0,\infty).$$
Noting that $N^n, O^n$ and $\barx^n(0)$ are mutually independent, we see that
$$M^n(t) := Y^n(t) -\barx^n(0) -|\alpha^n-\beta^n|t- O^n(t) - 2n(\gamma^n+\theta^n) \int_0^{t} Y^{n}(s)ds$$
is an $\{\clf^n_t\}$ martingale. Using Ito's formula, we have that
\bes\ba
& (Y^n(t) - O^n(t)) \exp\{-2n(\gamma^n+\theta^n)t\} = |\barx^n(0)| + \int_0^t \exp\{-2n(\gamma^n+\theta^n)s\} d M^n(s) \\
&\quad + \int_0^t 2n(\gamma^n+\theta^n) \exp\{-2n(\gamma^n+\theta^n)s\} O^n(s)ds + |\alpha^n - \beta^n| \int_0^t \exp\{-2n(\gamma^n+\theta^n)s\} ds,
\ea\ees
and so
\be\label{mart1}\ba
& \left(Y^n(t) + \frac{|\alpha^n-\beta^n|}{2n(\gamma^n+\theta^n)}\right)\exp\{-2n(\gamma^n+\theta^n)t\} - \left(|\barx^n(0)|+ \frac{|\alpha^n-\beta^n|}{2n(\gamma^n+\theta^n)}\right) \\
&=  \int_0^t \exp\{-2n(\gamma^n+\theta^n)s\}dO^n(s) + \int_0^t \exp\{-2n(\gamma^n+\theta^n)s\} d M^n(s).
\ea\ee
We observe that from the functional law of large numbers for renewal processes, $O^{n}\Go 0$ as $n\to\infty,$ and from the continuous mapping theorem, $\|O^n\| \Go 0$ as $n\to\infty.$ From Lemma 3.5 in \cite{budhirajaghosh06}, we have for some $c_1 \in (0,\infty)$ (independent of $n$ and $t$),
\bes
\E\left( \sup_{0\le u\le t} \left(|\hatn^n_s(u)|^2 + |\hatn^n_b(u) |^2\right)\right) \le c_1(t+1).
\ees
We have for $t\ge 0,$
\bes
\sup_{n\in\NN} \E(\|O^n\|^2_t)  \le \sup_{n\in\NN} \frac{2}{\sqrt{n}} \E\left( \sup_{0\le u\le t} \left(|\hatn^n_s(u)|^2 + |\hatn^n_b(u) |^2\right)\right) \le {2c_1(t+1)},
\ees
which implies the uniform integrability of $\{\|O^n\|_t: n\in\mathbb{N}\}.$
Thus we conclude that for $T\in [0,\infty),$
\be\label{con1}
\E\left(\sup_{0\le t\le T} O^n(t)\right) \to 0, \textt{as $n\to\infty.$}
\ee
 We next note that from \eqref{mart1} and \eqref{con1}, for any $t\in[0,\infty),$
\be\label{con2}
\E(Y^n(t)) \to \left(x_0+ \frac{|\alpha-\beta|}{2(\gamma+\theta)}\right) \exp\{2(\gamma+\theta)t\} - \frac{|\alpha-\beta|}{2(\gamma+\theta)}, \textt{as $n\to\infty.$}
\ee
From Doob's inequality and \eqref{con2}, for any $T\in[0,\infty),$
\be\label{con3}\ba
&\E\left(\sup_{0\le t\le T} \left|\int_0^t \exp\{-2n(\gamma^n+\theta^n)s\} d M^n(s) \right| \right)^2 \le 4 \E\left(\int_0^T \exp\{-2n(\gamma^n+\theta^n)s\} d M^n(s) \right)^2\\
& = 4 \E\left(\int_0^T\exp\{-4n(\gamma^n+\theta^n)s\}d[M^n,M^n]_s\right)  \le 4 \E([M^n,  M^n]_T) \\
& = 4 n^{-2}\E\left(N^n \left(n^2(\gamma^n+\theta^n) \int_0^T Y^{n}(s) ds \right)\right)  = 4 n^{-1} (n\gamma^n+ n\theta^n) \int_0^T \E(Y^{n}(s))ds \\
& \to 0, \textt{as $n\to\infty.$}
\ea\ee
Now from \eqref{mart1}, \eqref{con1}, and \eqref{con3}, for any $T\in [0,\infty),$
\be\label{con4}\ba
&\E\left(\sup_{0\le t\le T} \left|Y^n(t) - \left[\left(x_0+ \frac{|\alpha-\beta|}{2(\gamma+\theta)}\right) \exp\{2(\gamma+\theta)t\} -  \frac{|\alpha-\beta|}{2(\gamma+\theta)}\right]\right|\right) \\
& \le \E\left(\sup_{0\le t\le T} \left|Y^n(t) - \left[\left(\barx^n(0)+ \frac{|\alpha^n-\beta^n|}{2n(\gamma^n+\theta^n)}\right) \exp\{2n(\gamma^n+\theta^n)t\} - \frac{|\alpha^n-\beta^n|}{2n(\gamma^n+\theta^n)}\right] \right|\right)+ o(1) \\
& \le \exp\{2n(\gamma^n+\theta^n)T\} \E\left(\sup_{0\le t \le T}\int_0^t \exp\{-2n(\gamma^n+\theta^n)s\} d O^n(s)\right) \\
&\quad + \exp\{2n(\gamma^n+\theta^n)T\}\E\left( \sup_{0\le t \le T}\left|\int_0^t \exp\{-2n(\gamma^n+\theta^n)u\} d M^n(u)\right|\right) + o(1) \\
&\to 0,  \textt{as $n\to\infty.$}
\ea\ee

We next observe that
\bes\ba
&|\barx^n(t) - x(t)| \\
&\le |\barx^n(0) - x_0| + |\barn^n_s(t) -\alpha^n t - \barn^n_b(s)+\beta^n t| \\
& \quad + \left|  \barn_{sr}^n\left(n\theta^n \int_0^t
\barx^{n,+}(s)ds\right) - n\theta^n \int_0^t
\barx^{n,+}(s)ds \right|\\
&\quad+ \left|\barn_{br}^n\left(n\gamma^n\int_0^t
\barx^{n,-}(s)ds\right) - n\gamma^n\int_0^t
\barx^{n,-}(s)ds\right| \\
&\quad + \left|n\theta^n \int_0^t
\barx^{n,+}(s)ds - \theta \int_0^t
\barx^{n,+}(s)ds\right|  + \left|n\gamma^n \int_0^t
\barx^{n,-}(s)ds - \gamma \int_0^t
\barx^{n,-}(s)ds\right|\\
&\quad + \left|\theta \int_0^t
\barx^{n,+}(s)ds - \theta \int_0^t
x^{+}(s)ds\right|  + \left|\gamma \int_0^t
\barx^{n,-}(s)ds - \gamma \int_0^t
x^{-}(s)ds\right| \\
& \le |\barx^n(0) - x_0| + O^n(t) \\
& \quad  + \left|  \barn_{sr}^n\left(n\theta^n \int_0^t
\barx^{n,+}(s)ds\right) - n\theta^n \int_0^t
\barx^{n,+}(s)ds \right|\\
&\quad+ \left|\barn_{br}^n\left(n\gamma^n\int_0^t
\barx^{n,-}(s)ds\right) - n\gamma^n\int_0^t
\barx^{n,-}(s)ds\right| \\
&\quad + (\left|n\theta^n - \theta\right| + |n\gamma^n -\gamma|) \int_0^t
Y^{n}(s)ds \\
&\quad + (\theta + \gamma)\int_0^t \left|\barx^{n}(s) - x(s)\right|  ds.
\ea\ees
Gronwall's inequality yields that
\be\label{est1}\ba
\E\left(\sup_{0\le t\le T} |\barx^n(t) - x(t)|\right) &\le\E \left(|\barx^n(0) - x_0| + \sup_{0\le t\le T}O^n(t) \right. \\
& + \sup_{0\le t\le T}\left|  \barn_{sr}^n\left(n\theta^n \int_0^t
\barx^{n,+}(s)ds\right) - n\theta^n \int_0^t
\barx^{n,+}(s)ds \right|\\
&\quad+ \sup_{0\le t\le T}\left|\barn_{br}^n\left(n\gamma^n\int_0^t
\barx^{n,-}(s)ds\right) - n\gamma^n\int_0^t
\barx^{n,-}(s)ds\right| \\
&\quad + \left. (\left|n\theta^n - \theta\right| + |n\gamma^n -\gamma|) \int_0^T
Y^{n}(s)ds \right) e^{(\theta+\gamma)T}.
\ea\ee
Using the argument between \eqref{abandonment-bound} and \eqref{poisson-bound-1}, and from the results in \eqref{poisson-bound} and \eqref{poisson-bound-1}, we have for some $c_2, c_3\in (0, \infty),$
\be\label{con5}\ba
& \E\left(\sup_{0\le t\le T} \left|  \barn_{sr}^n\left(n\theta^n \int_0^t\barx^{n,+}(s)ds\right) - n\theta^n \int_0^t\barx^{n,+}(s)ds \right|\right) \\
&\le \frac{c_2(T+ \E(|\barx^n(0)|))}{\sqrt{n}} \to 0, \textt{as $n\to\infty,$}
\ea\ee
and
\be\label{con6}\ba
& \E\left(\sup_{0\le t\le T} \left|  \barn_{br}^n\left(n\theta^n \int_0^t\barx^{n,-}(s)ds\right) - n\gamma^n \int_0^t\barx^{n,-}(s)ds \right|\right) \\
& \le \frac{c_3(T+ \E(|\barx^n(0)|))}{\sqrt{n}} \to 0, \textt{as $n\to\infty.$}
\ea\ee
Applying \eqref{con1}, \eqref{con5}, \eqref{con6}, and the convergence $n\gamma^n \to \gamma, n\theta^n\to \theta$ to \eqref{est1},  \eqref{con} follows immediately. \ink

{\bf Proof of Lemma \ref{fluid-s}:} We first show (i) and (ii). Assume $\alpha \ge \beta.$ We consider the following three situations. \\
{\bf (a)} Let $x_0 > 0.$ Define $\tau_1 = \inf\{t\ge 0: x(t) \le 0\}.$ Then for $t \in [0, \tau_1), $ we have $x(t)\ge 0,$ and so
\be\label{eqn1:fluid-s}
x(t) = x_0 + (\alpha -\beta)t - \theta \int_0^t x(s)ds.
\ee
Solving the above equation, we have for $t\in[0,\tau_1),$
\be\label{eqn2:fluid-s}
x(t) = \left(x_0 - \frac{\alpha -\beta}{\theta}\right) e^{-\theta t} + \frac{\alpha - \beta}{\theta}.
\ee
If $\tau_1 < \infty,$ then $x(\tau_1) = \lim_{t \uparrow \tau_1} x(t) > 0,$ which contradicts the definition of $\tau_1$. Thus $\tau_1 = \infty,$
and so equation \eqref{eqn2:fluid-s} holds for all $t\in[0,\infty).$ \\
{\bf (b)} Let $x_0 = 0.$ We first assume $\alpha > \beta$ and note that
\bes
x'(0) = \alpha - \beta - \theta x^+_0 + \gamma x^-_0 = \alpha - \beta >0.
\ees
So there exists $\tau_2 > 0$ such that $x(t) > 0$ for $t\in(0, \tau_2].$ Define $\tilde x(t) = x(t+\tau_2), \ t\in[0,\infty).$ Then we have for $t\in[0,\infty),$
\bes
\tilde x(t) = \tilde x(0) + (\alpha -\beta) t + \int_0^t -\theta \tilde x^+(s) + \gamma \tilde x^-(s) ds.
\ees
Noting that $\tilde x(0) = x(\tau_2) >0$, and using the result in Part (a), we obtain that $\tilde x(t) > 0$ for all $t\in (0,\infty).$ Thus $x(t)\ge 0$ for all $t\in[0,\infty),$ and so equations \eqref{eqn1:fluid-s} and \eqref{eqn2:fluid-s} hold for all $t\in[0,\infty).$ If $\alpha = \beta,$ then
\be\label{zero}
x(t) = 0 \ \mbox{for all} \ t\in[0,\infty).
\ee
Otherwise, if \eqref{zero} fails, then there exists $0< t_1 < t_2 < \infty$ such that $x(t_1) = 0$ and $x(s) \neq 0$ for all $s\in (t_1, t_2].$ Without loss of generality, we assume $x(s)>0$ for $s\in(t_1, t_2].$ Then for $x\in (t_1, t_2],$
\bes
x(s) = x(t_1) - \theta \int_{t_1}^s x(u)du = - \theta \int_{t_1}^s x(u) du < 0,
\ees
which is a contradiction. \\
{\bf (c)} Let $x_0 < 0.$ We note that
\bes
x'(0) = \alpha - \beta - \theta x^+_0 + \gamma x^-_0 = \alpha - \beta - \gamma x_0>0.
\ees
Let $\tau_3 = \inf\{t\ge 0: x(t) \ge 0\}.$ Then for $t\in [0, \tau_3]$,
\bes
x(t) = x_0 + (\alpha - \beta)t - \gamma \int_0^t x(s)ds,
\ees
and so
\be\label{eqn3:fluid-s}
x(t) = \left(x_0 - \frac{\alpha -\beta}{\gamma}\right) e^{-\gamma t} + \frac{\alpha - \beta}{\gamma}.
\ee
From the fact that $x(\tau_3) =0,$ we have
\bes
\tau_3 = \gamma^{-1}\log\left(\frac{(\alpha -\beta) -\gamma x_0}{\alpha - \beta}\right) \in (0,\infty).
\ees
Define $\hat x(t) = x(t+\tau_3), t\in[0,\infty).$ We have for $t\in[0,\infty),$
\bes
\hat x(t) = \hat x(0) + (\alpha -\beta) t + \int_0^t \left(-\theta \hat x^+(s) + \gamma \hat x^-(s) \right) ds.
\ees
Noting that $\hat x(0) = x(\tau_3) = 0,$ and using the result in Part (b), we know that $\hat x(t) \ge 0$ for all $t\in[0,\infty).$ Hence $x(t) \ge 0$ for all $x\in[\tau_3, \infty),$ and equations \eqref{eqn1:fluid-s} and \eqref{eqn2:fluid-s} hold for $t\in[\tau_3, \infty).$ Combining this with \eqref{eqn3:fluid-s}, we obtain that
\bes
x(t) =
\begin{cases}
\left(x_0 - \frac{\alpha -\beta}{\gamma}\right) e^{-\gamma t} + \frac{\alpha - \beta}{\gamma}, & t\in[0, \tau_3], \\
 \frac{\alpha -\beta}{\theta}\left(1-e^{-\gamma t}\right), &  t\in[\tau_3,\infty).
 \end{cases}
\ees
At last, letting $y(t) = -x(t)$ and using the results in (i) and (ii), the results in (iii) and (iv) follow immediately.
\ink

{\bf Proof of Theorem \ref{diffusion-d}:} We first note that for $t\ge 0,$
\bes\ba
Z^n(t) & = \hatx^n(t) - \sqrt{n}x^n(t) \\
& = Z^n(0) + \hatn^n_s(t) - \hatn^n_b(t) - \hatn^n_{sr}\left(n\theta^n \int_0^t \barx^{n,+}(u) du \right) + \hatn^n_{br}\left(n\gamma^n \int_0^t \barx^{n,-}(u) du \right) \\
& \quad - n\theta^n \int_0^t \left(\hatx^{n,+}(u) - \sqrt{n}x^{n,+}(u) \right) du + n\gamma^n \int_0^t \left(\hatx^{n,-}(u) - \sqrt{n}x^{n,-}(u) \right) du. 
\ea\ees
Define for $t\geq 0$,
$$
\hatw^n(t) = \hatn^n_s(t) - \hatn^n_b(t) - \hatn^n_{sr} \left(n\theta^n\int_0^t\barx^{n,+}(u)du\right) +
\hatn^n_{br}\left(n\gamma^n\int_0^t\barx^{n,-}(u)du\right).
$$
We observe that, by the functional central limit
theorem for renewal processes (see Theorem 14.6 in Billingsley \cite{billingsley1999convergence}),
\be\label{fclt-renewal}
\hatn^n_s\Rightarrow W_s, \; \hatn_b^n\Rightarrow W_b,
\ee
where $W_s$ and $W_b$ are independent Brownian motions with zero
drifts and variances $\alpha^3\sigma^2$ and $\beta^3\varsigma^2$,
respectively. We also note that $\hatn^n_{sr}$ and $\hatn^n_{br}$
converge weakly to independent standard Brownian motions from the functional central
limit theorem for unit Poisson process. Further noting from Theorem
\ref{fluid-d} and Lemma \ref{fluid-s}, and using the random change of time theorem (see Section 3.14 in Billingsley \cite{billingsley1999convergence}), we obtain that
\be\label{fclt-renege}\ba
\hatn^n_{sr}\left(n\theta^n\int_0^\cdot \barx^{n,+}(u)du\right) &\Rightarrow  \int_0^\cdot \sqrt{\theta x^+(u)} d W_{sr}(u), \\
\hatn^n_{br}\left(n\gamma^n\int_0^\cdot \barx^{n,-}(u)du\right)& \Rightarrow \int_0^\cdot \sqrt{\gamma x^-(u)} dW_{br}(u), 
\ea\ee 
where $W_{sr}$ and $W_{br}$ are independent standard Brownian motions, which are independent of $W_s$ and $W_b$. 
Combining \eqref{fclt-renewal} and \eqref{fclt-renege}, we have
$$\hatw^n\Rightarrow \int_0^\cdot \sqrt{{\alpha ^3}{\sigma ^2} + {\beta ^3}{\varsigma ^2}  + \theta x^+(u) + \gamma x^-(u)} d W(u),$$ 
where $W$ is a standard Brownian motion. Let 
\[
M(t) = \int_0^t \sqrt{{\alpha ^3}{\sigma ^2} + {\beta ^3}{\varsigma ^2} + \theta x^+(u) + \gamma x^-(u)} d W(u), \ t\ge 0.
\]
There exists a random variable $Z(0)$ with law $\nu$ such that $(Z^n(0), \hatw^n) \Go \left(Z(0), M\right).$
By Skorohod representation theorem, without loss of generality, we assume that $(Z^n(0), \hatw^n)$ and $(Z(0), M)$ are defined on the same probability space and $(\hatx^n(0),
\hatw^n) \go (Z(0), M)$ almost surely and uniformly on compact sets of $[0,\infty)$. Define for $t\ge 0,$
$$\tilz^n(t) = Z(0) + M(t) -n\theta^n \int_0^t \tilz^{n,+}(s)ds + n\gamma^n \int_0^t \tilz^{n,-}(s)ds, $$
and
$$Z(t) = Z(0) + M(t) - \theta \int_0^t Z^{+}(s)ds + \gamma \int_0^t Z^{-}(s)ds.$$
From Lemma \ref{lip-mapping}, $\tilz^n$ and $Z$ are well-defined, and for $t\ge 0,$
\bes
\|\tilz^n-Z\|_t  \le (n\theta^n+n\gamma^n) \int_0^t \|\tilz^n-Z\|_s ds + |\theta - n\theta^n| \int_0^t Z^+(s) ds + |\gamma - n\gamma^n | \int_0^t Z^-(s) ds.
\ees
Using Gronwall's inequality,
\be\label{conv-zero-1}\ba
\|\tilz^n-Z\|_t & \le \left(|\theta - n\theta^n| \int_0^t Z^+(s) ds + |\gamma - n\gamma^n | \int_0^t Z^-(s) ds\right) e^{(n\theta^n+n\gamma^n)t} \\
& \to 0, \;\; \mbox{almost surely.} 
\ea\ee
Recall that for $t\ge 0,$
\bes\ba
Z^n(t)  = Z^n(0)  +\hatw^n(t) - n\theta^n\int_0^t\left(\hatx^{n,+}(s)-\sqrt{n}x^{n,+}(s)\right)ds + n\gamma^n\int_0^t\left(\hatx^{n,-}(s)-\sqrt{n}x^{n,-}(s)\right)ds.
\ea\ees
We then have that for $t\geq 0$, 
\bes\ba
\|Z^n -  \tilz^n\|_t & \le |Z^n(0)-Z(0)| + \|\hatw^n-M\|_t + n\theta^n \int_0^t |\hatx^{n,+}(s) - \sqrt{n}x^{n,+}(s) - \tilz^{n,+}(s)| ds \\
&\quad + n\gamma^n  \int_0^t |\hatx^{n,-}(s) - \sqrt{n}x^{n,-}(s) - \tilz^{n,-}(s)| ds \\ 
& \le |Z^n(0)-Z(0)| + \|\hatw^n-M\|_t + (n\theta^n+n\gamma^n)\int_0^t \|Z^n-\tilz^n\|_s ds. 
\ea\ees 
By Gronwall's inequality, 
\be\label{conv-zero-2}\ba
\|Z^n -  \tilz^n\|_t & \leq \left(|Z^n(0)-Z(0)| + \|\hatw^n-M\|_t
\right) e^{(n\theta^n+n\gamma^n)t} \\
& \to 0, \;\; \mbox{almost surely.}
\ea\ee
Combining \eqref{conv-zero-1} and \eqref{conv-zero-2}, the result follows immediately.  \ink

{\bf Proof of Theorem \ref{diff-ht}:} From \eqref{htc}, it is clear that $\alpha = \beta,$ where $\alpha = \lim_{n\to\infty} \alpha^n$ and $\beta = \lim_{n\to\infty} \beta^n.$ Then from Theorem \ref{fluid-d} and Lemma \ref{fluid-s}, $x\equiv 0.$ The rest of the proof is very similar to that of Theorem \ref{diffusion-d}. To show the convergence of $\hatx^n,$ we observe that for $t\ge 0,$
\bes
\hatx^n(t) = \hatx^n(0)  +\hatw^n(t) + \sqrt{n}(\alpha^n-\beta^n)t-
n\theta^n\int_0^t\hatx^{n,+}(s)ds +
n\gamma^n\int_0^t\hatx^{n,-}(s)ds, 
\ees
where 
$$\hatw^n(t) = \hatn^n_s(t) - \hatn^n_b(t) - \hatn^n_{sr} \left(n\theta^n\int_0^t\barx^{n,+}(u)du\right) +
\hatn^n_{br}\left(n\gamma^n\int_0^t\barx^{n,-}(u)du\right).$$
Same as the proof of Theorem \ref{diffusion-d}, we have that 
\bes
\hatn^n_s\Rightarrow W_s, \; \hatn_b^n\Rightarrow W_b,
\ees
where $W_s$ and $W_b$ are independent Brownian motions with zero
drifts and variances $\alpha^3\sigma^2$ and $\beta^3\varsigma^2$,
respectively. We also note that $\hatn^n_{sr}$ and $\hatn^n_{br}$
converge weakly to independent standard Brownian motions from functional central
limit theorem for unit Poisson process. Further noting that $\barx^n\Rightarrow 0$, and using the random change of time theorem, we have that
\bes\ba
\hatn^n_{sr}\left(n\theta^n\int_0^\cdot \barx^{n,+}(u)du\right) &\Rightarrow 0, \\
\hatn^n_{br}\left(n\gamma^n\int_0^\cdot \barx^{n,-}(u)du\right)& \Rightarrow 0. 
\ea\ees 
Combining the above convergences, we have $\hatw^n\Rightarrow \sqrt{{\alpha ^3}{\sigma ^2} + {\beta ^3}{\varsigma ^2}} W$, where $W$ is a standard Brownian motion. Furthermore,
there exists a random variable $Z(0)$ with law $\nu$ such that
$(\hatx^n(0), \hatw^n) \Go (Z(0), \sqrt{\alpha^3\sigma^2+\beta^3\varsigma^2} W)$.
By Skorohod representation theorem, without loss
of generality, we assume that $(\hatx^n(0), \hatw^n)$ and $(Z(0),
W)$ are defined on the same probability space and $(\hatx^n(0),
\hatw^n) \go (Z(0), W)$ almost surely and uniformly on compact sets of $[0,\infty)$. Define
$$\hatx(t) = \hatx(0) +  \sqrt{\alpha^3\sigma^2+\beta^3\varsigma^2} W(t) + ct -\theta \int_0^t \hatx^+(s)ds + \gamma \int_0^t \hatx^-(s)ds, \;\; t\geq 0.$$
From Lemma \ref{lip-mapping}, $\hatx$ is well-defined. 
We then have that for $t\geq 0$, 
\bes \ba
\|\hatx^n -  \hatx\|_t & \leq |\hatx^n(0)-\hatx(0)| +
\|\hatw^n-\sqrt{\alpha^3\sigma^2+\beta^3\varsigma^2}W\|_t + |\sqrt{n}(\alpha^n-\beta^n)-c|t\\
&\quad + (n\theta^n+n\gamma^n+\theta+\gamma)\int_0^t
\|\hatx^n-\hatx\|_s ds. 
\ea\ees 
By Gronwall's inequality, 
\bes\ba
\|\hatx^n -  \hatx\|_t & \leq \left(|\hatx^n(0)-\hatx(0)| + \|\hatw^n-\sqrt{\alpha^3\sigma^2+\beta^3\varsigma^2}W\|_t  + |\sqrt{n}(\alpha^n-\beta^n)-c|t
\right) e^{(n\theta^n+n\gamma^n+\theta+\gamma)t} \\
&  \to 0, \;\; \mbox{almost surely.}
\ea\ees 
The result follows immediately.  \ink

{\bf Proof of Theorem \ref{limitphi}:}
We first follow Section 5 of Chapter 15 in Karlin and Taylor \cite{karlin1981second} to construct a stationary density for $\hatx$. Denote by $\mu(x)$ the infinitesimal drift parameter $c - \theta {x^ + } + \gamma {x^ - }.$ We note that an antiderivative of $\frac{2\mu(x)}{a^2}$ is $ \frac{{2c}}{{{a^2}}}x - \frac{\theta }
{{{a^2}}}x^21\{ x \ge 0\}  - \frac{\gamma }{{{a^2}}}x^21\{ x < 0\}$. Define for $x\in\RR,$
\bes
s(x) = \exp \left\{ { \frac{{2c}}{{{a^2}}}x - \frac{\theta }
{{{a^2}}}x^21\{ x \ge 0\}  - \frac{\gamma }{{{a^2}}}x^21\{ x < 0\} }
\right\}.
\ees
We define a density function as follows:
\bes\ba
\psi(x)& =   \tilde C s(x) \\
& =  \begin{cases} \tilde C \exp \left\{  \frac{2c}{a^2}x - \frac{\theta } {a^2}x^2 \right\}, & x\ge 0\\
\tilde C \exp \left\{ \frac{2c}{a^2}x  - \frac{\gamma }{a^2}x^2 \right\}, & x<0
\end{cases} \\
& = \begin{cases}\frac{C}{\sqrt{\theta}}\exp \left\{ {\frac{{{c^2}}}{{\theta {a^2}}}}
\right\} \phi \left( x; \frac{c}{\theta}, \frac{a^2}{2\theta}\right), & x\ge 0\\
\frac{C}{\sqrt{\gamma}}\exp \left\{ {\frac{{{c^2}}}{{\gamma {a^2}}}}
\right\}\phi \left(x; \frac{c}{\gamma}, \frac{a^2}{2\gamma}\right), & x<0,
\end{cases}
\ea\ees
where
\bes
\tilde C = \frac{1}{\int_{-\infty}^\infty s(x) ds},
\ees
and
\bes
C = a\sqrt{\pi} \tilde C = \frac{1}{{ \frac{1}{\sqrt{\theta}} \exp \left\{
{\frac{{{c^2}}}{{\theta {a^2}}}} \right\}\left( {1 - \Phi \left(
{0;\frac{c}{\theta },\frac{{{a^2}}}{{2\theta }}} \right)} \right) +
\frac{1}{\sqrt{\gamma}}\exp \left\{ {\frac{{{c^2}}}{{\gamma {a^2}}}}
\right\} { \Phi \left( {0;\frac{c}{\gamma },\frac{{{a^2}}}{{2\gamma
}}} \right)} }}.
\ees
The uniqueness of the stationary distribution follows from the irreducibility, i.e. $a>0$ (see Peszat and Zabczyk \cite{peszat95}). 

We now study the limiting distribution for $Z.$ Recall $b=\sqrt{a^2 + |\alpha-\beta|}$. Define a time-homogeneous stochastic process $Z^*$ as follows.
\[
Z^*(t) = Z(0) + b W(t) -\theta\int_0^t Z^{*,+}(s) ds + \gamma \int_0^t Z^{*,-}(s) ds.
\]
We are going to show that the limiting distribution of $Z$ is the same as the stationary distribution of $Z^*.$ We first note the unique stationary distribution of $Z^*$ is given by $\psi(x, 0, 0, b).$ We then consider $Z-Z^*.$ Let $b(t) = \sqrt{{\alpha ^3}{\sigma ^2} + {\beta ^3}{\varsigma ^2}  + \theta x^+(t) + \gamma x^-(t)}, t\ge 0.$ We have that for $t\ge 0,$
\[
Z(t) - Z^*(t) = \int_0^t (b(s)-b) d W(s) + \int_0^t \left[-\theta (Z^+(s) - Z^{*,+}(s)) + \gamma(Z^-(s)-Z^{*,-}(s)) \right] ds.
\]
For $x, y \in \RR, $ let $g(x,y) = -\theta(x^+ - y^+) + \gamma (x^- - y^-)$, and we note that
\bes
g(x,y) = 
\begin{cases}
-\theta(x - y), & \textt{if} \ x>0, y>0,\\
-\theta x + \gamma y, & \textt{if} \ x>0, y\le 0, \\
-\gamma x +\theta y, & \textt{if} \ x\le 0, y>0,\\
-\gamma(x-y), & \textt{if} \ x\le 0, y\le 0.
\end{cases}
\ees
 In the following, we assume $\theta \le \gamma.$ (The case when $\theta > \gamma$ can be treated analogously.) Define
 \bes
h_1(z) = 
\begin{cases}
-\theta z, & \textt{if} \ z\ge0, \\
-\gamma z, & \textt{if} \ z < 0, 
\end{cases}
\ees
and
 \bes
h_2(z) = 
\begin{cases}
-\gamma z, & \textt{if} \ z\ge0, \\
-\theta z, & \textt{if} \ z < 0, 
\end{cases}
\ees
Then for all $x,y\in \RR,$
\be\label{ineqn}
h_2(x-y) \le g(x,y)\le h_1(x-y).
\ee 
Now consider the following stochastic integral equations.
\bes
V_1(t) = \int_0^t (b(s)-b) dW(s) + \int_0^t h_1(V_1(s))ds,
\ees
and
\bes
V_2(t) = \int_0^t (b(s)-b) dW(s) + \int_0^t h_2(V_2(s))ds,
\ees
Then from \cite{chanwilliams89}, $V_1(t) \to 0$ and $V_2(t)\to 0$ a.s. as $t\to\infty.$ Finally, from \eqref{ineqn}, 
\be\label{upper-lower-bound}
V_2(t) \le Z(t)-Z^*(t) \le V_1(t), \ t\ge 0.
\ee  
In fact, if for some $(t_1, t_2)$ such that $Z(t_1)-Z^*(t_1) = V_2(t_1)$ and $Z(t)-Z^*(t) < V_2(t)$ for $t\in(t_1, t_2).$ Then
\bes\ba
0 >  Z(t)-Z^*(t) - V_2(t) & = \int_0^t g(Z(s),Z^*(s)) ds - \int_0^t h_2(V_2(s)) ds \\
& \ge \int_0^t \left[ h_2(Z(s)-Z^*(s)) - h_2(V_2(s)) \right] ds  > 0,
\ea\ees
which is a contradiction. This shows $Z(t)-Z^*(t) \ge V_2(t)$ for $t\ge 0.$ Similarly, we can show $Z(t) -Z^*(t) \le V_1(t)$ for $t\ge 0.$ From \eqref{upper-lower-bound}, we have $Z(t)-Z^*(t) \to 0$ as $t\to\infty.$ The result in (i) follows. \ink

\subsection{Proof of Theorem \ref{limit-change}}

We will apply the following multiplicative Foster's criterion to show the positive recurrence of $\hatmbx^n$ for each $n\in \NN.$ Such criterion is introduced in Chapter 4 of \cite{bramson08}.

\begin{proposition}[Multiplicative Foster's Criterion]\label{foster} Suppose that $\{Y(t):t\ge0\}$ is a continuous time Markov process, such that for some positive $c, \epsilon,$ and $\kappa,$
\be\label{foster-1}
\E_y\left(|Y(c(|y|\vee\kappa))| \right) \le (1-\epsilon) (|y|\vee \kappa), \ \ \mbox{for all $y$.}
\ee
If 
\be\label{foster-2}
\{y: |y|\le \kappa\} \ \mbox{is a closed petite set,}
\ee
then $Y$ is positive Harris recurrent. 

\end{proposition}

Denote by $\mathbb{S}_n$ the state space of $\hatmbx^n$, i.e., $\mathbb{S}_n = (n^{-1/2} \mathbb{Z}) \times \RR_+ \times \RR_+,$ where $n^{-1/2} \mathbb{Z} = \{n^{-1/2}x: x\in\mathbb{Z}\}.$ For $z=(x,y_1,y_2) \in \hatmbx^n,$ define $|z| = \sqrt{x^2+y_1^2+y_2^2}.$ Let $ \mathcal{B}(\mathbb{S}_n)$ denote the Borel $\sigma$-field on $\mathbb{S}_n$ (the countable set $n^{-1/2} \mathbb{Z}$ is endowed with the discrete metric). A nonempty set $A \in \mathcal{B}(\mathbb{S}_n)$ is said to be petite if for some probability measure $a$ on $(0,\infty)$ and some nontrivial measure $\nu$, 
\bes
\nu(B) \le \int_0^\infty \PP_x\left(\hatmbx^n(t)\in B\right) a(dt)
\ees  
for all $x\in A$ and $B\in \mathcal{B}(\mathbb{S}_n).$ Conditions \eqref{foster-1} and \eqref{foster-2} are shown to be satisfied in Lemmas \ref{growth-bound} and \ref{petite}, respectively.

\begin{lemma}\label{growth-bound}
There exists $c_1 \in(0,\infty)$ such that for all $r\geq 0$ and $z = (x, y_1, y_2) \in \mathbb{S}_n$,
\be\label{stability-1}
\sup_{n\geq 1}\E_z\left(\left|\hatmbx^n(r|z|)\right|^2 \right) \le c_1(1+(r+1)|z|),
\ee
and moreover, 
\be\label{stability-2}
\lim_{|z|\to\infty}  \frac{\sup_{n\ge 1}\E_z\left(\left|\hatmbx^n(r|z|)\right|^2 \right) }{|z|^2}  = 0.
\ee
\end{lemma}
{\bf Proof:} We first consider $\hatx^n$ with initial value $\hatx^n(0)=x$. Recall from the proof of Theorem \ref{diff-ht} that for $t\ge 0,$
\be\label{queue-stability}
\hatx^n(t) = x  +\hatw^n(t) + \sqrt{n}(\alpha^n-\beta^n)t-
n\theta^n\int_0^t\hatx^{n,+}(s)ds +
n\gamma^n\int_0^t\hatx^{n,-}(s)ds, 
\ee
where 
$$\hatw^n(t) = \hatn^n_s(t) - \hatn^n_b(t) - \hatn^n_{sr} \left(n\theta^n\int_0^t\barx^{n,+}(u)du\right) +
\hatn^n_{br}\left(n\gamma^n\int_0^t\barx^{n,-}(u)du\right).$$
From Lemma 3.5 in \cite{budhirajaghosh06}, we have for some $c_1 \in (0,\infty)$ (independent of $n$ and $t$),
\be\label{renewal-bound}\ba
\E\left(\sup_{0\le u\le t} \left|\hatn^n_s(u) \right|^2 + \sup_{0\le u\le t} \left|\hatn^n_b(u) \right|^2\right) \le c_1(t+1).
\ea\ee
Define for $l\ge 0,$
\bes
\mathcal{G}^n_l = \sigma\{X^n(0), N^n_s(nv), N^n_b(nv), N^n_{sr}(nv), N^n_{br}(nv): v\in [0,l]\}.
\ees
Then $\hatn^n_{sr}$ is a $\{\mathcal{G}^n_l\}_{l\ge 0}$ square integrable martingale, and for any $t\ge 0,$ $T^n_1(t)\equiv n\theta^n\int_0^t \barx^{n,+}(v)dv$ is a $\{\mathcal{G}^n_l\}_{l\ge 0}$ stopping time. Using Burkholder-Davis-Gundy inequality (see Theorem 74 of Chapter IV in \cite{protter08}), for some $c_2\in (0,\infty)$ (independent of $n$ and $t$), 
\be\label{abandonment-bound}\ba
 \E\left(\sup_{0\le u \le t}\left|\hatn^n_{sr}\left(T^n_1(u)\right)\right|^2\right)  &  \le c_2 \E\left([\hatn^n_{sr},\hatn^n_{sr}](T^n_1(t))\right) \\
&  = c_2 \E(T^n_1(t)) = c_2n\theta^n \int_0^t \E( \barx^{n,+}(v)) dv.
\ea\ee
We next observe that for $t\ge 0,$ 
\bes
X^{n,+}(t) \le X^{n,+}(0) + N^n_s(t) - N^n_{sr}\left(\theta^n \int_0^t X^{n,+}(v) dv \right),
\ees
and so
\bes
\E(\barx^{n,+}(t)) \le \E(\barx^{n,+}(0)) + \E(\barn^n_s(t)) - n\theta^n \int_0^t \E(\barx^{n,+}(v)) dv.
\ees
From \eqref{renewal-bound}, there exists $c_3 \in (0,\infty)$ such that for $t\ge0,$
\[
\E(\barn^n_s(t)) = \frac{1}{\sqrt{n}}\E(|\hatn^n_s(t)|) + \alpha^n t \le c_3(t+1).
\]
Define for $t\ge 0$,
\bes
y^n_1(t) = \frac{x}{\sqrt{n}} + c_3(t+1) - n\theta^n \int_0^t y^n_1(v) dv.
\ees
Using the property of ordinary differential inequalities, we have 
\[
\E_z(\barx^{n,+}(t)) \le y^n_1(t), \ \ t\ge 0.
\]
Solving the ODE for $y^n_1$, we have for $t\ge 0,$
\be\label{fluid-ht-bound}
 \E_z(\barx^{n,+}(t)) \le y^n_1(t) = \left(\frac{x}{\sqrt{n}} + c_3 - \frac{c_3}{n\theta^n} \right)e^{-n\theta^n t} + \frac{c_3}{n\theta^n}.
\ee
Applying \eqref{fluid-ht-bound} to \eqref{abandonment-bound}, we have for some $c_4\in (0,\infty)$ (independent of $n$ and $t$),
\be\label{poisson-bound}
\E_z\left(\sup_{0\le u \le t}\left|\hatn^n_{sr}\left(\int_0^u \barx^{n,+}(v) dv \right)\right|^2\right) \le c_4 (t+|z|), \ t\ge 0. 
\ee
Using the similar argument, for some $c_5\in (0,\infty)$ (independent of $n$ and $t$),
\be\label{poisson-bound-1}
\E_z\left(\sup_{0\le u \le t}\left|\hatn^n_{br}\left(\int_0^u \barx^{n,-}(v) dv \right)\right|^2\right) \le c_5 (t+|z|), \ t\ge 0. 
\ee
From \eqref{renewal-bound}, \eqref{poisson-bound}, and \eqref{poisson-bound-1}, we have  
\bes
\E_z\left[\left(\sup_{0\le u\le t}\left|\hatw^n(u)\right|\right)^2\right] = \E_z\left(\sup_{0\le u\le t}\left|\hatw^n(u)\right|^2\right) \le  (c_1+c_4+c_5) (t+|z|+1), \ t\ge 0.
\ees
Define for $t\ge 0,$
\be\label{center}
\tilde x^n(t) = x + \sqrt{n}(\alpha^n-\beta^n)t- n\theta^n\int_0^t\tilde x^{n,+}(s)ds + n\gamma^n\int_0^t\tilde x^{n,-}(s)ds.
\ee
Let $\phi(x) = -n\theta^n x^+ + n\gamma^n x^-, x\in \RR$ in Lemma \ref{lip-mapping}. Noting that $\phi$ is Lipschitz continuous with Lipschitz constant $\sup_{n\in\NN} \max\{n\theta^n,  n\gamma^n\},$ the Lipschitz constant for the mapping $\mathcal{M}^\phi$ in independent of $t$ and $n$. More precisely, there exists $\kappa\in(0,\infty)$ such that for $t\ge 0$,
\bes
\|\mathcal{M}^\phi(x_1) - \mathcal{M}^\phi(x_2) \|_t \le \kappa \|x_1 - x_2 \|_t. 
\ees 
Thus we have for $t\ge 0,$
\be\label{queue-bound}\ba
\E_z\left(\sup_{0\le u\le t} |\hatx^n(u) - \tilde x(u)|^2\right) & \le \kappa  \E_z\left(\sup_{0\le u\le t}\left|\hatw^n(u)\right|^2\right) \\
&  \le  \kappa {(c_1+c_4+c_5) (t+|z|+1)}.
\ea\ee
Solving \eqref{center} (similar to Lemma \ref{fluid-s}), we have 
\be\label{center-estimate}
\sup_{0\le u\le t} |\tilde x(u)| \le |z| e^{-\min\{n\theta^n, n\gamma^n\} t} +  \frac{\sqrt{n}|\alpha^n -\beta^n|}{\min\{n\theta^n, n\gamma^n\}}. 
\ee
Let $L_1 = \inf_{n\in\NN}\min\{n\theta^n, n\gamma^n\}$ and $L_2 = \sup_{n\in\NN} \frac{\sqrt{n}|\alpha^n -\beta^n|}{\min\{n\theta^n, n\gamma^n\}}.$ 
Combining \eqref{queue-bound} and \eqref{center-estimate}, we have
\bes
\E_z\left(\sup_{0\le u\le t} |\hatx^n(u)|^2\right) \le  \kappa {(c_1+c_4+c_5) (t+|z|+1)} + 2|z|^2 e^{-2L_1 t} + 2L_2^2,
\ees 
and so
\be\label{queue-bound-new}
\E_z\left( |\hatx^n(r|z|)|^2\right) \le \kappa {(c_1+c_4+c_5) ((r+1)|z|+1)} + 2|z|^2 e^{-2L_1 r|z|} + 2L_2^2.
\ee
We next focus on $\hata^n_s$ and $\hata^n_b.$ For $t\ge 0,$ 
\bes
\E_z\left[ (\hata^n_s(t))^2\right] \le \frac{1}{n} \E_z\left[ \left(U^n_{N^n_s(nt)+1)}\right)^2\right] \le \frac{1}{n} \E_z \left(\sum_{k=1}^{N^n_s(nt)+1}(U^n_k)^2\right), 
\ees
and from Wald's identity, there exists $c_6\in (0,\infty)$ such that for $r\ge 0,$
\be\label{residual-1}\ba
\sup_{n\ge 1}\E_z\left[(\hata^n_s(r|z|))^2 \right] &\le  \sup_{n\ge 1} \left(\frac{1}{n} \E_z\left[(U^n_1)^2\right] \left[1+ \E_z(N^n_s(nr|z|))\right] \right)\\
& \le c_6 (r|z|+1).
\ea\ee
Similarly, there exists $c_7\in (0,\infty)$ such that for $r\ge 0,$
\be\label{residual-2}
\sup_{n\ge 1}\E_z\left[ (\hata^n_b(r|z|))^2 \right] \le   c_7 (r|z|+1).
\ee
It is clear that \eqref{stability-1} and \eqref{stability-2} follow, on combining \eqref{queue-bound-new}, \eqref{residual-1}, and \eqref{residual-2}. \ink

\begin{lemma}\label{petite}
Fix $n\in\NN.$ Assume that $\alpha^n \le \beta^n$, and that there exist $ 0 < t_b^n < t_s^n < \infty,$ such that for any $\kappa > 0, $ 
\be\label{assum-petite}
\PP(U_1^n \in (t_s^n -\kappa, t_s^n + \kappa)) > 0, \  \mbox{and} \ \  \PP(V_1^n \in (t_b^n -\kappa, t_b^n + \kappa)) > 0. 
\ee
Then for $R\in (0,\infty),$ the set $B_n = \{z\in \mathbb{S}_n: |z| \le R\}$ is petite. 
\end{lemma}
{\bf Proof:} The proof idea is similar to those of Propositions 3.7 and 3.8 in Chapter 4 of \cite{bramson08} and Lemma 3.7 in \cite{MeynDown94}. For notation convenience, we drop $n$ from all quantities. 
Without loss of generality, assume $X^{+}(0)>0$. Choose $K\in\NN$ and $\delta\in (0,\infty)$ such that 
\[
Kt_b < t_s \le (K+1)t_b,
\]
and
\[
t_b >3\delta, \ \ K(t_b + \delta) < t_s - \delta.
\]
Denote by $G_s(t)$ and $G_b(t)$ the number of sellers and buyers abandoning the system by time $t$, and let $P_s$ and $P_b$ be exponential random variables with means $1/\theta$ and $1/\gamma$, respectively. Define the following events
\bes\ba
& E_1 = \left\{Kt_b -\frac{\delta}{2} \le \sum_{k=1}^K V_k \le K(t_b+\delta) \right\}, \\
& E_2 =  \left\{t_b - \frac{\delta}{2} < V_{K+1} \le t_b+\delta \right\}, \\
& E_3 =  \left\{t_s - \delta < U_1 \le t_s+\delta \right\},\\
& E_4 = \left\{ G_b\left( \frac{K(t_b+\delta)+(t_s-\delta)}{2} \right) = (K-X^{+}(0))^+\right\}, \\
& E_5 = \left\{ G_s\left( \frac{K(t_b+\delta)+(t_s-\delta)}{2} \right) = (K-X^{+}(0))^-\right\}.
\ea\ees
Then for $|z|\le R,$
\bes\ba
&\PP_z(E_1) \ge \left[\PP_z\left(t_b - \frac{\delta}{2K} \le V_1 \le t_b+\delta\right)\right]^K > 0, \\
& \PP_z(E_2) = \PP_z\left(t_b - \frac{\delta}{2} < V_{1} \le t_b+\delta \right) > 0, \\
& \PP_z(E_3) = \PP_z\left(t_s - \delta < U_{1} \le t_s+\delta \right) > 0,\\
&\PP_z(E_4|E_1\cap E_2 \cap E_3) \ge \left(\PP_z\left(P_b \le \frac{(t_s-\delta)-K(t_b+\delta)}{2}\right)\right)^K > 0, \\
& \PP_z(E_5|E_1\cap E_2 \cap E_3) \ge \left(\PP_z\left(P_s \le \frac{(t_s-\delta)+K(t_b+\delta)}{2}\right)\right)^R > 0.
\ea\ees
Noting that $E_1, E_2$ and $E_3$ are independent, and $E_4$ and $E_5$ are independent, there exists $\epsilon\in(0,1)$ such that
\bes
\PP(E_1\cap E_2\cap E_3 \cap E_4 \cap E_5) = \PP( E_4 \cap E_5|E_1\cap E_2\cap E_3 ) \PP(E_1\cap E_2\cap E_3)  \ge \epsilon. 
\ees
We note that for $t\in [\frac{(t_s-\delta)+K(t_b+\delta)}{2}, t_s-\delta],$
\bes
 \PP\left(\left. X(t) = 0 \right |  E_1\cap E_2\cap E_3 \cap E_4 \cap E_5\right) = 1.
\ees
For $B_1, B_2 \in \mathcal{B}(\RR_+), $ we have for $t\in [\frac{(t_s-\delta)+K(t_b+\delta)}{2}, t_s-\delta],$
\bes\ba
& \PP(X(t)=0, A_s(t) \in B_1, A_b(t)\in B_2) \\
& \ge \PP(A_s(t) \in B_1, A_b(t)\in B_2, E_1\cap E_2\cap E_3 \cap E_4 \cap E_5) \\
& = \PP\left(U_1 - t \in B_1, \sum_{k=1}^{K+1}V_k - t \in B_2, E_1\cap E_2\cap E_3 \cap E_4 \cap E_5 \right)  \\
& \ge \epsilon \PP\left(U_1 \in (B_1 + t) \cap(t_s-\delta, t_s+\delta)\right) \\
& \quad \quad \quad \times \PP\left( \sum_{k=1}^{K+1}V_k \in B_2+t, \sum_{k=1}^K V_k \in \left(Kt_b-\frac{\delta}{2}, K(t_b+\delta)\right), V_{K+1} \in \left(t_b-\frac{\delta}{2}, t_b +\delta \right)\right),
\ea\ees
where $B_1+t = \{y+t: y\in B_1\}$ and $B_2+t = \{y+t: y\in B_2\}$.
For $t\in [\frac{(t_s-\delta)+K(t_b+\delta)}{2}, t_s-\delta],$ let $a$ be a probability measure concentrated on $t$, and let $\nu$ be a measure on $\mathcal{B}(\mathbb{Z}\times\RR_+\times\RR_+)$ such that for any $B_0 \subset \mathbb{Z}, B_1, B_2\in \mathcal{B}(\RR_+),$
\bes\ba
\nu(B_0\times B_1\times B_2) & = \epsilon 1_{\{0\in B_0\}} \PP(U_1 \in (B_1+t)\cap (t_s-\delta, t_s+\delta)) \\
&\times \PP\left(\sum_{k=1}^{K+1}V_k \in B_2+t, \sum_{k=1}^K V_k \in (Kt_b-\delta, K(t_b+\delta)), V_{K+1} \in (t_b-\delta, t_b +\delta)\right).
\ea\ees
Clearly, $\nu(\mathbb{Z}\times\RR_+\times\RR_+) > 0,$ and so $\nu$ is nontrivial. 
Finally, for $B \in \mathcal{B}(\mathbb{Z}\times\RR_+\times\RR_+),$ we have 
\bes
\int_0^\infty \PP_z(\mathbb{X}^n(t)\in B) a(dt) \ge \nu(B).
\ees
This shows the lemma. \ink

Recall that $\Pi^n$ is a stationary distribution of $\hatmbx^n,$ and $\pi^n$ is the first-coordinate marginal distribution of $\Pi^n$. The following two lemmas will be used to show the tightness of $\pi^n.$ The proofs are the same as those of Theorems 3.4 and 3.5 in \cite{BudhirajaLee09}, and so we omit them here. For $\rho\in(0,\infty)$ and a compact set $C\subset \RR$, let
\be\label{stopping3}
\tau^n_C(\rho) = \inf\{t\geq \rho: \hatmbx^n(t)\in C\}.
\ee

\begin{lemma}\label{tightness-1}
For some $c_2, \rho \in (0,\infty)$ and a compact set $C\subset \mathbb{S}_n$, 
\bes
\sup_{n\in\NN} \E_z\left( \int_0^{\tau^n_C(\rho)} (1+ |\hatmbx^n(t)|) dt \right) \le c_2(1+|z|^2), \ z\in \mathbb{S}_n.
\ees
\end{lemma}

\begin{lemma}\label{tightness-2}
Let $f: \mathbb{S}_n \rightarrow \RR_+$ be a measurable map. Define for $\rho\in(0,\infty)$ and a compact set $C\subset \mathbb{S}_n$, 
\bes
G^n(z)= \E_z\left(\int_0^{\tau^n_C(\rho)} f(\hatmbx^n(t))dt\right), \ z\in \mathbb{S}_n.
\ees
Assume 
\be\label{cond}
\sup_{n\ge 1}  G^n(z) \;\mbox{is finite for all $z\in\mathbb{S}_n$, and uniformly bounded on $C$}.
\ee
Then there exists a $\kappa\in(0,\infty)$ such that, for all $n\in\NN, t\in[\rho,\infty)$ and $z\in \mathbb{S}_n$,
\bes
\frac{1}{t}\E_z\left[G^n(\hatmbx^n(t))\right]+ \frac{1}{t}\int_0^t \E_z\left[f(\hatmbx^n(s))\right]ds \leq \frac{1}{t}G^n(z)+ \kappa.
\ees
\end{lemma}

{\bf Proof of Theorem \ref{limit-change}:} We first show the positive recurrence of $\hatmbx^n$ for each $n\in\NN.$ We first note that condition \eqref{foster-1} follows from Lemma \ref{growth-bound}. Next without loss of generality, assume $\alpha^n \le \beta^n.$ We consider the following three cases.  
\bi
\item[\rm (1)] Assume that $\alpha^n \le \beta^n,$ and that one of $U_1^n$ and $V_1^n$ can take at least two positive values. Denote by $F_s^n$ and $F_b^n$ the distribution functions of $U^n_1$ and $V^n_1.$ For $\eta_1, \eta_2 \in [0,1],$ define
\bes
\tilde{t}_s^n = \sup\{t\ge 0: F_s^n(t) < 1-\eta_1\}, \ \tilde{t}_b^n = \inf\{t\ge 0: F_b^n(t) \ge \eta_2\}.
\ees
Noting that $\E(U_1^n) \ge \E(V_1^n),$ we can choose $\eta_1, \eta_2 \in [0,1]$ such that $\tilde{t}_s^n > \tilde{t}_b^n > 0$ and set 
\bes
t_s^n = \tilde{t}_s^n, \ \ t_b^n = \tilde{t}_b^n. 
\ees
From Lemma \ref{petite}, for $R\in (0,\infty),$ the set $B_n = \{z\in \mathbb{S}_n: |z| \le R\}$ is petite. Then the positive recurrence of $\hatmbx^n$ follows immediately from Proposition \ref{foster}.

%Then the positive recurrence of $\hatmbx^n$ follows immediately from Proposition \ref{foster} and Lemmas \ref{growth-bound} and \ref{petite}. 

\item[\rm (2)] Assume that $\alpha^n < \beta^n$, and that $\PP(U_1^n = 1/\alpha^n) = \PP(V_1^n = 1/\beta^n) =1.$ We can set 
\bes
t_s^n =1/\alpha^n, \ \ t_b^n = 1/\beta^n. 
\ees
Again from Lemma \ref{petite}, for $R\in (0,\infty),$ the set $B_n = \{z\in \mathbb{S}_n: |z| \le R\}$ is petite. Then the positive recurrence of $\hatmbx^n$ follows immediately from Proposition \ref{foster}.

%Then again the positive recurrence of $\hatmbx^n$ follows immediately from Proposition \ref{foster} and Lemmas \ref{growth-bound} and \ref{petite}. 
\item[\rm (3)] Assume that $\alpha^n = \beta^n$, and that $\PP(U_1^n = 1/ \alpha^n) = \PP(V_1^n = 1/\beta^n) =1.$ Then $N^n_s(t) = N^n_b(t),$ and $\hatmbx^n$ is a positive recurrent birth and death process.
\ei
Finally, the convergence of $\pi^n$ can be shown in the same way as those of Theorems 3.2 and 3.1 in \cite{BudhirajaLee09}, given the above Proposition \ref{growth-bound} and Lemmas \ref{tightness-1} and \ref{tightness-2}. \ink

\section{Numerical examples: Diffusion models, tables and figures}

\subsection{Two diffusion models} \label{diff-model}

We apply Theorems \ref{diffusion-d} and \ref{diff-ht} to derive two diffusion models for a double-ended queue with general parameters $\alpha, \beta, \sigma^2, \varsigma^2, \theta, \gamma$. 

{\bf Model I.} Consider a sequence of double-ended queues, indexed by $n\in\NN$ under Assumption \ref{main-d} and the heavy traffic condition \eqref{htc}. From Theorem \ref{diff-ht}, we have for large $N\in\NN$, 
\bes
\hatx^N \overset{d}{\approx} \hatx,
\ees
where as in \eqref{diff-ht-eqn},
\begin{equation*}
\hatx(t) =  \hatx(0)  + \sqrt{\alpha^3\sigma^2+\beta^3\varsigma^2}  W(t) + ct - \theta \int_0^t \hatx^{+}(u) du  + \gamma \int_0^t \hatx^{-}(u) du, \ t\ge0.
\end{equation*}
Fix such $N\in\NN$. Letting $s=Nt$, we have that
\bes
X^N \overset{d}{\approx} \hatx_1^N,
\ees
where 
\begin{equation*}
\hatx_1^N(s) = {X^N}(0) + \sqrt{\alpha^3\sigma^2+\beta^3\varsigma^2}  W(s) + \frac{c s}{\sqrt N } - \frac{\theta}{N}\int_0^s \hatx_1^{N,+}(u)du  + \frac{\gamma }{N}\int_0^s \hatx_1^{N,- }(u)du, \ s\ge0.
\end{equation*}
From Assumptions \ref{main-d} and the heavy traffic condition \eqref{htc}, we have that 
\bes
\sqrt{\alpha^3\sigma^2+\beta^3\varsigma^2} \approx \sqrt{(\alpha^N)^3(\sigma^N)^2+(\beta^N)^3(\varsigma^N)^2}, \  \frac{c}{\sqrt N } \approx \alpha^N-\beta^N, \ \frac{\theta}{N} \approx \theta^N, \ \frac{\gamma}{N} \approx \gamma^N,
\ees
and so
\bes
X^N \overset{d}{\approx} \hatx_1^N \overset{d}{\approx} \hatx_2^N,
\ees
where 
\begin{equation*}\ba
\hatx_2^N (s) & = {X^N}(0) + \sqrt {{{\left( {{\alpha ^N}} \right)}^3}
{{\left( {{\sigma ^N}} \right)}^2} + {{\left( {{\beta ^N}}
\right)}^3}{{\left( {{\varsigma ^N}} \right)}^2}} W(s) + \left(
{{\alpha ^N} - {\beta ^N}} \right)s \\
& \quad - {\theta ^N}\int_0^s \hatx_2^{N,+}(v)dv  + {\gamma ^N}\int_0^s \hatx_2^{N,-}(v)dv, \ s\ge 0.
\ea\end{equation*}
Thus for a double-ended queue with parameters satisfying Assumption \ref{main-d} and the heavy traffic condition \eqref{htc}, i.e. the arrival rates $\alpha, \beta$ are close, and the reneging rates $\theta, \gamma$ are very small comparing with $\alpha, \beta$ and $|\alpha-\beta|$, the dynamics of the queue length process $\{X(t): t\geq0\}$ can be approximated by an asymmetric O-U process
\begin{equation}\label{approx-model-1}
X_1(t) = X(0) + \sqrt{\alpha^3\sigma^2+\beta^3\varsigma^2}  W(t) + \left( {\alpha  - \beta }
\right)t - \theta \int_0^t X_1^+(u)du  + \gamma \int_0^t X_1^-(u)du, \ t\ge 0.
\end{equation}
As our first diffusion model, we use $X_1$ to approximate the queue length process with general parameters. 
%%More precisely, for a double-ended queue with parameters $\alpha, \beta, \sigma^2, \varsigma^2, \theta$, and $\gamma,$ the dynamics of the queue length process $\{X(t): t\geq0\}$ will be approximated by the asymmetric O-U process $\{X_1(t): t\ge 0\}$ defined in \eqref{approx-model-1}. 

{\bf Model II.} The second diffusion model can be obtained in the similar way from Theorems \ref{fluid-d} and \ref{diffusion-d}. To make it precise, consider a sequence of double-ended queues, indexed by $n\in\NN$ under Assumption \ref{main-d}. From Theorem \ref{diffusion-d}, for large enough $N \in \NN,$ we have 
\bes
Z^N = \hatx^N - \sqrt{N} x^N \overset{d}{\approx} Z,
\ees
where as in Theorem \ref{diffusion-d},
\bes\ba
Z(t)& = Z(0) + \int_0^{t} \sqrt{{\alpha ^3}{\sigma ^2} + {\beta ^3}{\varsigma ^2} + \theta x^+(u) + \gamma x^-(u) } dW(u) - \theta \int_0^t Z^+(u)du + \gamma \int_0^t Z^-(u) du, \ t\ge 0. 
\ea\ees
Fix such $N$ and let $s=Nt.$ Then we have
\bes
X^N(s) - \sqrt{N} x^N(s/N) \overset{d}{\approx} Z_1^N(s), \ s\ge 0,
\ees
where 
\bes\ba
Z_1^N(s)  & =  \sqrt{N} Z^N(0)+ \int_0^{s} \sqrt{{\alpha ^3}{\sigma ^2} + {\beta ^3}{\varsigma ^2} + \theta x^+(u/N) + \gamma x^-(u/N) } dW(u) \\
& \quad - \frac{\theta}{N} \int_0^s Z_1^{N,+}(u)du + \frac{\gamma}{N} \int_0^s Z_1^{N,-}(u) du. 
\ea\ees
Using Assumption \ref{main-d}, we have 
\bes
X^N(s) - \sqrt{N} x^N(s/N) \overset{d}{\approx} Z_1^N \overset{d}{\approx} Z_2^N,
\ees
where 
\bes\ba
Z_2^N(s) &  = \sqrt{N} Z^N(0)+ \int_0^{s} \sqrt{(\alpha^N)^3(\sigma^N)^2 + (\beta^N)^3(\varsigma^N)^2 + \theta^N N x^{N,+}(u/N) + \gamma^N N x^{N,-}(u/N) } dW(u) \\
& \quad - \theta^N \int_0^s Z_2^{N,+}(u)du + \gamma^N \int_0^s Z_2^{N,-}(u) du, \ s\ge 0.
\ea\ees
We next observe that 
\bes\ba
Nx^N(t/N) & = Nx^N(0) + (\alpha^N - \beta^N) t - N^2 \theta^N \int_0^{t/N} x^{N,+}(s)ds + N^2\gamma^N \int_0^{t/N} x^{N,-}(s)ds \\
& = Nx^N(0) + (\alpha^N - \beta^N) t - \theta^N \int_0^{t} N x^{N,+}(s/N)ds + \gamma^N \int_0^{t} N x^{N,-}(s/N)ds, \ t\ge 0.
\ea\ees
Define
\bes
x^N_2(t) = Nx^N(0) + (\alpha^N - \beta^N) t - \theta^N \int_0^t x^N_2(s)ds + \gamma^N \int_0^N x^N_2(s)ds, \ t\ge 0,
\ees 
Then $Nx^N(t/N) = x^N_2(t), t\ge 0,$ and we have 
\bes
X^N - x_2^N \overset{d}{\approx} Z_2^N,
\ees
and $Z_2$ can be rewritten as follows:
\bes\ba
Z_2^N(s) &  = X^N(0) - \sqrt{N}x^N(0) + \int_0^{s} \sqrt{(\alpha^N)^3(\sigma^N)^2 + (\beta^N)^3(\varsigma^N)^2 + \theta^N x^{N,+}_2(u) + \gamma^N x^{N,-}_2(u) } dW(u) \\
& \quad - \theta^N \int_0^s Z_2^{N,+}(u)du + \gamma^N \int_0^s Z_2^{N,-}(u) du, \ s\ge0.
\ea\ees
Thus for a double-ended queue with parameters satisfying Assumption \ref{main-d}, i.e. the reneging rates $\theta, \gamma$ are much small comparing with the arrival rates $\alpha, \beta$, the dynamics of the queue length process $\{X(t): t\geq0\}$ can be approximated by 
\be\label{approx-model-2}
X_2(t) = x_2(t) + Z_2(t), \ t\ge 0,
\ee
where 
\be\label{f-approx}
x_2(t) = (\alpha - \beta) t - \theta \int_0^t x^+_2(s)ds + \gamma \int_0^t x^-_2(s)ds, \ t\ge 0,
\ee
and 
\be\label{d-approx}
Z_2(t)= X(0) + \int_0^t \sqrt{{\alpha ^3}{\sigma ^2} + {\beta ^3}{\varsigma ^2} + \theta x^+_2(u) + \gamma x^-_2(u) } dW(u)  - \theta\int_0^t   Z^+_2(u) du +\gamma\int_0^t  Z^-_2(u)du, \ t\ge 0.
\ee
As our second diffusion approximate model, \eqref{approx-model-2} is used in Section \ref{sec:NE} to approximate the dynamics of the queue length process $\{X(t): t\geq0\}$ for a double-ended queue with parameters $\alpha, \beta, \sigma^2, \varsigma^2, \theta$, and $\gamma$.

\subsection{Tables}

\begin{table}[H]
\center \scalebox{0.8}{
\begin{tabular}
{|c|c|c c c c|}\hline \multicolumn{2}{|c|}{Exponential distribution} & \multicolumn{4}{|c|}{$L_1$} \\[2ex]
\hline ($\alpha$, $\beta$) & ($\theta$, $\gamma$) & $L_1^s$ & $L_1^p$ & $L_1^{d,1}$ & $L_1^{d,2}$\\[2ex]
\hline \hline

\multirow{6}{*}{(1, 1)}&(1, 1) & 0.0001 & 0 & 0 & 0\\

&{} & $\pm$0.0024 & NA & NA & NA\\

&(0.1, 0.1) & -0.0178 & 0 & 0 & 0\\

&{} & $\pm$0.0243 & NA & NA& NA\\

&(0.01, 0.01) & 0.1234 & 0 & 0 & 0\\

&{} & $\pm$0.2084 & NA & NA & NA\\ \hline

\multirow{6}{*}{(1, 1.5)} &(1, 1.5) & -0.2352 & -0.2343 & -0.2161 & -0.3333\\

& {}& $\pm$0.0022 & 0.41\% & 0.98\% & 41.7\%\\

&(0.1, 0.15) & -3.248 & -3.2532 & -3.2251 & -3.3333\\

& {}& $\pm$0.0192 & 0.16\% & 0.44\%& 2.63\%\\

&(0.01, 0.015) & -33.1485 & -33.3332 & -33.3327 & -33.3333\\

&{} & $\pm$0.1754 & 0.56\% & 0.56\% &0.56\% \\ \hline

\multirow{6}{*}{(1, 2)} &(1, 2) & -0.3876 & -0.3858 & -0.3178 & -0.5000\\

&{} & $\pm$0.002 & 0.47\% & 0.04\% &29\% \\

&(0.1, 0.2) & -4.9779 & -4.9719 & -4.9776 & -5.0000\\

&{} & $\pm$0.0157 & 0.12\% & 0.01\% & 0.45\%\\

&(0.01, 0.02) & -49.9609 & -50 & -50 & -50\\

&{} & $\pm$0.142 & 0.08\% & 0.08\% & 0.08\%\\ \hline
\end{tabular}}
\caption{The first moment of the stationary distribution when the arrival process is a Poisson process}
\end{table}

\begin{table}[H]
\center \scalebox{0.8}{
\begin{tabular}
{|c|c|c c c c|}\hline \multicolumn{2}{|c|}{Uniform distribution} & \multicolumn{4}{|c|}{$L_1$} \\[2ex]
\hline ($\alpha$, $\beta$) & ($\theta$, $\gamma$) & $L_1^s$ & $L_1^p$ & $L_1^{d,1}$ & $L_1^{d,2}$\\[2ex]
\hline \hline

\multirow{6}{*}{(1, 1)}&(1, 1) & 0.0004 & 0 & 0 & 0\\

&{} & $\pm$0.0017 & NA & NA & NA\\

&(0.1, 0.1) & -0.0009 & 0 & 0 & 0\\

&{} & $\pm$0.0141 & NA & NA & NA\\

&(0.01, 0.01) & -0.1309 & 0 & 0 & 0\\

&{} & $\pm$0.1231 & NA & NA & NA\\ \hline

\multirow{6}{*}{(1, 1.5)} &(1, 1.5) & -0.2736 & -0.2343 & -0.2979 & -0.3333\\

& {}& $\pm$0.0015 & 14.39\% & 8.87\% & 21.82\%\\

&(0.1, 0.15) & -3.3315 & -3.2532 & -3.3280 & -3.3333\\

& {}& $\pm$0.0114 & 2.35\% & 0.10\% & 0.054\%\\

&(0.01, 0.015) & -33.4634 & -33.3332 & -33.3333 & -33.3333\\

&{} & $\pm$0.1132 & 0.39\% & 0.39\% & 0.39\%\\ \hline

\multirow{6}{*}{(1, 2)} &(1, 2) & -0.4375 & -0.3858 & -0.4714 & -0.5000\\

&{} & $\pm$0.0013 & 11.82\% & 7.76\% & 14.28\%\\

&(0.1, 0.2) & -4.9946 & -4.9719 & -4.9998 & -5.0000\\

&{} & $\pm$0.0109 & 0.45\% & 0.10\% & 0.11\%\\

&(0.01, 0.02) & -50.0716 & -50 & -50 & -50\\

&{} & $\pm$0.1036 & 0.14\% & 0.14\% & 0.14\%\\ \hline
\end{tabular}}
\caption{The first moment of the stationary distribution when the inter-arrival times follow  Uniform
distribution}
\end{table}

\begin{table}[H]
\center \scalebox{0.8}{
\begin{tabular}
{|c|c|c c c c|}\hline \multicolumn{2}{|c|}{Erlang distribution} & \multicolumn{4}{|c|}{$L_1$} \\[2ex]
\hline ($\alpha$, $\beta$) & ($\theta$, $\gamma$) & $L_1^s$ & $L_1^p$ & $L_1^{d,1}$ & $L_1^{d,2}$\\[2ex]
\hline \hline

\multirow{6}{*}{(1, 1)}&(1, 1) & 0.0117 & 0 & 0 & 0\\

&{} & $\pm$0.0024 & NA & NA & NA\\

&(0.1, 0.1) & 0.0505 & 0 & 0 & 0\\

&{} & $\pm$0.0186 & NA & NA & NA\\

&(0.01, 0.01) & 0.0807 & 0 & 0 & 0\\

&{} & $\pm$0.1848 & NA & NA & NA\\ \hline

\multirow{6}{*}{(1, 1.5)} &(1, 1.5) & -0.2654 & -0.2343 & -0.2804 & -0.3333\\

& {}& $\pm$0.002 & 10.74\% & 6.84\% & 25.58\%\\

&(0.1, 0.15) & -3.2975 & -3.2532 & -3.3165 & -3.3333\\

& {}& $\pm$0.0155 & 1.34\% & 0.57\% & 1.08\%\\

&(0.01, 0.015) & -33.1629 & -33.3332 & -33.3333 & -33.3333\\

&{} & $\pm$0.1613 & 0.51\% & 0.51\% & 0.51\%\\ \hline

\multirow{6}{*}{(1, 2)} &(1, 2) & -0.4285 & -0.3858 & -0.4493 & -0.5000\\

&{} & $\pm$0.0018 & 9.96\% & 4.87\% & 16.69\%\\

&(0.1, 0.2) & -4.9832 & -4.9719 & -4.9983 & -5.0000\\

&{} & $\pm$0.015 & 0.23\% & 0.30\% & 0.34\%\\

&(0.01, 0.02) & -50.089 & -50 & -50 & -50\\

&{} & $\pm$0.1507 & 0.18\% & 0.18\% & 0.18\%\\ \hline
\end{tabular}}
\caption{The first moment of the stationary distribution when the inter-arrival times follow Erlang
distribution}
\end{table}

\begin{table}[H]
\center \scalebox{0.8}{
\begin{tabular}
{|c|c|c c c c|}\hline \multicolumn{2}{|c|}{Hyper-exponential distribution} & \multicolumn{4}{|c|}{$L_1$} \\[2ex]
\hline ($\alpha$, $\beta$) & ($\theta$, $\gamma$) & $L_1^s$ & $L_1^p$ & $L_1^{d,1}$ & $L_1^{d,2}$\\[2ex]
\hline \hline

\multirow{6}{*}{(1, 1)}&(1, 1) & 0.0022 & 0 & 0 &0\\

&{} & $\pm$0.0032 & NA & NA&NA\\

&(0.1, 0.1) & -0.0169 & 0 & 0&0\\

&{} & $\pm$0.0321 & NA & NA&NA\\

&(0.01, 0.01) & 0.016 & 0 & 0&0\\

&{} & $\pm$0.3177 & NA & NA&NA\\ \hline

\multirow{6}{*}{(1, 1.5)} &(1, 1.5) & -0.2039 & -0.2343 & -0.1735&-0.3333\\

& {}& $\pm$0.0028 & 14.89\% & 14.92\% & 63.46\%\\

&(0.1, 0.15) & -3.1406 & -3.2532 & -3.1368 &-3.3333\\

& {}& $\pm$0.0271 & 3.59\% & 0.12\% & 6.13\%\\

&(0.01, 0.015) & -33.2392 & -33.3332 & -33.3261& 33.3333\\

&{} & $\pm$0.237 & 0.28\% & 0.26\% & 0.28\%\\ \hline

\multirow{6}{*}{(1, 2)} &(1, 2) & -0.3383 & -0.3858 & -0.2866& -0.5\\

&{} & $\pm$0.0026 & 14.04\% & 15.26\% & 47.80\%\\

&(0.1, 0.2) & -4.8819 & -4.9719 & -4.8822 & -5\\

&{} & $\pm$0.0214 & 1.84\% & 0.01\% & 2.42\%\\

&(0.01, 0.02) & -50.1134 & -50 & -50 & -50\\

&{} & $\pm$0.1959 & 0.23\% & 0.23\% & 0.23\%\\ \hline
\end{tabular}}
\caption{The first moment of the stationary distribution when the inter-arrival times follow  hyper-exponential
distribution}
\end{table}

\begin{table}[H]
\center \scalebox{0.8}{
\begin{tabular}
{|c|c|c c c c|}\hline \multicolumn{2}{|c|}{Exponential distribution} & \multicolumn{4}{|c|}{$L_2$} \\[2ex]
\hline ($\alpha$, $\beta$) & ($\theta$, $\gamma$) & $L_2^s$ & $L_2^p$ & $L_2^{d,1}$ & $L_2^{d,2}$\\[2ex]
\hline \hline

\multirow{6}{*}{(1, 1)}&(1, 1) & 1.409 & 1.4104 & 1 & 1\\

&{} & $\pm$0.0042 & 0.10\% & 29.03\% & 29.03\%\\

&(0.1, 0.1) & 11.3894 & 11.3045 & 10 & 10\\

&{} & $\pm$0.0838 & 0.74\% & 12.20\% & 12.20\%\\

&(0.01, 0.01) & 103.2893 & 104.0397 & 100 &100\\

&{} & $\pm$2.2995 & 0.73\% & 3.18\% &3.18\% \\ \hline

\multirow{6}{*}{(1, 1.5)} &(1, 1.5) & 1.4354 & 1.4372 & 1.3194 & 1.7052\\

& {}& $\pm$0.0038 & 0.12\% & 8.1\% & 18.8\%\\

&(0.1, 0.15) & 21.2369 & 21.2498 & 21.9505 & 27.0518\\

& {}& $\pm$0.1458 & 0.06\% & 3.36\% & 27.38\%\\

&(0.01, 0.015) & 1218.2624 & 1211.1069 & 1219.4 & 1290.5\\

&{} & $\pm$12.3607 & -0.59\% & 0.09\% & 5.9\%\\ \hline

\multirow{6}{*}{(1, 2)} &(1, 2) & 1.4828 & 1.4841 & 1.7014 & 2.6287\\

&{} & $\pm$0.0036 & 0.09\% & 14.74\% & 77.28\%\\

&(0.1, 0.2) & 34.8606 & 34.956 & 37.3703 & 48.7868\\

&{} & $\pm$0.1677 & 0.27\% & 7.2\% &40\%\\

&(0.01, 0.02) & 2601.2009 & 2600 & 2625 & 2737.9\\

&{} & $\pm$15.2948 & 0.05\% & 0.9\% & 5.25\%\\ \hline
\end{tabular}}
\caption{The second moment of the stationary distribution when arrival process is a Poisson process}
\end{table}

\begin{table}[H]
\center \scalebox{0.8}{
\begin{tabular}
{|c|c|c c c c|}\hline \multicolumn{2}{|c|}{Uniform distribution} & \multicolumn{4}{|c|}{$L_2$} \\[2ex]
\hline ($\alpha$, $\beta$) & ($\theta$, $\gamma$) & $L_2^s$ & $L_2^p$ & $L_2^{d,1}$ & $L_2^{d,2}$\\[2ex]
\hline \hline

\multirow{6}{*}{(1, 1)}&(1, 1) & 0.8254 & 1.4104 & 0.3333 & 0.3333\\

&{} & $\pm$0.002 & 70.87\% & 59.62\% & 59.62\%\\

&(0.1, 0.1) & 4.3492 & 11.3045 & 3.3333& 3.3333\\

&{} & $\pm$0.0336 & 159.92\% & 23.36\% & 23.36\%\\

&(0.01, 0.01) & 34.6831 & 104.0397 & 33.3333 & 33.3333\\

&{} & $\pm$0.7472 & 199.97\% & 3.89\% & 3.89\%\\ \hline

\multirow{6}{*}{(1, 1.5)} &(1, 1.5) & 0.8961 & 1.4372 & 0.3993 & 0.6779\\

& {}& $\pm$0.002 & 60.38\% & 55.45\% & 24.35\%\\

&(0.1, 0.15) & 15.775 & 21.2498 & 13.8778 & 16.7789\\

& {}& $\pm$0.0952 & 34.71\% & 12.03\% & 6.34\%\\

&(0.01, 0.015) & 1148.5144 & 1211.1069 & 1138.8889 & 1167.8\\

&{} & $\pm$7.7166 & 5.45\% & 0.84\% & 1.68\%\\ \hline

\multirow{6}{*}{(1, 2)} &(1, 2) & 1.0102 & 1.4841 & 0.5019 & 1.0429\\

&{} & $\pm$0.0021 & 46.91\% & 50.32\% & 3.24\%\\

&(0.1, 0.2) & 30.1933 & 34.956 & 27.4992 & 32.9289\\

&{} & $\pm$0.1226 & 15.77\% & 8.92\% & 9.06\%\\

&(0.01, 0.02) & 2551.7944 & 2600 & 2525 & 2579.3\\

&{} & $\pm$10.7995 & 1.89\% & 1.05\% & 1.07\%\\ \hline
\end{tabular}}
\caption{The second moment of the stationary distribution when the inter-arrival times follow  Uniform
distribution}
\end{table}

\begin{table}[H]
\center \scalebox{0.8}{
\begin{tabular}
{|c|c|c c c c|}\hline \multicolumn{2}{|c|}{Erlang distribution} & \multicolumn{4}{|c|}{$L_2$} \\[2ex]
\hline ($\alpha$, $\beta$) & ($\theta$, $\gamma$) & $L_2^s$ & $L_2^p$ & $L_2^{d,1}$ & $L_2^{d,2}$\\[2ex]
\hline \hline

\multirow{6}{*}{(1, 1)}&(1, 1) & 0.9304 & 1.4104 & 0.5000 & 0.5000\\

&{} & $\pm$0.0064 & 51.95\% & 46.26\% & 46.26\%\\

&(0.1, 0.1) & 6.0528 & 11.3045 & 5.0000 & 5.0000\\

&{} & $\pm$0.1485 & 86.76\% & 17.39\% & 17.39\%\\

&(0.01, 0.01) & 48.0992 & 104.0397 & 50.0000 & 50.0000 \\

&{} & $\pm$4.4079 & 116.3\% & 3.95\% & 3.95\%\\ \hline

\multirow{6}{*}{(1, 1.5)} &(1, 1.5) & 0.9857 & 1.4372 & 0.5526 & 0.8550\\

& {}& $\pm$0.0056 & 45.80\% & 43.94\% & 13.26\%\\

&(0.1, 0.15) & 16.5479 & 21.2498 & 15.2503 & 18.5501\\

& {}& $\pm$0.1894 & 28.41\% & 7.84\% & 12.1\%\\

&(0.01, 0.015) & 1158.5 & 1211.1069 & 1152.8 & 1185.5\\

&{} & $\pm$16.5071 & 4.54\% & 0.05\% & 2.33\%\\ \hline

\multirow{6}{*}{(1, 2)} &(1, 2) & 1.0728 & 1.4841 & 0.6375 & 1.2411\\

&{} & $\pm$0.0052 & 38.34\% & 40.2\% & 15.69\\

&(0.1, 0.2) & 31.4929 & 34.956 & 28.7436 & 34.9112\\

&{} & $\pm$0.2351 & 11\% & 8.73\% & 10.85\%\\

&(0.01, 0.02) & 2542.1 & 2600 & 2567.3 & 2599.1\\

&{} & $\pm$20.576 & 2.28\% & 1\% & 2.24\%\\ \hline
\end{tabular}}
\caption{The second moment of the stationary distribution when the inter-arrival times follow  Erlang
distribution}
\end{table}

\begin{table}[H]
\center \scalebox{0.8}{
\begin{tabular}
{|c|c|c c c c|}\hline \multicolumn{2}{|c|}{Hyper-exponential distribution} & \multicolumn{4}{|c|}{$L_2$} \\[2ex]
\hline ($\alpha$, $\beta$) & ($\theta$, $\gamma$) & $L_2^s$ & $L_2^p$ & $L_2^{d,1}$ & $L_2^{d,2}$\\[2ex]
\hline \hline

\multirow{6}{*}{(1, 1)}&(1, 1) & 1.9943 & 1.4104 & 2 & 2\\

&{} & $\pm$0.0063 & 29.28\% & 0.29\% &0.29\% \\

&(0.1, 0.1) & 20.8656 & 11.3045 & 20 & 20\\

&{} & $\pm$0.1625 & 45.82\% & 4.15\% &4.15\%\\

&(0.01, 0.01) & 205.774 & 104.0397 & 200 & 200\\

&{} & $\pm$4.7111 & 49.44\% & 2.81\% &2.81\%\\ \hline

\multirow{6}{*}{(1, 1.5)} &(1, 1.5) & 1.9962 & 1.4372 & 2.0092 & 2.4491\\

& {}& $\pm$0.0057 & 28.01\% & 0.65\% & 22.69\% \\

&(0.1, 0.15) & 29.4329 & 21.2498 & 28.0112 & 34.4908\\

& {}& $\pm$0.1921 & 27.80\% & 4.83\% & 17.18\%\\

&(0.01, 0.015) & 1307.8225 & 1211.1069 & 1277.5943 & 1344.9\\

&{} & $\pm$16.5066 & 7.40\% & 2.31\% & 2.84\%\\ \hline

\multirow{6}{*}{(1, 2)} &(1, 2) & 2.0048 & 1.4841 & 2.0252 & 3.0251\\

&{} & $\pm$0.0048 & 25.97\% & 1.02\% & 50.89\%\\

&(0.1, 0.2) & 41.5526 & 34.956 & 39.8704 & 52.7513\\

&{} & $\pm$0.2282 & 15.88\% & 4.05\% & 26.95\%\\

&(0.01, 0.02) & 2660.7665 & 2600 & 2649.9986 & 2777.5\\

&{} & $\pm$20.4434 & 2.28\% & 0.40\% & 4.39\%\\ \hline
\end{tabular}}
\caption{The second moment of the stationary distribution when the inter-arrival times follow 
hyper-exponential distribution}
\end{table}

\subsection{Figures}

\begin{figure}[H]\center
\begin{tabular}{ccc}
\includegraphics[width=5cm]{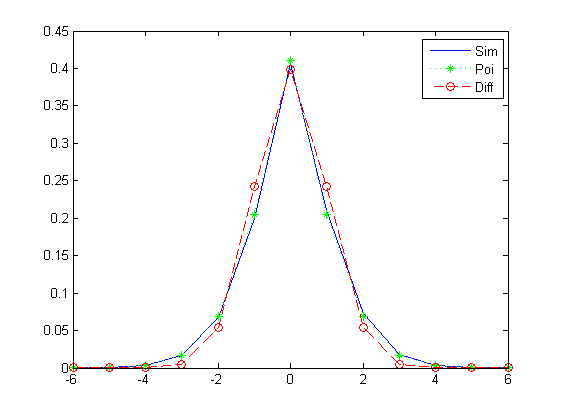} & \includegraphics[width=5cm]{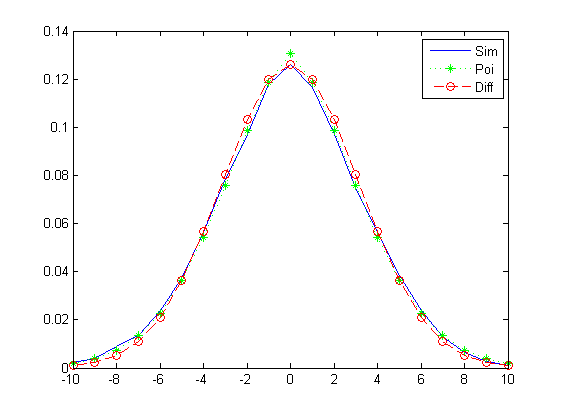} &
\includegraphics[width=5cm]{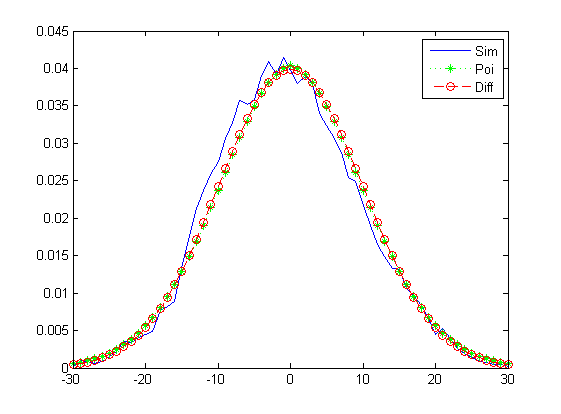}\\
  (a) $(\alpha, \beta, \theta, \gamma) = (1,1,1,1)$ & (b) $(\alpha, \beta, \theta, \gamma) = (1,1,0.1,0.1)$ &
  (c) $(\alpha, \beta, \theta, \gamma) = (1,1,0.01,0.01)$ \\[6pt]
  \includegraphics[width=5cm]{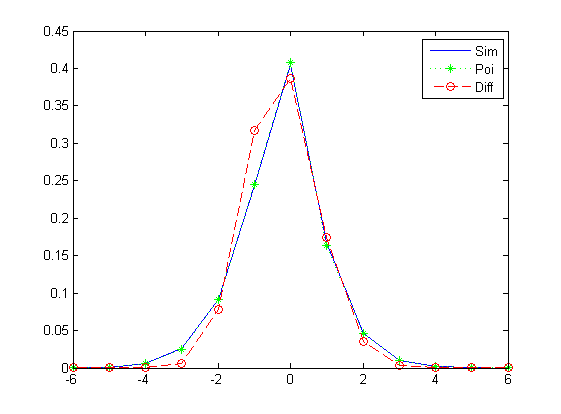} & \includegraphics[width=5cm]{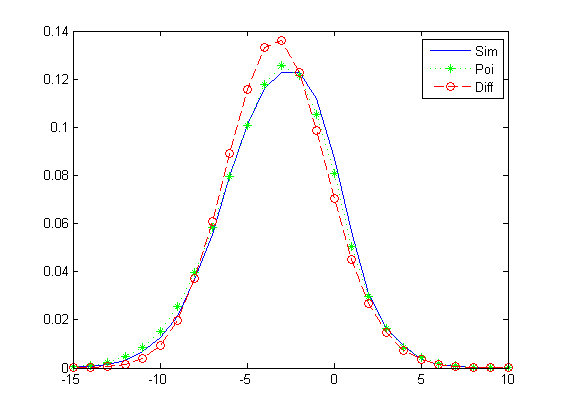} &
\includegraphics[width=5cm]{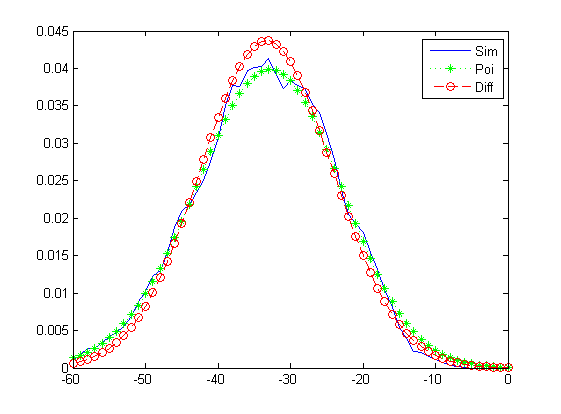}\\
  (a) $(\alpha, \beta, \theta, \gamma) = (1,1.5,1,1.5)$ & (b) $(\alpha, \beta, \theta, \gamma) = (1,1.5,0.1,0.15)$ &
  (c) $(\alpha, \beta, \theta, \gamma) = (1,1.5,0.01,0.015)$ \\[6pt]
  \includegraphics[width=5cm]{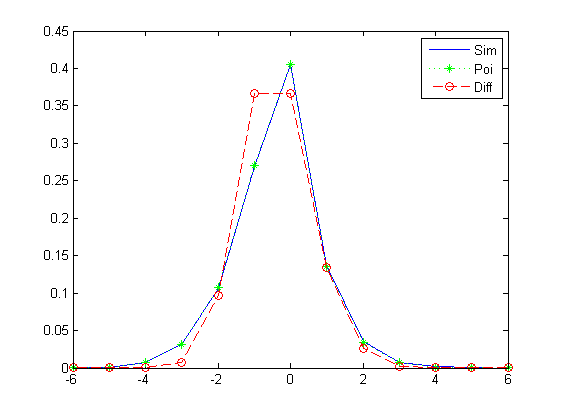} & \includegraphics[width=5cm]{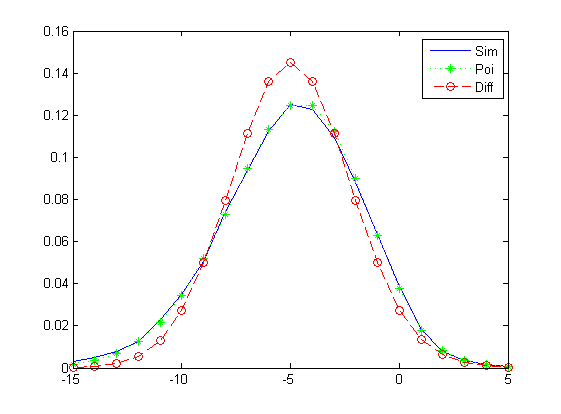} &
\includegraphics[width=5cm]{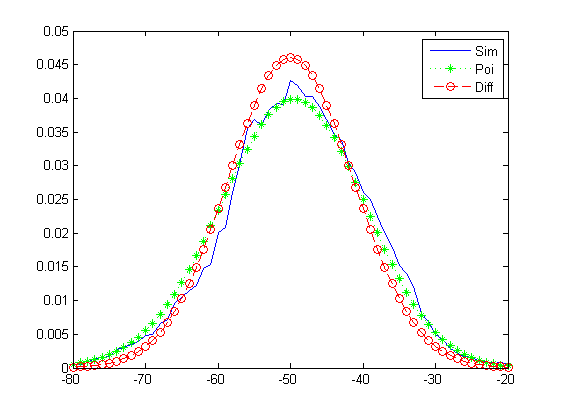}\\
  (a) $(\alpha, \beta, \theta, \gamma) = (1,2,1,2)$ & (b) $(\alpha, \beta, \theta, \gamma) = (1,2,0.1,0.2)$ &
  (c) $(\alpha, \beta, \theta, \gamma) = (1,2,0.01,0.02)$ \\[6pt]
\end{tabular}
\caption{Density functions by simulation method, Poisson approximation, and heavy traffic diffusion approximation, when
inter-arrival times follow exponential distribution} \label{fig:expplot}
\end{figure}

\begin{figure}[H]\center
\begin{tabular}{ccc}
\includegraphics[width=5cm]{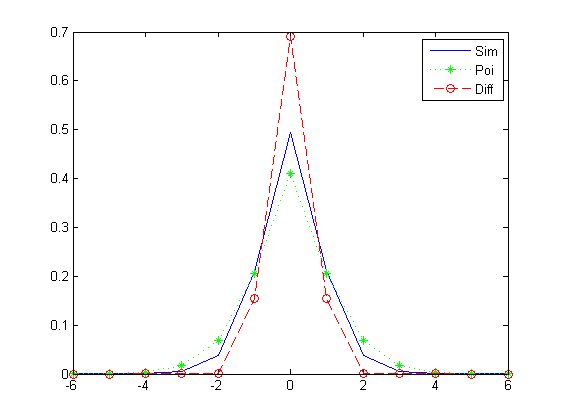} & \includegraphics[width=5cm]{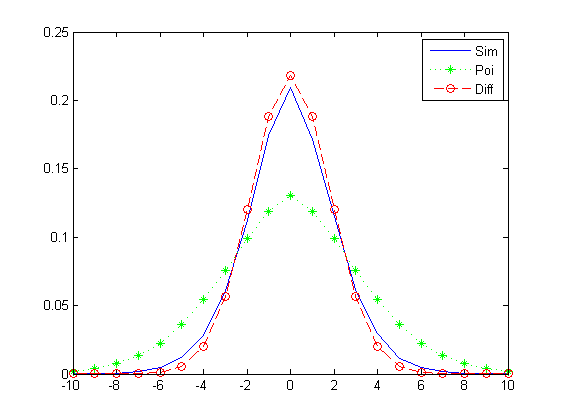} &
\includegraphics[width=5cm]{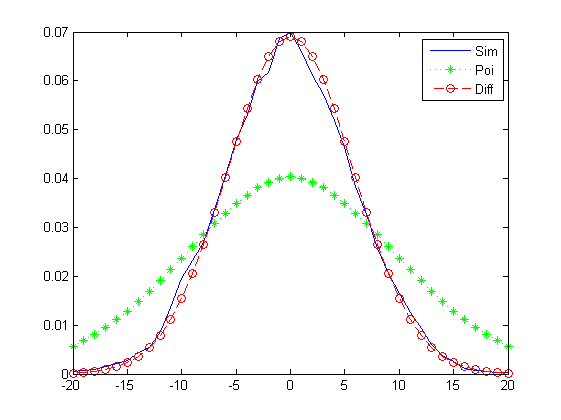}\\
  (a) $(\alpha, \beta, \theta, \gamma) = (1,1,1,1)$ & (b) $(\alpha, \beta, \theta, \gamma) = (1,1,0.1,0.1)$ &
  (c) $(\alpha, \beta, \theta, \gamma) = (1,1,0.01,0.01)$ \\[6pt]
  \includegraphics[width=5cm]{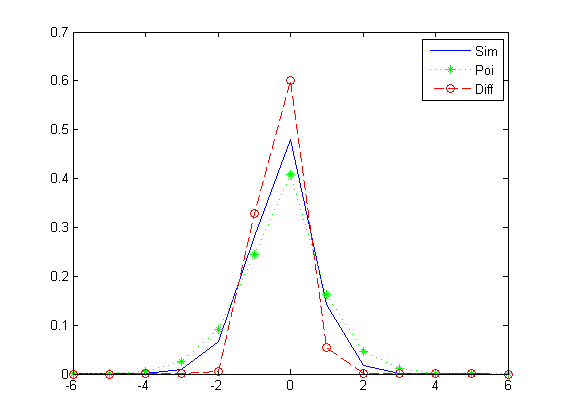} & \includegraphics[width=5cm]{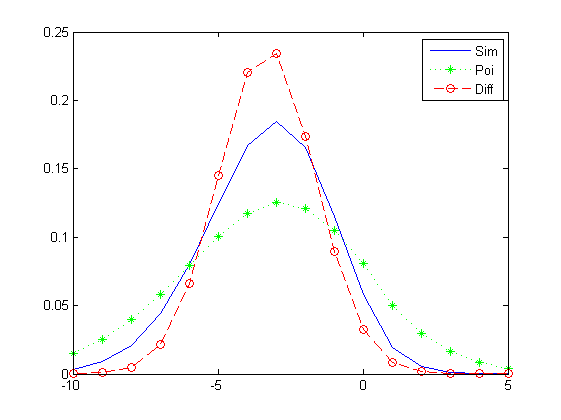} &
\includegraphics[width=5cm]{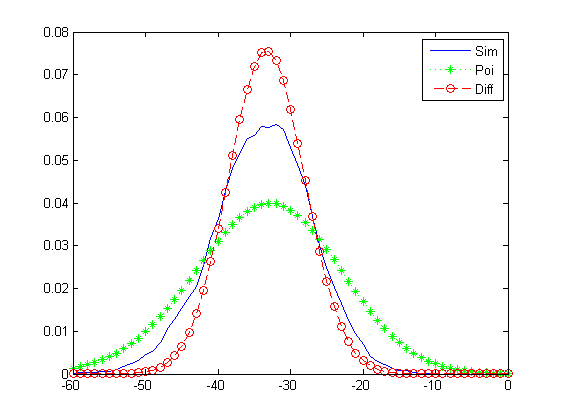}\\
  (a) $(\alpha, \beta, \theta, \gamma) = (1,1.5,1,1.5)$ & (b) $(\alpha, \beta, \theta, \gamma) = (1,1.5,0.1,0.15)$ &
  (c) $(\alpha, \beta, \theta, \gamma) = (1,1.5,0.01,0.015)$ \\[6pt]
  \includegraphics[width=5cm]{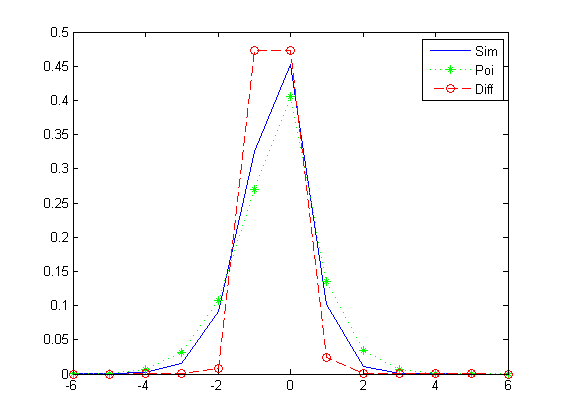} & \includegraphics[width=5cm]{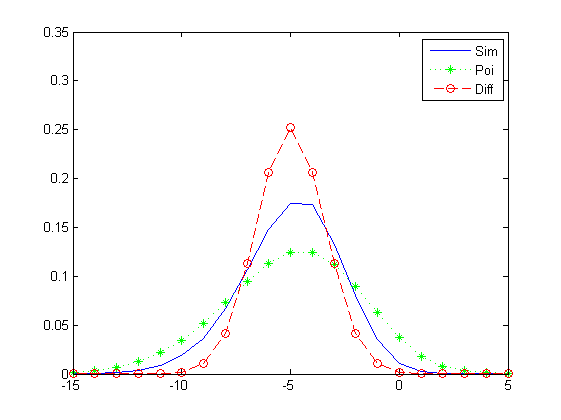} &
\includegraphics[width=5cm]{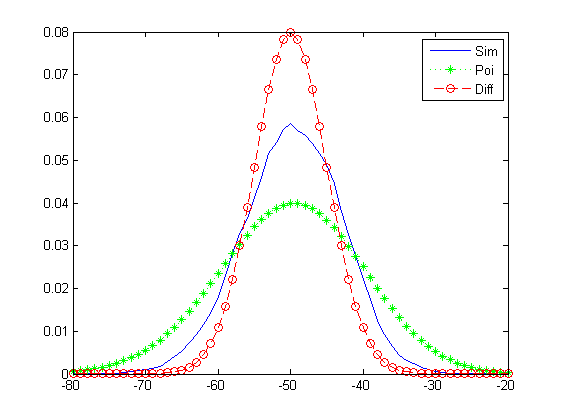}\\
  (a) $(\alpha, \beta, \theta, \gamma) = (1,2,1,2)$ & (b) $(\alpha, \beta, \theta, \gamma) = (1,2,0.1,0.2)$ &
  (c) $(\alpha, \beta, \theta, \gamma) = (1,2,0.01,0.02)$ \\[6pt]
\end{tabular}
\caption{Density functions by simulation method, Poisson approximation, and heavy traffic diffusion approximation, when
inter-arrival times follow uniform distribution} \label{fig:unifplot}
\end{figure}

\begin{figure}[H]\center
\begin{tabular}{ccc}
\includegraphics[width=5cm]{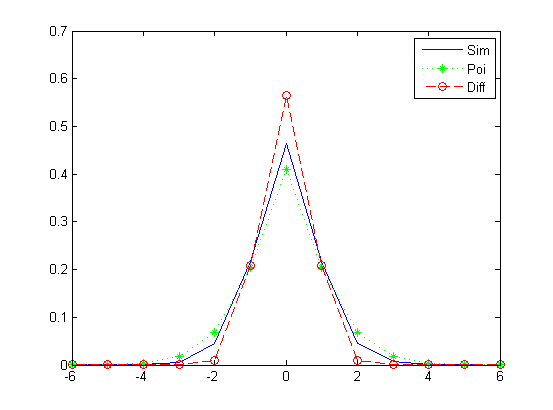} & \includegraphics[width=5cm]{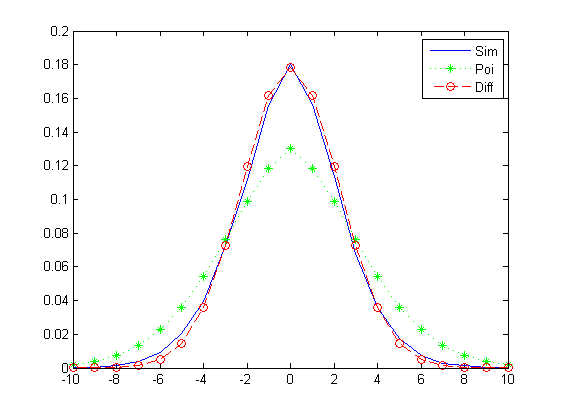} &
\includegraphics[width=5cm]{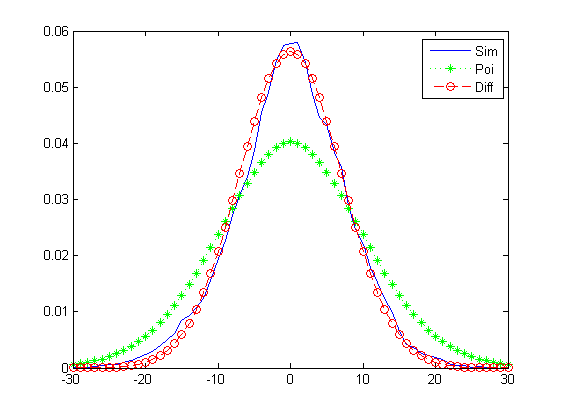}\\
  (a) $(\alpha, \beta, \theta, \gamma) = (1,1,1,1)$ & (b) $(\alpha, \beta, \theta, \gamma) = (1,1,0.1,0.1)$ &
  (c) $(\alpha, \beta, \theta, \gamma) = (1,1,0.01,0.01)$ \\[6pt]
  \includegraphics[width=5cm]{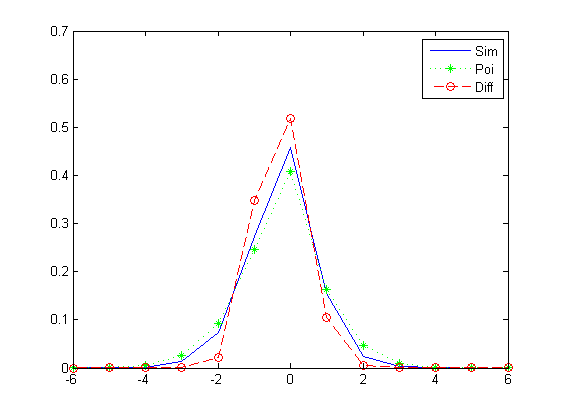} & \includegraphics[width=5cm]{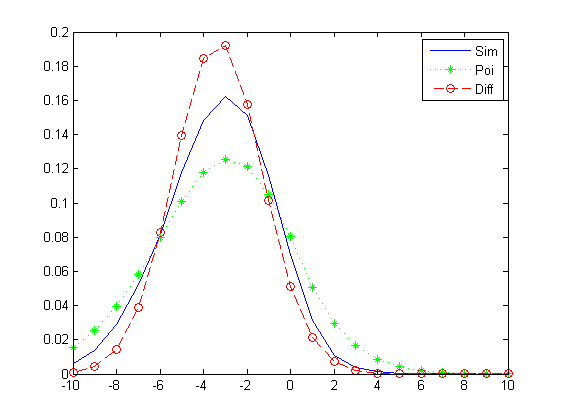} &
\includegraphics[width=5cm]{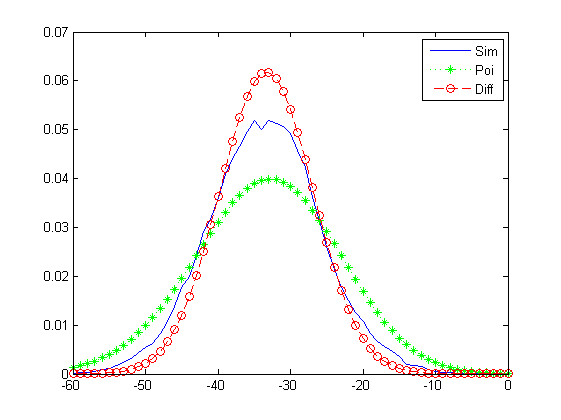}\\
  (a) $(\alpha, \beta, \theta, \gamma) = (1,1.5,1,1.5)$ & (b) $(\alpha, \beta, \theta, \gamma) = (1,1.5,0.1,0.15)$ &
  (c) $(\alpha, \beta, \theta, \gamma) = (1,1.5,0.01,0.015)$ \\[6pt]
  \includegraphics[width=5cm]{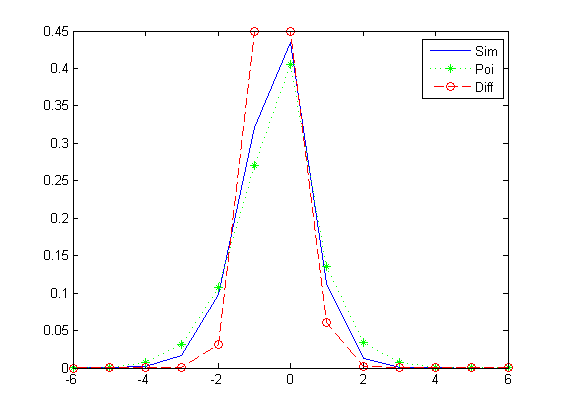} & \includegraphics[width=5cm]{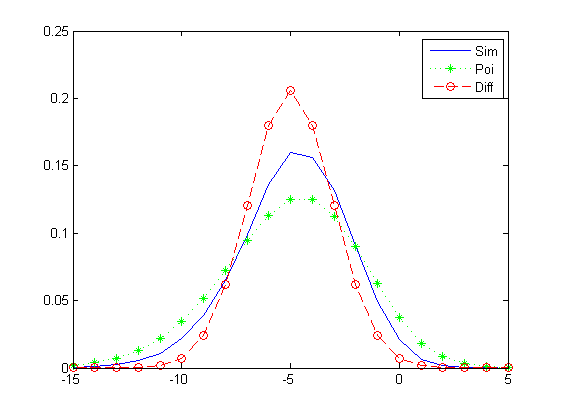} &
\includegraphics[width=5cm]{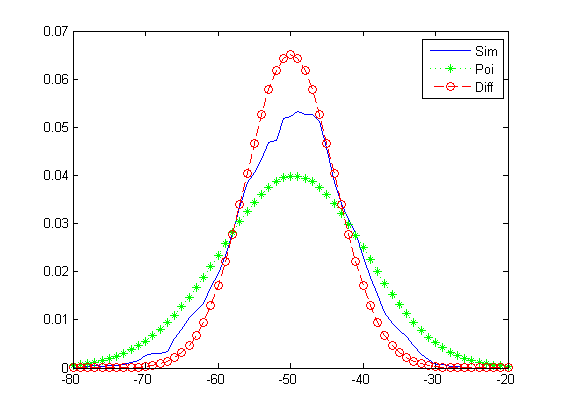}\\
  (a) $(\alpha, \beta, \theta, \gamma) = (1,2,1,2)$ & (b) $(\alpha, \beta, \theta, \gamma) = (1,2,0.1,0.2)$ &
  (c) $(\alpha, \beta, \theta, \gamma) = (1,2,0.01,0.02)$ \\[6pt]
\end{tabular}
\caption{Density functions by simulation method, Poisson approximation, and heavy traffic diffusion approximation, when
inter-arrival times follow  Erlang distribution} \label{fig:erlangplot}
\end{figure}

\begin{figure}[H]\center
\begin{tabular}{ccc}
\includegraphics[width=5cm]{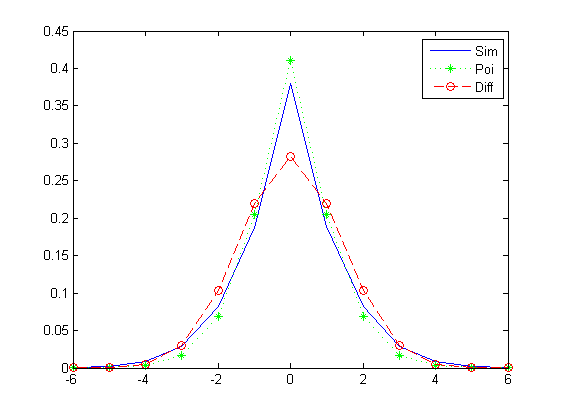} & \includegraphics[width=5cm]{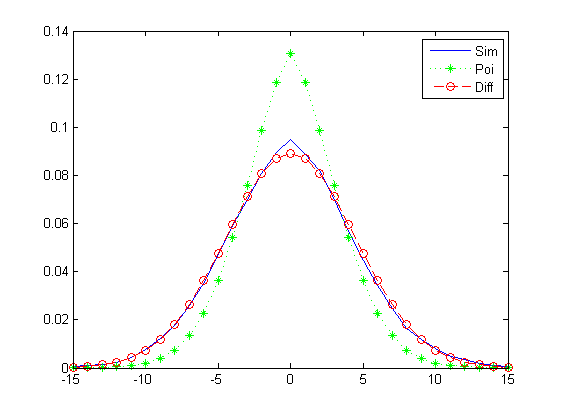} &
\includegraphics[width=5cm]{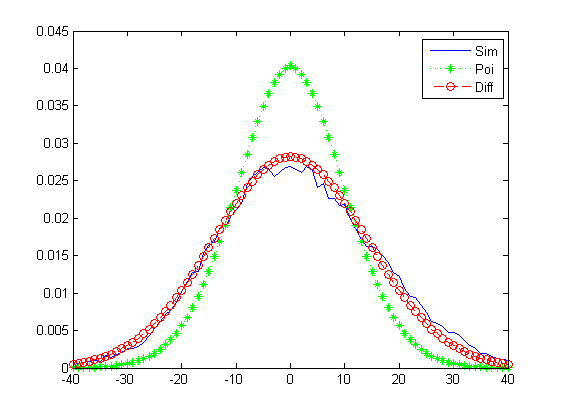}\\
  (a) $(\alpha, \beta, \theta, \gamma) = (1,1,1,1)$ & (b) $(\alpha, \beta, \theta, \gamma) = (1,1,0.1,0.1)$ &
  (c) $(\alpha, \beta, \theta, \gamma) = (1,1,0.01,0.01)$ \\[6pt]
  \includegraphics[width=5cm]{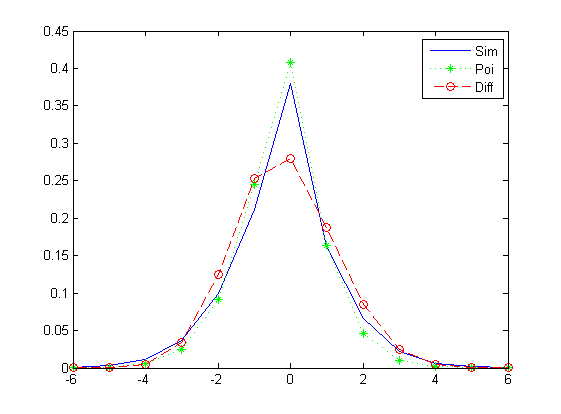} & \includegraphics[width=5cm]{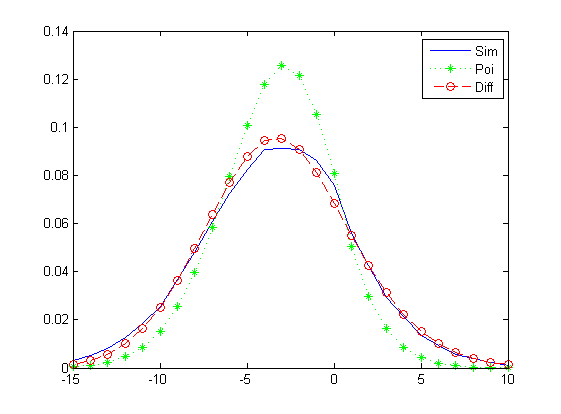} &
\includegraphics[width=5cm]{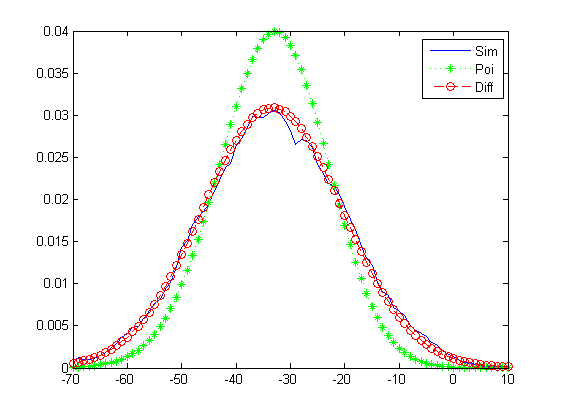}\\
  (a) $(\alpha, \beta, \theta, \gamma) = (1,1.5,1,1.5)$ & (b) $(\alpha, \beta, \theta, \gamma) = (1,1.5,0.1,0.15)$ &
  (c) $(\alpha, \beta, \theta, \gamma) = (1,1.5,0.01,0.015)$ \\[6pt]
  \includegraphics[width=5cm]{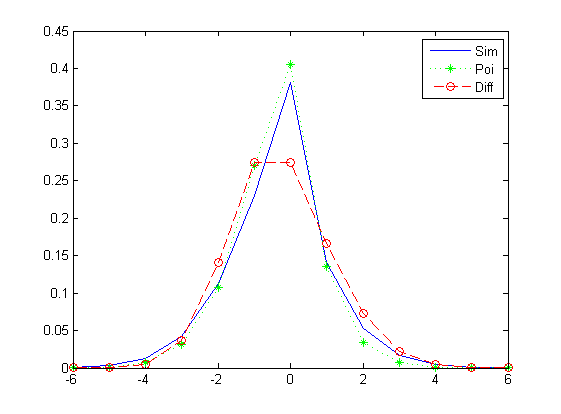} & \includegraphics[width=5cm]{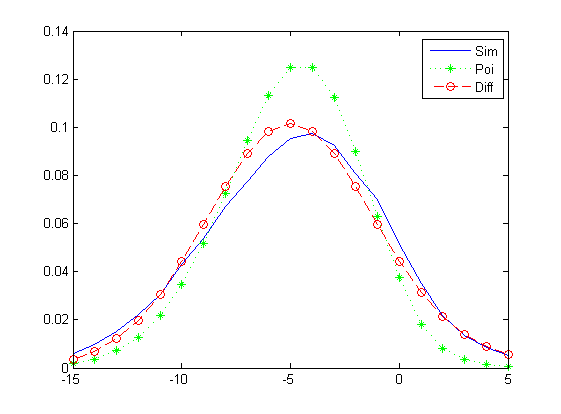} &
\includegraphics[width=5cm]{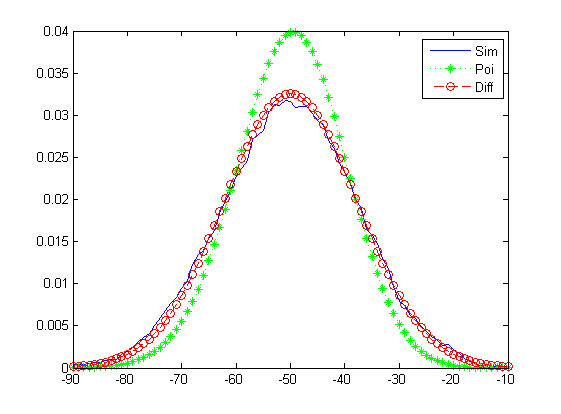}\\
  (a) $(\alpha, \beta, \theta, \gamma) = (1,2,1,2)$ & (b) $(\alpha, \beta, \theta, \gamma) = (1,2,0.1,0.2)$ &
  (c) $(\alpha, \beta, \theta, \gamma) = (1,2,0.01,0.02)$ \\[6pt]
\end{tabular}
\caption{Density functions by simulation method, Poisson approximation, and heavy traffic diffusion approximation, when
inter-arrival times follow hyper-exponential distribution} \label{fig:hypexpplot}
\end{figure}

%%%%%%%%%%%%%%%%%
\end{document}